\documentclass[12pt,twoside]{amsart}
\usepackage{}

\usepackage{amsthm}
\usepackage{amsmath}
\usepackage{amssymb}
\usepackage{amsfonts}
\usepackage{latexsym}
\usepackage{mathabx}
\usepackage{enumitem}
\usepackage{mathrsfs}\usepackage{stmaryrd}
\usepackage[all]{xy} \SelectTips{eu}{}
\usepackage{hyperref}
\usepackage{eucal}
\usepackage{xcolor}
\usepackage{tikz,tikz-cd}
\usetikzlibrary{arrows}
\usetikzlibrary{matrix}

\newtheorem*{intsetup*}{Setup}
\newtheorem*{intNotations*}{Notations and Conventions}
\newtheorem{intthm}{Theorem}[]

\newtheorem{intque}{Question}[]

\newtheorem*{intexa*}{Example}

\newcommand{\numberseries}{\bfseries}   
\newlength{\thmtopspace}                
\newlength{\thmbotspace}                
\newlength{\thmheadspace}               
\newlength{\thmindent}                  
\newtheoremstyle{bfupright head,slanted body}
                {\thmtopspace}{\thmbotspace}
                {\slshape}{\thmindent}{\bfseries}{.}{\thmheadspace}
                {{\numberseries \thmnumber{#2\;}}\thmnote{#3}}

\newtheoremstyle{bfupright head,upright body}
                {\thmtopspace}{\thmbotspace}
                {\upshape}{\thmindent}{\bfseries}{.}{\thmheadspace}
                {{\numberseries \thmnumber{#2\;}}\thmnote{#3}}

\newtheoremstyle{fixed bf head,slanted body}
                {\thmtopspace}{\thmbotspace}{\slshape}
                {\thmindent}{\bfseries}{.}{\thmheadspace}
                {{\numberseries \thmnumber{#2\;}}\thmname{#1}\thmnote{ (#3)}}

\newtheoremstyle{fixed bf head,upright body}
                {\thmtopspace}{\thmbotspace}{\upshape}
                {\thmindent}{\bfseries}{.}{\thmheadspace}
                {{\numberseries \thmnumber{#2\;}}\thmname{#1}\thmnote{ (#3)}}

\newtheoremstyle{numbered paragraph}
                {\thmtopspace}{\thmbotspace}{\upshape}
                {\thmindent}{\upshape}{}{\thmheadspace}
                {{\numberseries \thmnumber{#2.}}}

\theoremstyle{bfupright head,slanted body}
\newtheorem{res}{}[section]             \newtheorem*{res*}{}

\theoremstyle{bfupright head,upright body}
\newtheorem{bfhpg}[res]{}               \newtheorem*{bfhpg*}{}

\theoremstyle{slanted body, fixed bf head}
\newtheorem{thm}[res]{Theorem}          \newtheorem*{thm*}{Theorem}
\newtheorem{prp}[res]{Proposition}      \newtheorem*{prp*}{Proposition}
\newtheorem{cor}[res]{Corollary}        \newtheorem*{cor*}{Corollary}
\newtheorem{lem}[res]{Lemma}            \newtheorem*{lem*}{Lemma}
         \newtheorem*{que*}{Question}
\theoremstyle{upright body ,fixed bf head}
\theoremstyle{definition}
\newtheorem{setup}[res]{Setup}           \newtheorem*{setup*}{Setup}
\newtheorem{dfn}[res]{Definition}       \newtheorem*{dfn*}{Definition}
\newtheorem{rmk}[res]{Remark}           \newtheorem*{rmk*}{Remark}
\newtheorem{exa}[res]{Example}           \newtheorem*{exa*}{Example}
           \newtheorem*{exer*}{Exercise}

\theoremstyle{numbered paragraph}
\newtheorem{ipg}[res]{}

\newlength{\thmlistleft}        
\newlength{\thmlistright}       
\newlength{\thmlistpartopsep}   
\newlength{\thmlisttopsep}      
\newlength{\thmlistparsep}      
\newlength{\thmlistitemsep}     

\newcounter{eqc}
\newenvironment{eqc}{\begin{list}{\upshape (\textit{\roman{eqc}})}%
    {\usecounter{eqc}%
      \setlength{\leftmargin}{\thmlistleft}%
      \setlength{\labelwidth}{\thmlistleft}%
      \setlength{\rightmargin}{\thmlistright}%
      \setlength{\partopsep}{\thmlistpartopsep}%
      \setlength{\topsep}{\thmlisttopsep}%
      \setlength{\parsep}{\thmlistparsep}%
      \setlength{\itemsep}{\thmlistitemsep}}}%
  {\end{list}}%

\newcounter{prt}
\newenvironment{prt}{\begin{list}{\upshape (\alph{prt})}%
    {\usecounter{prt}%
      \setlength{\leftmargin}{\thmlistleft}%
      \setlength{\labelwidth}{\thmlistleft}%
      \setlength{\rightmargin}{\thmlistright}%
      \setlength{\partopsep}{\thmlistpartopsep}%
      \setlength{\topsep}{\thmlisttopsep}%
      \setlength{\parsep}{\thmlistparsep}%
      \setlength{\itemsep}{\thmlistitemsep}}}%
  {\end{list}}%

\newcounter{rqm}
\newenvironment{rqm}{\begin{list}{\upshape (\arabic{rqm})}%
    {\usecounter{rqm}%
      \setlength{\leftmargin}{\thmlistleft}%
      \setlength{\labelwidth}{\thmlistleft}%
      \setlength{\rightmargin}{\thmlistright}%
      \setlength{\partopsep}{\thmlistpartopsep}%
      \setlength{\topsep}{\thmlisttopsep}%
      \setlength{\parsep}{\thmlistparsep}%
      \setlength{\itemsep}{\thmlistitemsep}}}%
  {\end{list}}%

  {\end{list}}%

\newenvironment{prf*}[1][Proof]{%
  \begin{proof}[\bf #1]
    \setcounter{equation}{0}
    }
  {\end{proof}
}

\newcommand{\pgref}[1]{\ref{#1}}

\renewcommand{\eqref}[1]{(\pgref{eq:#1})}

\newcommand{\thmcite}[2][?]{\cite[Theorem~#1]{#2}}
\newcommand{\prpcite}[2][?]{\cite[Proposition~#1]{#2}}

\newcommand{\corcite}[2][?]{\cite[Corollary~#1]{#2}}
\newcommand{\lemcite}[2][?]{\cite[Lemma~#1]{#2}}

\newcommand{\dfncite}[2][?]{\cite[Definition~#1]{#2}}
\newcommand{\rmkcite}[2][?]{\cite[Remark~#1]{#2}}
\newcommand{\exacite}[2][?]{\cite[Example~#1]{#2}}
\newcommand{\exercite}[2][?]{\cite[Exercise~#1]{#2}}
\newcommand{\eqclbl}[1]{{\upshape(\textit{#1})}}
\newcommand{\proofofimp}[3][:]{\mbox{\eqclbl{#2}$\!\implies\!$\eqclbl{#3}#1}}
\numberwithin{equation}{res}

\def\urltilda{\kern -.15em\lower .7ex\hbox{\~{}}\kern .04em}

\newcommand{\coker}{\mbox{\rm coker}}
\newcommand{\calC}{\mathcal{C}}
\newcommand{\calW}{\mathcal{W}}
\newcommand{\calF}{\mathcal{F}}
\newcommand{\calA}{\mathcal{A}}
\newcommand{\calB}{\mathcal{B}}
\newcommand{\calD}{\mathcal{D}}
\newcommand{\calS}{\mathcal{S}}
\newcommand{\calM}{\mathcal{M}}

\newcommand{\calE}{\mathcal{E}}
\newcommand{\calX}{\mathcal{X}}
\newcommand{\calY}{\mathcal{Y}}
\newcommand{\calT}{\mathcal{T}}
\newcommand{\sfC}{\mathsf{C}}
\newcommand{\sfD}{\mathsf{D}}
\newcommand{\Mon}{\mathsf{Mon}}
\newcommand{\Epi}{\mathsf{Epi}}

\newcommand{\PGF}{\mathsf{PGF}}

\newcommand{\sfW}{\mathsf{W}}

\newcommand{\Proj}{\mathsf{Proj}}
\newcommand{\id}{\mathrm{id}}

\newcommand{\Mod}{\mathsf{Mod}}
\newcommand{\Ind}{\mathrm{Ind}}

\newcommand{\Ho}{\mathrm{Ho}}
\newcommand{\St}{\mathrm{St}}
\newcommand{\Id}{\mathrm{Id}}
\newcommand{\xra}[2][]{\xrightarrow[#1]{\:#2\:}}
\newcommand{\Rop}{R^{\sf op}}
\newcommand{\op}{\sf op}

\newcommand{\Bop}{B^{\sf op}}
\newcommand{\cok}{\mbox{\rm coker}}
\newcommand{\kernel}{\mbox{\rm ker}}
\newcommand{\Hom}{\operatorname{Hom}}
\newcommand{\Aut}{\operatorname{Aut}}
\newcommand{\Ext}{\operatorname{Ext}}
\newcommand{\Cofib}{\sf{Cofib}}
\newcommand{\Fib}{\sf{Fib}}
\newcommand{\Weq}{\sf{Weq}}

   \def\soft#1{\leavevmode\setbox0=\hbox{h}\dimen7=\ht0\advance
    \dimen7 by-1ex\relax\if t#1\relax\rlap{\raise.6\dimen7
    \hbox{\kern.3ex\char'47}}#1\relax\else\if T#1\relax
    \rlap{\raise.5\dimen7\hbox{\kern1.3ex\char'47}}#1\relax
    \else\if d#1\relax\rlap{\raise.5\dimen7\hbox{\kern.9ex
    \char'47}}#1\relax\else\if D#1\relax\rlap{\raise.5\dimen7
    \hbox{\kern1.4ex\char'47}}#1\relax\else\if l#1\relax
    \rlap{\raise.5\dimen7\hbox{\kern.4ex\char'47}}#1\relax
    \else\if L#1\relax\rlap{\raise.5\dimen7\hbox{\kern.7ex
    \char'47}}#1\relax\else\message{accent \string\soft
    \space #1 not defined!}#1\relax\fi\fi\fi\fi\fi\fi}

\begin{document}

\title[Transfer of abelian model structures to equivariant categories]%
{Transfer of abelian model structures to equivariant categories and homotopy squares}

\author[Z.X. Di]{Zhenxing Di}
\address{Zhenxing Di: School of Mathematical Sciences,
Huaqiao University, Quanzhou 362021, China}
\email{dizhenxing@163.com}

\author[L.P. Li]{Liping Li}
\address{Liping Li: Department of Mathematics, Hunan Normal University, Changsha 410081, China}
\email{lipingli@hunnu.edu.cn}

\author[L. Liang]{Li Liang}
\address{Li Liang: Department of Mathematics, Gansu Center for Fundamental Research in Complex Systems Analysis and Control, Lanzhou Jiaotong University, Lanzhou 730070, China}
\email{lliangnju@gmail.com}
\urladdr{https://sites.google.com/site/lliangnju}

\author[G.L. Tang]{Guoliang Tang}
\address{Guoliang Tang: School of Mathematical Sciences,
Zhejiang Normal University, Jinhua 321004, China}
\email{tangguoliang970125@163.com}

\author[R.M. Zhu]{Rongmin Zhu}
\address{Rongmin Zhu (Corresponding Author): School of Mathematical Sciences,
Huaqiao University, Quanzhou 362021, China}
\email{rongminzhu@hotmail.com}

\thanks{Z.X. Di was partly supported by NSF of China (Grant No. 12471034),
Scientific Research Fund of Fujian Province (Grant No. 605-52525002) and
Scientific Research Fund of Huaqiao University (Grant No. 605-50Y22050); L.P. Li was partly supported by NSF of China (Grant No. 12171146); L. Liang was partly supported by NSF of China (Grant No. 12271230), NSF of Gansu Province (Grant No. 26JRRA059) and the Foundation for Innovative Fundamental Research Group Project of Gansu Province (Grant No. 25JRRA805); R.M. Zhu was partly supported by NSF of China (Grant No. 12201223).}



\renewcommand{\thefootnote}{\alph{footnote}}

\begin{abstract}
Let $G$ be a finite group acting on a Grothendieck category $\calA$ with enough projectives, such that $|G|$ is invertible in $\calA$.
We prove a general lifting theorem for abelian model structures from $\calA$ to its equivariant category $\calA^G$, and establish a triangle equivalence up to retracts between the corresponding homotopy categories.
We also construct a commutative square whose horizontal functors are triangle equivalences and whose vertical comparison functors are triangle equivalences up to retracts.
This square relates derived functors on the lifted equivariant model categories to the equivariantizations of the derived functors on the original homotopy categories.
In the module category setting, we illustrate the above results using the PGF Hovey triples, and apply them to homotopy squares induced by a Frobenius bimodule and by a stable equivalence of adjoint type.
\end{abstract}

\maketitle
	\tableofcontents

\section*{Introduction}\label{Introduction}
\noindent
Abelian model structures provide a natural framework for developing homotopy theory in abelian categories.
By Hovey's correspondence \cite{HOVEY2002}, an abelian model structure on an abelian category $\calA$ can be represented by a Hovey triple $\calM=(\calC,\calW,\calF)$, determined by two compatible complete cotorsion pairs.
The homotopy category $\Ho(\calM)$ associated with any Hovey triple is naturally triangulated; see \cite[Section 6]{2025GillBook}.
When $\calM$ is hereditary, this triangulated category can moreover be described in terms of the stable category of the Frobenius category of bifibrant objects.
Many model structures arising in Gorenstein homological algebra are hereditary; see e.g. \cite{2016Gill,HOVEY2002,2020Sto}.

On the other hand, group actions on categories provide a categorical way to incorporate symmetry into homological constructions.
Equivariant categories occur naturally in representation theory, algebraic geometry, and homological algebra; see e.g. \cite{JV94,2015Chen,2017Chen,2010DGNO,2014Elag,2019Sun}.
A fundamental example arises from a finite group $G$ acting on a ring $R$.
In this case, there is a canonical isomorphism $\Mod(R)^G\cong\Mod(RG)$, where $RG$ is the skew group ring.
Thus, equivariantization provides a categorical approach to study homological structures over skew group rings.

These two frameworks lead to a natural problem.
Suppose that a finite group $G$ acts on a Grothendieck category $\calA$ equipped with a Hovey triple $\calM=(\calC,\calW,\calF)$.
There are then two possible constructions.
One may first lift $\calM$ to a triple $\calM^G=(\calC^G,\calW^G,\calF^G)$ on the equivariant category $\calA^G$ and form the homotopy category $\Ho(\calM^G)$.
Alternatively, one may first form the triangulated homotopy category $\Ho(\calM)$ and then consider its equivariant category $\Ho(\calM)^G$.
This raises the following questions:
\begin{intque}\label{quea}
When does $\calM$ induce an abelian model structure on $\calA^G$?
Does $\Ho(\calM)^G$ admit a canonical triangulated structure?
How are $\Ho(\calM^G)$ and $\Ho(\calM)^G$ related?
\end{intque}

The classical transfer theory for model structures along adjunctions does not directly answer the first question, since a transferred model structure need not remain abelian.
Related transfer problems for abelian model structures have been studied by Dalezios and Psaroudakis \cite{2023DP}.
In particular, they establish sufficient conditions for right- and left-lifting hereditary abelian model structures along adjoint functors, formulated in terms of complete cotorsion pairs and suitable acyclicity conditions, with applications to lifting recollements of abelian categories to recollements of the associated homotopy categories.
Our approach is complementary and is tailored to the Frobenius setting.
More precisely, we develop a transfer theory for Hovey triples along Frobenius functors.
Recall that a functor $F$ is called Frobenius if it fits into a Frobenius pair $(F, H)$, that is, both $(F,H)$ and $(H,F)$ are adjoint pairs; see \cite{SE02,1999CJG,Morita65}.
The symmetry of a Frobenius pair allows us to control simultaneously the two cotorsion pairs defining a Hovey triple and to formulate the transfer conditions directly in terms of the classes of the Hovey triple.
In particular, Theorem \ref{HOVEY case1} gives a characterization of when a Hovey triple is right-induced along a faithful Frobenius functor, reformulating the general right acyclicity condition in terms of the composite $H F$ and the unit-counit structure of the adjunction.

The equivariant setting fits naturally into this framework.
Indeed, the induction and forgetful functors form a Frobenius pair $\Ind : \calA \rightleftarrows \calA^G : U$; see \cite[Lemma 4.6]{2010DGNO}.
When $|G|$ is invertible in $\calA$, the counit $\Ind \circ U \to \id_{\calA^G}$ splits.
This allows the lifting problem to be reduced to the $G$-invariance of the classes in the original Hovey triple.

Our first main result gives a positive answer to Question \ref{quea}.
More precisely, we have the following result.

\begin{intthm} \label{retract (Q^G,W^G,R^G)}
Let $\mathcal A$ be a Grothendieck category with enough projectives, let $G$ be a finite group acting on $\mathcal A$ such that $|G|$ is invertible in $\mathcal A$, and let $\mathcal M = (\mathcal C, \mathcal W, \mathcal F)$ be a cofibrantly generated hereditary Hovey triple in $\mathcal A$.
If any two of the three classes $\mathcal C$, $\mathcal W$ and $\mathcal F$ are $G$-invariant, then the following statements hold:
\begin{rqm}
\item
$
\mathcal M^G = (\mathcal C^G, \mathcal W^G, \mathcal F^G)
$
forms a cofibrantly generated hereditary Hovey triple in $\mathcal A^G$.

\item The equivariant category $\Ho(\mathcal M)^G$ of the homotopy category $\Ho(\mathcal M)$ is triangulated and admits a unique canonical triangulated structure.

\item The induced functor
\[
\Ho(\gamma_{\mathcal A}^G):
\Ho(\mathcal M^G) \to \Ho(\mathcal M)^G
\]
is a triangle equivalence up to retracts. Moreover, $\Ho(\gamma_{\mathcal A}^G)$ is a triangle equivalence if the stable category
$\St_{\omega^G}(\mathcal C^G \cap \mathcal F^G)$
is idempotent complete.
\end{rqm}
\end{intthm}

To illustrate Theorem \ref{retract (Q^G,W^G,R^G)}, we consider the Hovey triple
\[
\mathfrak{PGF}(R)=(\PGF(R),\PGF(R)^\bot,\Mod(R))
\]
in $\Mod(R)$; see Example \ref{exa of Goren}.
Then Theorem \ref{retract (Q^G,W^G,R^G)}, together with the canonical isomorphism $\Mod(R)^G \cong \Mod(RG)$, gives a comparison functor
\[
\Theta : \Ho(\mathfrak{PGF}(RG)) \to \Ho(\mathfrak{PGF}(R))^G,
\]
which is a triangle equivalence up to retracts; see Corollary \ref{cor of PGF}.
This example illustrates concretely how the lifting and comparison results of Theorem \ref{retract (Q^G,W^G,R^G)} behave in the module-theoretic setting. Similar illustrations may be obtained from Gorenstein projective, Gorenstein injective, and Gorenstein flat Hovey triples whenever the corresponding model structures exist.

Now suppose we have a Quillen adjunction $(F, H)$ between two categories $\mathcal A$ and $\mathcal B$, each equipped with its own Hovey triple.
This gives derived functors $\mathbb{L}F$ and $\mathbb{R}H$ between their homotopy categories.
Assume further that both $\mathcal A$ and $\mathcal B$ carry a $G$-action, and that the Hovey triples satisfy the invariance conditions of Theorem \ref{retract (Q^G,W^G,R^G)}.
Then the model structures lift to the equivariant categories $\mathcal A^G$ and $\mathcal B^G$.
On one hand, the lifted adjunction $(F^G, H^G)$ yields a total left derived functor
\[
\mathbb{L}(F^G) : \operatorname{Ho}(\mathcal M_{\mathcal A}^G) \to \operatorname{Ho}(\mathcal M_{\mathcal B}^G);
\]
see Proposition \ref{Quillen adj L(F^G)}.
On the other hand, the derived functor $\mathbb{L}F$ inherits a $G$-functor structure and hence induces an equivariant functor
\[
\mathbb{L}(F)^G : \operatorname{Ho}(\mathcal M_{\mathcal A})^G \to \operatorname{Ho}(\mathcal M_{\mathcal B})^G
\]
on the corresponding equivariant categories; see Proposition \ref{G-functor for LF}.
A natural question arises:

\begin{intque}\label{queb}
When are the two derived functors $\mathbb{L}(F^G)$ and $\mathbb{L}(F)^G$ compatible via the comparison functors $\Theta_{\mathcal A}$ and $\Theta_{\mathcal B}$?
Equivalently, when does the square
\[
\xymatrix@C=1.5CM{
   \operatorname{Ho}(\mathcal M_{\mathcal A}^G)  \ar[d]_-{\operatorname{Ho}(\gamma_{\mathcal A}^G)}  \ar[r]^-{\mathbb{L}(F^G)}
&  \operatorname{Ho}(\mathcal M_{\mathcal B}^G)  \ar[d]^-{\operatorname{Ho}(\gamma_{\mathcal B}^G)}             \\
   \operatorname{Ho}(\mathcal M_{\mathcal A})^G                                 \ar[r]^-{\mathbb{L}(F)^G}
&  \operatorname{Ho}(\mathcal M_{\mathcal B})^G                                            }
\]
commute up to natural isomorphism?
Moreover, when do the horizontal arrows become triangle equivalences?
\end{intque}

To address these questions, we need two types of conditions.
Commutativity of the diagram is not automatic; it requires additional compatibility between the adjunction $(F,H)$, the group actions, and the Hovey triples.
To determine when the horizontal arrows become equivalences, we rely on a general criterion from Theorem \ref{Ho(M_A)=Ho(M_B)}.
Our second main theorem provides sufficient conditions that guarantee both the commutativity of the square and, under stronger hypotheses, the triangle equivalence of the horizontal arrows.

\begin{intthm} \label{thmB homotopy square}
Let $\calA$ and $\calB$ be Grothendieck categories with enough projectives.
Let $G$ be a finite group acting both on $\calA$ and $\calB$ via $(\rho^\calA, \theta^\calA)$ and $(\rho^\calB, \theta^\calB)$, respectively, such that $|G|$ is invertible in both categories.
Let $\calM_\calA = (\calC_\calA, \calW_\calA, \calF_\calA)$ and $\calM_\calB = (\calC_\calB, \calW_\calB, \calF_\calB)$ be cofibrantly generated hereditary Hovey triples in $\calA$ and $\calB$, respectively, such that in each triple two of the three classes are $G$-invariant.
Suppose that $(F, H)$ is a Quillen adjunction between $\calA$ and $\calB$ with $F$ or $H$ a $G$-functor.
If $(H,F)$ is an adjoint pair and $F(\calF_\calA \cap \calW_\calA) \subseteq \calW_\calB$, then there exist a total left derived functor $\mathbb{L}(F^G)$ and an equivariant functor $\mathbb{L}(F)^G$ such that the diagram
\[
\xymatrix@C=1.2CM{
   \Ho(\calM_\calA^G)
\ar[d]_-{\Ho(\gamma_{\calA}^G)}^-{ }
\ar[r]^-{\mathbb{L}(F^G)}
&  \Ho(\calM_\calB^G)
\ar[d]^-{\Ho(\gamma_\calB^G)}_-{ }             \\
   \Ho(\calM_\calA)^G
\ar[r]^-{\mathbb{L}(F)^G}
&  \Ho(\calM_\calB)^G                                            }
\]
commutes up to natural isomorphism.
If furthermore, $F$ and $H$ are faithful, $H(\calC_\calB \cap \calW_\calB) \subseteq \calW_\calA$, $\cok(\eta_X) \in \calW_\calA$ for every $X \in \calC_\calA$, and $\ker(\varepsilon_Y) \in \calW_\calB$ for every $Y \in \calF_\calB$, then both $\mathbb{L}(F^G)$ and $\mathbb{L}(F)^G$ are triangle equivalences.
\end{intthm}

We illustrate Theorem \ref{thmB homotopy square} in module categories equipped with PGF Hovey triples.
Two natural sources of the required two-sided adjunctions are Frobenius bimodules and stable equivalences of adjoint type.

\begin{rqm}
\item \emph{Frobenius bimodules.}
A classical result \cite{1999CJG} shows that Frobenius pairs between module categories arise essentially from tensor functors associated with Frobenius bimodules.
When the module categories are equipped with their PGF Hovey triples and the additional hypotheses concerning the group actions, faithfulness, and the unit and counit are satisfied, Theorem \ref{thmB homotopy square} yields a commutative square relating the derived tensor functors before and after equivariantization.
Thus, the corresponding triangle equivalence between the PGF homotopy categories lifts to the equivariant level; see Corollary \ref{homotopy square for F-bimodules}.

\item \emph{Stable equivalences of adjoint type.}
Motivated by the notion introduced by Xi \cite[3.2]{2008Xi}, we consider stable equivalences of adjoint type in the generalised sense, defined by a pair of bimodules whose associated tensor functors are mutually adjoint on both sides.
This setting supplies the two-sided Quillen adjunctions, faithfulness, and the required conditions on the unit and counit in a natural way.
Consequently, Theorem \ref{thmB homotopy square} again produces a commutative homotopy square for the PGF Hovey triples.
In particular, a stable equivalence of adjoint type between the base rings lifts to the corresponding skew group rings; see Corollary \ref{homotopy square for adjoint type}.
\end{rqm}

These two illustrations clarify the content of Theorem \ref{thmB homotopy square} in a concrete module-theoretic setting.
Under the stated hypotheses, taking total derived functors is compatible, up to natural isomorphism, with passing to equivariant objects.
Moreover, triangle equivalences between the PGF homotopy categories of the base rings can be transported to the corresponding equivariant categories, or equivalently, to the module categories over the skew group rings.

The paper is organised as follows.
Section~\ref{sec1} recalls the necessary background on cotorsion pairs, abelian model structures, Hovey triples, group actions, and equivariant categories.
Section~\ref{section Frobenius pair} develops the transfer theory for Hovey triples along Frobenius functors, which provides the technical foundation for the lifting results.
Section~\ref{Proofs of Theorems A} proves Theorem \ref{retract (Q^G,W^G,R^G)}, establishing the lifting of cofibrantly generated hereditary Hovey triples to equivariant categories and comparing the resulting homotopy categories.
Section~\ref{Proofs of Theorem B} studies the compatibility of derived functors with equivariantization and proves Theorem \ref{thmB homotopy square}.

Section~\ref{Abelian model structures on Mod(RG) int} illustrates the main results using PGF Hovey triples.
In Subsection~\ref{Abelian model structures on Mod(RG)}, we illustrate Theorem \ref{retract (Q^G,W^G,R^G)} by comparing the PGF homotopy category over a skew group ring $RG$ with the equivariant category of the PGF homotopy category over $R$.
In Subsections~\ref{Homotopy squares via bimodule} and \ref{Homotopy squares via stable equivalence}, we illustrate Theorem \ref{thmB homotopy square} through homotopy squares induced by a Frobenius bimodule and by a stable equivalence of adjoint type, respectively.
The appendix contains the proofs of Lemmas \ref{admiss for Ho(M)} and \ref{admiss for St(C cap F)}.

\section{Preliminaries} \label{sec1}
\noindent In this preliminary section, we mainly fix some notation, recall relevant notions and collect some necessary facts.
{\sf Throughout the paper}, all functors are covariant, all rings are nonzero associative rings with identity, and all modules are unitary. For a ring $R$, an $R$-module means a left $R$-module, and we denote by $\Mod(R)$ the category of left $R$-modules; right $R$-modules are viewed as modules over the opposite ring $R^{\op}$.

\subsection{Preliminaries on (abelian) model structures and Hovey triples}
\label{Preliminaries}
In this section we mainly collect the necessary background on model structures with emphasis on abelian model structures and Hovey triples. We recall some standard facts about homotopy categories associated with hereditary Hovey triples, which will be used throughout the paper.

\begin{bfhpg}[\bf Weak factorization systems]
Let $l: A \to B$ and $r: C \to D$ be morphisms in a category $\calM$. Recall that $l$ is said to have the {\it left lifting property} with respect to $r$ (equivalently, $r$ has the {\it right lifting property} with respect to $l$) if for every pair of morphisms $f: A \to C$ and $g: B \to D$ with $rf = gl$, there exists a morphism $t: B \to C$ such that $f = tl$ and $g = rt$. For a class $\sfC$ of morphisms in $\calM$, denote by $\sfC^{\boxslash}$ the class of all morphisms $r$ in $\calM$ having the right lifting property with respect to every morphism in $\sfC$. The class $^{\boxslash}\sfC$ is defined dually.
Following Bousfield \cite{Bo77}, a pair $(\sfC,\sfD)$ of classes of morphisms in $\calM$ is called a {\it weak factorization system} if $\sfC^{\boxslash}=\sfD$, $\sfC={^{\boxslash}\sfD}$, and every morphism $\alpha$ in $\calM$ can be decomposed as $\alpha = f c$ with $c\in\sfC$ and $f\in\sfD$.
\end{bfhpg}

The concepts of (closed) model structures and model categories were introduced by Quillen \cite{1967Quillen} to develop homotopy theory.

\begin{bfhpg}[\bf Model structures]
A (closed) \emph{model category} is a category $\calM$ having finite limits and colimits, equipped with a \emph{model structure}, that is, a triple $(\Cofib(\calM), \Weq(\calM), \Fib(\calM))$ of classes of morphisms satisfying the following conditions:
\begin{prt}
\item $(\Cofib(\calM), \Weq(\calM)\cap\Fib(\calM))$ and $(\Cofib(\calM)\cap\Weq(\calM), \Fib(\calM))$ are weak factorization systems;
\item $\Weq(\calM)$ is closed under retracts and satisfies the 2-out-of-3 property: if two of the three morphisms $\alpha$, $\beta$, and $\beta\alpha$ are weak equivalences, then so is the third.
\end{prt}

A morphism in $\Cofib(\calM)$ (resp., $\Weq(\calM)$, $\Fib(\calM)$) is called a {\it cofibration} (resp., {\it weak equivalence}, {\it fibration}). A morphism in $\Cofib(\calM)\cap\Weq(\calM)$ (resp., $\Weq(\calM)\cap\Fib(\calM)$) is a {\it trivial cofibration} (resp., {\it trivial fibration}). An object is \emph{cofibrant} if the map from the initial object to it is a cofibration; dually, it is \emph{fibrant} if the map from it to the terminal object is a fibration. An object is \emph{trivial} if the map from the initial object to it is a weak equivalence, or equivalently, if the map from it to the terminal object is a weak equivalence.

Let $\mathcal{I}$ be a set of morphisms in a category $\mathcal{M}$. Denote by $\mathrm{Cell}(\mathcal{I})$ the class of relative $\mathcal{I}$-cell complexes, that is, transfinite compositions of pushouts of morphisms in $\mathcal{I}$, and by $^{\oplus}\mathrm{Cell}(\mathcal{I})$ the class of retracts of such complexes. A model category $\mathcal{M}$ is \emph{cofibrantly generated} if there are sets $\mathcal{I}$ and $\mathcal{J}$ of morphisms such that
\[
\Cofib(\mathcal{M}) = {^{\oplus}\mathrm{Cell}(\mathcal{I})}
\quad\text{and}\quad
\Cofib(\mathcal{M})\cap\Weq(\mathcal{M}) = {^{\oplus}\mathrm{Cell}(\mathcal{J})};
\]
equivalently, $\Fib(\mathcal{M}) = \mathcal{J}^{\boxslash}$ and $\Fib(\mathcal{M})\cap\Weq(\mathcal{M}) = \mathcal{I}^{\boxslash}$. The sets $\mathcal{I}$ and $\mathcal{J}$ are then called the {\it generating cofibrations} and {\it generating trivial cofibrations}, respectively.
\end{bfhpg}

In this paper, we focus mainly on model structures on abelian categories that are compatible with the abelian structure, namely, abelian model structures.

\begin{bfhpg}[\bf Abelian model structures]\label{thm Hovey tripls}
Recall that a model structure on an abelian category $\calA$ is called \emph{abelian} if the following conditions hold:
\begin{prt}
\item every cofibration is a monomorphism with a cofibrant cokernel; and
\item every fibration is an epimorphism with a fibrant kernel.
\end{prt}

By \prpcite[4.2]{HOVEY2002}, this is equivalent to Hovey's original definition.
An \emph{abelian model category} is a bicomplete abelian category equipped with an abelian model structure.
\end{bfhpg}

\begin{bfhpg}[\bf Cotorsion pairs]\label{def of Cot pairs}
As an analogue of torsion pairs, cotorsion pairs are defined using the Ext functor instead of the Hom functor.
Recall that, following Enochs and Jenda \cite{rha}, a pair $(\calC,\calD)$ of classes of objects in an abelian category $\calA$ is called a \emph{cotorsion pair} if $\calC^{\bot}=\calD$ and $^{\bot}\calD=\calC$, where
\begin{align*}
\calC^{\bot} &= \{ M \in \calA \mid \Ext_{\calA}^{1}(C,M)=0 \text{ for all } C\in\calC \}, \\
^{\bot}\calD &= \{ M \in \calA \mid \Ext_{\calA}^{1}(M,D)=0 \text{ for all } D\in\calD \}.
\end{align*}

A cotorsion pair $(\calC,\calD)$ is \emph{complete} if for every $M\in\calA$ there exist short exact sequences
\[ 0\to D\to C\to M\to 0 \quad \textrm{and} \quad 0\to M\to D'\to C'\to 0 \]
with $D,D'\in\calD$ and $C,C'\in\calC$.

For any class $\calS$ of objects in $\calA$, $(^\bot(\calS^\bot),\calS^\bot)$ forms a cotorsion pair, called the cotorsion pair \emph{cogenerated} by $\calS$. Similarly, $(^\bot\calS,(^\bot\calS)^\bot)$ forms a cotorsion pair, called the cotorsion pair \emph{generated} by $\calS$. In particular, by Saor\'{\i}n and \v{S}\v{t}ov\'{\i}\v{c}ek \corcite[2.15(3)]{2013Sto}, if $\calA$ is a Grothendieck category with enough projectives and $\calS$ is a set, then $(^\bot(\calS^\bot),\calS^\bot)$ is a complete cotorsion pair in $\calA$.

A cotorsion pair $(\calC,\calD)$ in $\calA$ is called \emph{resolving} if $\calC$ is closed under kernels of epimorphisms. Dually, one defines \emph{coresolving} cotorsion pairs. A cotorsion pair is \emph{hereditary} if it is both resolving and coresolving.

Let $(\calC,\calD)$ be a complete cotorsion pair in $\calA$. Becker \cite{Be14} shows that $(\calC,\calD)$ is hereditary if and only if it is resolving, equivalently, if and only if it is coresolving.
\end{bfhpg}

\begin{rmk}\label{exact terminology}
We sometimes consider cotorsion pairs in exact categories. The standard reference for exact categories is \cite{2010ExaCat}. The definition of cotorsion pairs in exact categories is the same as in abelian categories, with the only difference being the terminology. When the underlying category is exact, we still use the terminology of short exact sequences.
\end{rmk}

The following result, now known as \emph{Hovey's correspondence}, is a central result in abelian model category theory. Recall that a class of objects in $\mathcal{A}$ is called \emph{thick} if it is closed under direct summands, extensions, kernels of epimorphisms, and cokernels of monomorphisms.

\begin{thm}[Hovey's correspondence]\label{Hovey cor}
Let $\calA$ be an abelian category. There is a bijective correspondence between
\begin{rqm}
\item abelian model structures $(\Cofib(\calA),\Weq(\calA),\Fib(\calA))$ on $\calA$, and
\item triples $(\calC,\calW,\calF)$ of classes of objects in $\calA$ such that both $(\calC,\calW\cap\calF)$ and $(\calC\cap\calW,\calF)$ are complete cotorsion pairs, and $\calW$ is thick.
\end{rqm}
\end{thm}

\begin{rmk}\label{rmk for Hovey cor}
In Hovey \thmcite[2.2]{HOVEY2002}, it is additionally assumed that $\calA$ is bicomplete. Cui, Rong, and Zhang \cite{2025ZhangPu} pointed out that this condition is not necessary. This fact is also implicitly acknowledged in references such as Gillespie \cite{Gill2011} and \cite{HOWTO}.
\end{rmk}

\begin{bfhpg}[\bf Hovey triples]
By Hovey's correspondence, an abelian model structure on an abelian category $\calA$ can be represented succinctly by a triple of classes of objects in $\calA$ that satisfies the conditions in \ref{Hovey cor}(2). Such a triple is therefore often referred to as an abelian model structure and is called a \emph{Hovey triple}.

A Hovey triple $(\calC,\calW,\calF)$ is said to be \emph{cofibrantly generated} if the corresponding abelian model structure on $\calA$ is cofibrantly generated; when $\calA$ is a Grothendieck category with enough projectives, this is equivalent to both cotorsion pairs $(\calC\cap\calW,\calF)$ and $(\calC,\calW\cap\calF)$ being cogenerated by sets (see \cite[Proposition 1.2.7 and Corollary 1.2.2]{Be14}). A Hovey triple $(\calC,\calW,\calF)$ is \emph{hereditary} if both complete cotorsion pairs $(\calC\cap\calW,\calF)$ and $(\calC,\calW\cap\calF)$ are hereditary.
\end{bfhpg}

Hereditary Hovey triples enjoy the following property, taken from \prpcite[4.2]{2016Gill}.

\begin{prp} \label{C cap F Frobenius}
Let $\calM = (\calC, \calW, \calF)$ be a hereditary Hovey triple in an abelian category $\calA$.
Then $\calC \cap \calF$ is a Frobenius category (with the canonical exact structure) whose projective-injective objects are precisely those in the core $\omega = \calC \cap \calW \cap \calF$.
\end{prp}

The following result, taken from \cite[Lemma 3.3]{2025GillBook}, and its dual will be used frequently.

\begin{lem}\label{fundamental Lemma}
Let $(\calX,\calY)$ be a complete cotorsion pair in an exact category $\calD$ with $\omega=\calX\cap\calY$. Let $\delta$ be a class of objects in $\calD$ closed under biproducts and such that $\omega\subseteq\delta\subseteq\calX$. Consider the solid diagram
\[
\xymatrix{
   Y_A
 \ar@{-->}[d]_-{c} \ar@{>->}[r]
 & X_A
 \ar@{-->}[d]_-{h} \ar@{->>}[r]
 & A
 \ar[d]_-{f}  \\
   Y_B
 \ar@{>->}[r]
 & X_B
 \ar@{->>}[r]
 & B
 }
\]
whose rows are short exact sequences with $X_A\in\calX$ and $Y_B\in\calY$. Then there are dashed arrows $c$ and $h$ making the diagram commute. The lifting pair $(c,h)$ is unique in the stable category $\St_\delta(\calD)$ and depends only on the $\delta$-equivalence class of $f$. More precisely,
for any $f':A\to B$ and corresponding $(c',h')$, if $[f]_\delta=[f']_\delta$, then $[h]_\delta=[h']_\delta$ and $[c]_\delta=[c']_\delta$.
\end{lem}

\begin{bfhpg}[Homotopy categories]\label{tri structure of Ho(A)}
Let $\calM=(\calC,\calW,\calF)$ be a Hovey triple in an abelian category $\calA$.
Recall from \thmcite[2.6]{2016Gill} that there is a category $\Ho(\calM)$, called the \emph{homotopy category} associated to $\calM$, whose objects are those of $\calA$ and whose morphisms are given by
\[
\Hom_{\Ho(\calM)}(X,Y)=\Hom_\calA(RQX,RQY)/\sim,
\]
where $RQ$ is the bifibrant replacement functor, and $f\sim g$ iff $g-f$ factors through an object in the core $\calC\cap\calW\cap\calF$.
Thus, every morphism in $\Ho(\calM)$ is represented by a morphism $h:RQX\to RQY$ in $\calA$, so that $f=[h]$.
For a more explicit decomposition of such morphisms, see \cite[Proposition 5.15]{2025GillBook}.

There is a canonical localization functor $\gamma_{\calA}:\calA\to\Ho(\calM)$, which is the identity on objects and sends each morphism $f$ to $[RQf]$.
The following localization theorem, due to Gillespie \cite[Theorem 5.21]{2025GillBook}, states that this functor has the usual universal property of localization.
We will use this result frequently in the sequel.
\end{bfhpg}

\begin{thm}\label{loc thm}
Let $\calM=(\calC,\calW,\calF)$ be a Hovey triple in an abelian category $\calA$.
Then the functor $\gamma_{\calA}$ sends every weak equivalence in $\calA$ to an isomorphism in $\Ho(\calM)$.
Moreover, for any category $\calB$ and any functor $F:\calA\to\calB$ that sends weak equivalences to isomorphisms, there is a unique functor $\Ho(F):\Ho(\calM)\to\calB$ making the diagram
\[
\xymatrix{
  \calA
\ar[d]_{F} \ar[r]^-{\gamma_{\calA}}
& \Ho(\calM)
\ar@{-->}[dl]^(0.4){\Ho(F)} \\
  \calB
}
\]
commute.
\end{thm}

\begin{rmk} \label{tri for Ho(M)}
The homotopy category $\Ho(\calM)$ is a triangulated category. In \cite[Section 6]{2025GillBook}, Gillespie gives a detailed description of the triangulated category structure of $\Ho(\calM)$. The shift functor $\mathsf{\Sigma}: \Ho(\calM) \to \Ho(\calM)$ is defined via the Localization Theorem \ref{loc thm}.
More precisely, $\mathsf{\Sigma}$ is the unique functor such that the diagram
\[
\xymatrix{
   \calA
\ar[d]_{\gamma_{\calA} \mathsf{\Sigma}}  \ar[r]^{\gamma_{\calA} \quad}
&  \Ho(\calM)
\ar@{-->}[dl]^(0.4){\mathsf{\Sigma}}   \\
   \Ho(\calM)                      }
\]
commutes, where the functor $\gamma_{\calA} \mathsf{\Sigma}$ is defined as follows:
\begin{itemize}
\item
for an object $X \in \calA$, $\gamma_{\calA} \mathsf{\Sigma} (X) = \mathsf{\Sigma} X$, which is taken from a short exact sequence $X \rightarrowtail W_X \twoheadrightarrow \mathsf{\Sigma} X$ with $W_X \in \calW$;

\item
for any morphism $f \in \Hom_\calA(X, Y)$, $\gamma_{\calA} \mathsf{\Sigma} (f) = \gamma_{\calA} (\mathsf{\Sigma} f) = [RQ \mathsf{\Sigma} f]$, where $\mathsf{\Sigma} f$ is the unique morphism making the following diagram commute
\[
\xymatrix{
  X
\ar[d]_-{f} \ar@{>->}[r]
& W_X
\ar[d] \ar@{->>}[r]
& \mathsf{\Sigma} X
\ar[d]^-{\mathsf{\Sigma} f} \\
  Y
\ar@{>->}[r]
& W_Y
\ar@{->>}[r]
& \mathsf{\Sigma} Y                }
\]
\end{itemize}
For any morphism $f \in \Hom_\calA(X, Y)$, the \emph{standard exact triangle}
\[
\xymatrix@C=0.7cm{
  X
\ar[r]^-{\gamma_{\calA}(f)}
& Y
\rightarrow
 C_f
\rightarrow
\mathsf{\Sigma} X
}
\]
in $\Ho(\calM)$ is obtained from the following pushout diagram:
\[
\xymatrix{
  X
\ar[d]_{f}  \ar@{>->}[r]
& W_X
\ar[d] \ar@{->>}[r]
&  \mathsf{\Sigma} X
\ar@{=}[d]  \\
  Y
\ar@{>->}[r]
& C_f
\ar@{->>}[r]
&  \mathsf{\Sigma} X            }
\]
The class of \emph{exact triangles} in $\Ho(\calM)$ is precisely the class of all triangles
that are isomorphic to standard exact triangles.
\end{rmk}

The following result, taken from \thmcite[4.3]{2016Gill}, describes the relationship between the stable category and the homotopy category.

\begin{thm} \label{hered Hovey correspondence}
Let $\calM = (\calC, \calW, \calF)$ be a hereditary Hovey triple in an abelian category $\calA$.
Then the inclusion $\calC \cap \calF \hookrightarrow \calA$ induces a triangle equivalence
\[
\St_\omega(\calC \cap \calF) \hookrightarrow \Ho(\calM).
\]
Its inverse is the bifibrant replacement functor.
\end{thm}

\subsection{Preliminaries on group actions and equivariant categories}
\label{sec:equivariant}
In this section we recall the basic notions of group actions on categories, invariant classes, equivariant objects, equivariant morphisms, and the equivariant category. We also discuss the induction functor, the equivariant functor, and some of their properties, which will be used later in the construction of equivariant model structures and Hovey triples.

\begin{bfhpg}[\bf Group actions on categories]\label{GrAct on Cat}
Let $G$ be a finite group and $\calA$ a category.
Recall that a (right) \emph{action} $(\rho,\theta)$ of $G$ on $\calA$ consists of the following data:
\begin{itemize}
\item for each $g\in G$, an auto-equivalence $\rho_g:\calA\to\calA$; and
\item for each $g,h\in G$, a natural isomorphism $\theta_{g,h}:\rho_g\rho_h \overset{\cong}{\longrightarrow} \rho_{hg}$,
\end{itemize}
subject to the \emph{2-cocycle condition}, i.e., for all $g,h,k\in G$, the diagram
\[
\xymatrix@R=0.75cm@C=1.5cm{
  \rho_g \rho_h \rho_k
\ar[d]_-{\theta_{g, h} \rho_k}
\ar[r]^-{\rho_g \theta_{h, k}}
& \rho_g \rho_{kh}
\ar[d]^-{\theta_{g, kh}}  \\
  \rho_{hg} \rho_k
\ar[r]_-{\theta_{hg,k}}
& \rho_{khg}}
\]
commutes. The natural isomorphisms $\theta_{g,h}$ are called \emph{composition isomorphisms}.

Every such action admits a natural isomorphism $u:\rho_e \to \Id_{\calA}$, called the \emph{unit} of the action (where $e$ is the identity element of $G$).
If each $\rho_g$ is an automorphism and each $\theta_{g,h}$ is the identity natural transformation (for all $g,h\in G$), the action is called \emph{strict}.

Left actions are defined analogously. The assignments
\[
\lambda_g=\rho_{g^{-1}} \quad\text{and}\quad \gamma_{g,h}=\theta_{g^{-1},h^{-1}}, \qquad \forall\,g,h\in G,
\]
translate a right action into a left action.
\end{bfhpg}

The following example is well known; see e.g. \cite[Example 2.6]{2017Chen} or \cite[Subsection 2.3]{2025Equ}.

\begin{exa} \label{G for RMOD}
Let $R$ be a ring, $G$ a finite group, and $\rho: G \to \Aut(R)^{\op}$ a group homomorphism.
For each $g \in G$, there exists an automorphism
\[
(-)^g: \Mod(R^{\op}) \to \Mod(R^{\op})
\]
sending an $\Rop$-module $M$ to $M^g$; the \emph{twisted} module $M^g$ equals $M$ as an abelian group, but with a new $R^{\op}$-action $\cdot_g$ given by
\[
m \cdot_g r = m \cdot_R \rho(g^{-1})(r),
\]
where $m \in M$, $r \in R$ and $\cdot_R$ denotes the original $R$-module action.
For a morphism $f: M \to N$ in $\Mod(\Rop)$, we set $f^g = f : M^g \to N^g$.
It is easy to check that $(M^h)^g = M^{hg}$; see \cite[Subsection 2.3]{2025Equ}.
Then one has $(-)^g \circ (-)^h = (-)^{hg}$, giving a strict right $G$-action on $\Mod(R^{\mathrm{op}})$.

Similarly, for an $R$-module $N$, we define $N^g = N$ as an abelian group with the new $R$-action $\cdot_g$ given by
\[
r \cdot_g n = \rho(g^{-1})(r) \cdot_R n.
\]
The same computation then shows that $(-)^g \circ (-)^h = (-)^{hg}$, yielding a strict right $G$-action on $\Mod(R)$.
\end{exa}

Let $\calA$ be a category and $G$ a finite group acting on $\calA$ via $(\rho,\theta)$.
Recall that a class $\calC$ of objects in $\calA$ is $G$-\emph{invariant} if $\rho_g\calC\subseteq\calC$ for each $g\in G$.

\begin{lem} \label{D^G G-invariant}
Let $\calA$ be an abelian category and $G$ a finite group acting on $\calA$ via $(\rho,\theta)$.
Suppose that $(\calC,\calD)$ is a cotorsion pair in $\calA$.
Then $\calC$ is $G$-invariant if and only if $\calD$ is $G$-invariant.
\end{lem}

\begin{prf*}
Suppose that $\calC$ is $G$-invariant.
Since $\rho_g:\calA\to\calA$ is an auto-equivalence for each $g\in G$, we have the isomorphism
$\Ext_{\calA}^1(C,\rho_gD) \cong \Ext_{\calA}^1(\rho_{g^{-1}}C,D)$,
where $C\in\calC$ and $D\in\calD$.
The latter vanishes because $\rho_{g^{-1}}C\in\calC$.
Thus, $\rho_gD\in \calC^\bot=\calD$, and hence, $\calD$ is $G$-invariant.
The converse is proved similarly.
\end{prf*}

By Lemma \ref{D^G G-invariant}, we obtain the following result on Hovey triples.

\begin{lem} \label{(Q,W,R) G-invariant}
Let $\calA$ be an abelian category, $G$ a finite group acting on $\calA$ via $(\rho,\theta)$, and $(\calC,\calW,\calF)$ a Hovey triple in $\calA$.
If any two of the three classes $\calC$, $\calW$, $\calF$ are $G$-invariant, then so is the third.
\end{lem}

\begin{prf*}
If $\calC$ and $\calW$ are $G$-invariant, then so is $\calC\cap\calW$. Hence, $\calF$ is $G$-invariant by Lemma \ref{D^G G-invariant}, since $(\calC\cap\calW,\calF)$ is a cotorsion pair. The case where $\calW$ and $\calF$ are $G$-invariant is analogous.

Now, suppose $\calC$ and $\calF$ are $G$-invariant. We show that $\calW$ is $G$-invariant, i.e., $\rho_g W \in \calW$ for all $W\in\calW$ and $g\in G$.
Since $(\calC,\calW\cap\calF)$ is a complete cotorsion pair, there is a short exact sequence $0\to W\to E\to Q\to 0$ with $E\in\calW\cap\calF$ and $Q\in\calC$.
Because $\calW$ is thick, we also have $Q\in\calW$.
As $\rho_g$ is exact, applying it gives short exact sequence
$0\to \rho_g W\to \rho_g E\to \rho_g Q\to 0$.
By Lemma \ref{D^G G-invariant}, both $\calC\cap\calW$ and $\calW\cap\calF$ are $G$-invariant as $(\calC\cap\calW,\calF)$ and $(\calC,\calW\cap\calF)$ are cotorsion pairs. Thus, $\rho_g E\in\calW\cap\calF$ and $\rho_g Q\in\calC\cap\calW$.
Since $\calW$ is thick, it follows that $\rho_g W\in\calW$, as desired.
\end{prf*}

\begin{bfhpg}[\bf The $G$-functor] \label{G fun}
Let $G$ be a finite group acting on two categories $\calA$ and $\calB$, with actions $(\rho^\calA,\theta^\calA)$ and $(\rho^\calB,\theta^\calB)$, respectively.
Following \cite{2025Equ}, a $G$-\emph{functor}
\[
(F,\sigma) : (\calA,\rho^\calA,\theta^\calA) \to (\calB,\rho^\calB,\theta^\calB)
\]
consists of a functor $F : \calA \to \calB$ together with a family of natural isomorphisms
\[
\sigma=\{\sigma_g : F\circ\rho^\calA_g \rightarrow \rho^\calB_g\circ F\}_{g\in G}
\]
making the following diagrams commute
\[
\xymatrix@R=1cm@C=1cm{
   F \rho^\calA_g \rho^\calA_h
\ar[r]^-{\sigma_g \rho^\calA_h}
\ar[d]_-{F \theta^\calA_{g,h}}
&  \rho^\calB_g F \rho^\calA_h
\ar[r]^-{\rho^\calB_g \sigma_h}
&  \rho^\calB_g \rho^\calB_h F
\ar[d]^-{\theta^\calB_{g,h} F }      \\
   F \rho^\calA_{hg}
\ar[rr]^-{\sigma_{hg}}
&&  \rho^\calB_{hg}F                 }
\quad\text{and}\quad
\xymatrix{
   F \rho^\calA_e
\ar[rr]^-{\sigma_e}
\ar[dr]_-{Fu_{\calA}}
&& \rho^\calB_e F
\ar[dl]^-{u_{\calB} F}               \\
&  F                              }
\]
where $u_\calA$ and $u_\calB$ are the units of the $G$-actions on $\calA$ and $\calB$, respectively.

For simplicity, we shall often write $F$ instead of $(F,\sigma)$ for a $G$-functor, whenever the family $\sigma$ is clear from the context.
\end{bfhpg}

\begin{exa} \label{inclusion G-functor}
Let $\calC$ be a $G$-invariant class of objects in $\calA$.
Then the $G$-action $(\calA, \rho, \theta)$ on $\calA$ naturally restricts to a $G$-action $(\calC, \rho, \theta)$ on $\calC$.
In this case, the inclusion functor $\iota:\calC\hookrightarrow\calA$ is a $G$-functor in the natural way.
\end{exa}

\begin{bfhpg}[\bf The $G$-natural transformation] \label{G natural trans}
Let $(F_i, \sigma^i) : (\calA, \rho^{\calA}, \theta^{\calA}) \to (\calB, \rho^{\calB}, \theta^{\calB})$, $i = 1, 2$, be two $G$-functors.
Following \cite{2025Equ}, a natural transformation $\eta : F_1 \Rightarrow F_2$ is called a \emph{$G$-natural transformation} if, for every $g \in G$, the following diagram is commutative:
\[
\xymatrix{
F_1 \rho^{\calA}_g
\ar[r]^-{\eta \rho^{\calA}_g}
\ar[d]_-{\sigma^1_g}
&
F_2 \rho^{\calA}_g
\ar[d]^-{\sigma^2_g}
\\
\rho^{\calB}_g F_1
\ar[r]_-{\rho^{\calB}_g \eta}
&
\rho^{\calB}_g F_2
}
\]
In this case, we write $\eta : (F_1, \sigma^1) \Rightarrow_G (F_2, \sigma^2)$.
\end{bfhpg}

\begin{bfhpg}[\bf The equivariant category]\label{equi cat}
Let $\calA$ be a category endowed with a $G$-action $(\rho,\theta)$.
Following \cite{2025Equ}, a $G$-{\it equivariant object} $(X,\phi)$ of $\calA$ consists of an object $X\in\calA$ together with a family of isomorphisms
\[
\phi=\{\phi_g:\rho_g X \overset{\cong}{\longrightarrow} X\}_{g\in G},
\]
called the {\it equivariant structure} or {\it linearization} of $X$, such that for all $g,h\in G$ the diagram
\begin{equation*} \label{key comm sq1} \tag{\ref{equi cat}.1}
\begin{gathered}
\xymatrix{
   \rho_g \rho_h X
\ar[d]_-{\theta_{g,h,X}}
\ar[r]^-{\rho_g \phi_h}
&  \rho_g X
\ar[d]^-{\phi_g}                                           \\
   \rho_{hg}X
\ar[r]^-{\phi_{hg}}
&  X                                                                          }
\end{gathered}
\end{equation*}
commutes.
When there is no risk of confusion, a $G$-equivariant object $(X,\phi)$ will be denoted simply by $X$.
The commutative diagram (\ref{key comm sq1}) will be referred to as the \emph{cocycle condition} for the linearizations.

Given two $G$-equivariant objects $(X,\phi)$ and $(Y,\psi)$, a morphism $f:(X,\phi)\to(Y,\psi)$ is a morphism $f:X\to Y$ in $\calA$ that is compatible with the linearizations, i.e., such that for every $g\in G$ the diagram
\begin{equation*} \label{key comm sq2} \tag{\ref{equi cat}.2}
\begin{gathered}
\xymatrix{
   \rho_g X
\ar[d]_-{\phi_g}
\ar[r]^-{\rho_g f}
&  \rho_g Y
\ar[d]^-{\psi_g}                      \\
   X
\ar[r]^-{f}
&  Y                                              }
\end{gathered}
\end{equation*}
commutes.

The {\it equivariant category}, denoted by $\calA^G$, is defined to have $G$-equivariant objects as objects and morphisms as above. This category is also called the \emph{equivariantization} of $\calA$ with respect to the $G$-action.
\end{bfhpg}

\begin{rmk} \label{equiv is abe}
If $\calA$ is additive, then so is the equivariant category $\calA^G$.
Moreover, when $\calA$ is abelian, $\calA^G$ is also abelian; we refer to \rmkcite[2.6]{2025Equ} for explicit constructions of the zero object, kernels, and cokernels in $\calA^G$.
Note that a sequence $0 \to (A,\alpha) \to (B,\beta) \to (C,\gamma) \to 0$ in $\calA^G$ is exact if and only if the underlying sequence $0 \to A \to B \to C \to 0$ in $\calA$ is exact.
\end{rmk}

\begin{exa}[\bf Skew group algebra]\label{exa skew}
Let $R$ be a ring, $G$ a finite group, and $\rho: G \to \operatorname{Aut}(R)^{\op}$ a group homomorphism.
The \emph{skew group algebra} $RG$ is defined as a free $R$-module with basis $G$.
Its elements are formal sums $\sum_{g\in G} a_g g$, $a_g\in R$, with multiplication given by
\[
(a g)(b h) = (a\rho(g^{-1})(b))(gh)
\]
for $a,b\in R$ and $g,h\in G$.
When $G$ acts trivially on $R$ (i.e., $\rho(g)=\Id_R$ for all $g\in G$), this reduces to the usual group algebra $R[G]$.

An $(RG)^{\op}$-module is an $R^{\op}$-module $N$  equipped with a right $G$-action satisfying
\[
(n\cdot_R r)g = (n g)\cdot_R \rho(g)(r).
\]
By Example \ref{G for RMOD}, there is a strict right $G$-action on both $\Mod(R^{\op})$ and $\Mod(R)$, giving the equivariant categories $\Mod(R^{\op})^G$ and $\Mod(R)^G$.
It is known that an $(RG)^{\op}$-module can be viewed as a $G$-equivariant object of $\Mod(R^{\op})$.
Indeed, by \exacite[2.12]{2019Sun}, \cite[Example 2.6]{2017Chen} or \cite[Subsection 2.3]{2025Equ}, there is a canonical isomorphism
\[
\Phi:\Mod((RG)^{\op}) \overset{\cong}{\longrightarrow} \Mod(R^{\op})^G,
\]
defined as follows:
\begin{rqm}
\item[$\bullet$] for an $(RG)^{\op}$-module $N$, $\Phi(N)=(N,\psi)$, where each
\[
\psi_g:N^g \overset{\cong}{\longrightarrow} N,\qquad n\mapsto n g,
\]
is an isomorphism of right $R$-modules;
\item[$\bullet$] for a morphism $f$ in $\Mod((RG)^{\op})$, set $\Phi(f)=f$.
\end{rqm}

Similarly, there is a canonical isomorphism
\[
\Psi:\Mod(RG) \overset{\cong}{\longrightarrow} \Mod(R)^G
\]
given by $\Psi(M)=(M,\phi)$ for an $RG$-module $M$, where each
\[
\phi_g : M^g \overset{\cong}{\longrightarrow} M,\qquad m\mapsto g^{-1} m,
\]
is an isomorphism of left $R$-modules, and $\Psi(f)=f$ for morphisms $f$.
\end{exa}

{\sf Throughout the paper}, unless otherwise stated, we regard an $RG$-module as a $G$-equivariant object $(M,\phi)$ of $\Mod(R)$, with $M\in\Mod(R)$ and $\phi_g:M^g\to M$ an $R$-isomorphism for each $g\in G$.
Similarly, an $(RG)^{\op}$-module is viewed as a $G$-equivariant object $(N,\psi)$ of $\Mod(R^{\op})$, with $\psi_g:N^g\to N$ an $R^{\op}$-isomorphism.

For a family $\{(X^i,\phi^i)\}_{i\in I}$ of objects in $\calA^G$, the coproduct is given by the pair $(\coprod_{i\in I}X^i,\Phi)$, where $\coprod_{i\in I}X^i$ is the coproduct of $\{X^i\}_{i\in I}$ and, for each $g\in G$, $\Phi_g$ is the composite of the canonical isomorphisms
\[
\rho_g(\coprod_{i\in I}X^i) \longrightarrow \coprod_{i\in I}\rho_g X^i \longrightarrow \coprod_{i\in I}X^i.
\]
Similarly, colimits of direct systems in $\calA^G$ are constructed in the evident way.
These characterizations make it straightforward to verify the following result concerning the $\textsf{AB}$ axioms.

\begin{lem} \label{prp of AB345}
Let $\calA$ be an abelian category and $G$ a finite group acting on $\calA$.
If $\calA$ satisfies the axioms $\textsf{AB3}$, $\textsf{AB4}$, or $\textsf{AB5}$, then so does the equivariant category $\calA^G$.
\end{lem}

\begin{bfhpg}[\bf The induction functor]\label{frobenius pair (Ind,w)}
Let $\calA$ be a category with small coproducts and $G$ a finite group acting on $\calA$.
There are two adjoint pairs $(\Ind, U)$ and $(U, \Ind)$ between $\calA$ and $\calA^G$. Here $U: \calA^G\to\calA$ is the \emph{forgetful} functor and $\Ind:\calA\to\calA^G$ is the \emph{induction} functor sending an object $X$ to the induced linearization
\[
\Ind(X)=\bigl(\bigoplus_{h\in G}\rho_h X,\; \theta(X)\bigr),
\]
where for each $g\in G$, the isomorphism
\[
\theta(X)_g : \rho_g\bigl(\bigoplus_{h\in G}\rho_h X\bigr)
\rightarrow
              \bigoplus_{h\in G}\rho_h X
\]
is induced diagonally by the isomorphisms
$\theta_{g, hg^{-1},X}:\rho_g\rho_{hg^{-1}}X \to \rho_h X$.
For details, we refer the reader to \cite[Lemma 3.8]{2014Elag} or \cite[Lemma 4.6]{2010DGNO}.
\end{bfhpg}

\begin{rmk} \label{watts thm}
In the module category, the induction functor admits a more explicit description.
Indeed, let $R$ be a ring, $G$ a finite group, and $\rho:G\to\operatorname{Aut}(R)^{\op}$ a group homomorphism.
Then the induction functor $\Ind:\Mod(R)\to\Mod(RG)$ is naturally isomorphic, by the Eilenberg--Watts theorem, to the tensor functor $RG\otimes_R -$.
In particular, $\Ind(X)\cong RG\otimes_R X$ for any $R$-module $X$.
Dually, for any $R^{\op}$-module $Y$, we have $\Ind(Y)\cong Y\otimes_R RG$.
\end{rmk}

The induction functor can be used to produce a family of generators for the equivariant category.
For the reader's convenience, we provide a proof of the following result.

\begin{lem} \label{prp of generator}
Let $\calA$ be an abelian category satisfying axiom $\textsf{AB3}$, and let $G$ be a finite group acting on $\calA$.
If $\calA$ has a family of generators, then so does $\calA^G$.
The same holds for a family of projective generators.
\end{lem}

\begin{prf*}
Let $\calS$ be a family of generators of $\calA$.
We claim that
\[
\Ind(\calS)=\{\Ind(S)\mid S\in\calS\}
\]
is a family of generators of $\calA^G$.
Indeed, let $f:(X,\phi)\to(Y,\psi)$ be a non-zero morphism in $\calA^G$.
Then $f:X\to Y$ is non-zero in $\calA$.
Since $\calS$ generates $\calA$, there exist $S\in\calS$ and a morphism $g:S\to X$ such that $f\circ g\neq 0$.
Since $(\Ind,U)$ is an adjoint pair, we have a natural isomorphism
\[
\Hom_{\calA^G}(\Ind(S),(X,\phi)) \cong \Hom_{\calA}(S,X),
\]
so $g$ corresponds to a morphism $\widetilde{g}:\Ind(S)\to(X,\phi)$.
Then the composition $f\circ\widetilde{g}:\Ind(S)\to(Y,\psi)$ is non-zero in $\calA^G$.

If the objects of $\calS$ are projective, then $\Ind(\calS)$ is a family of projective generators of $\calA^G$, since $\Ind$, being left adjoint to the exact functor $U$, preserves projective objects.
\end{prf*}

As an immediate consequence of Lemmas \ref{prp of AB345} and \ref{prp of generator}, we have:

\begin{prp} \label{prp of Grothendieck}
Let $\calA$ be an abelian category and $G$ a finite group acting on $\calA$.
If $\calA$ is a Grothendieck category, then so is $\calA^G$.
Moreover, if $\calA$ is a Grothendieck category with a family of projective generators, then $\calA^G$ is as well.
\end{prp}

\begin{bfhpg}[\bf The equivariant functor] \label{equi fun}
Suppose that a finite group $G$ acts on two categories $\calA$ and $\calB$, with actions $(\rho^\calA,\theta^\calA)$ and $(\rho^\calB,\theta^\calB)$, respectively.
By \cite[Lemma 2.9]{2025Equ}, a $G$-functor $(F,\sigma):(\calA,\rho^\calA,\theta^\calA)\to(\calB,\rho^\calB,\theta^\calB)$ induces an \emph{equivariant functor} $F^G:\calA^G\to\calB^G$ between the corresponding equivariant categories, defined as follows:
\begin{rqm}
\item[$\bullet$] for an object $(X,\phi)$ in $\calA^G$, $F^G(X,\phi)=(FX,\overline{\phi})$, where
\[
\overline{\phi}=\{\overline{\phi}_g:\rho^\calB_g FX
\overset{(\sigma_{g,X})^{-1}}{\longrightarrow} F\rho^\calA_g X
\overset{F\phi_g}{\longrightarrow} FX\}_{g\in G};
\]
\item[$\bullet$] for a morphism $f:(X,\phi)\to(Y,\psi)$ in $\calA^G$, $F^G(f)=F(f)$.
\end{rqm}
For more details, we refer to the proof of \cite[Lemma 2.9]{2025Equ}. We mention that if the functor $F$ is additive, then so is $F^G$.
It is straightforward to check that $F^G$ commutes with the forgetful functor $U$, i.e., the diagram
\[
\xymatrix@C=1cm{
   \calA^G
\ar[r]^-{F^G}
\ar[d]_-{U}
&  \calB^G
\ar[d]^-{U} \\
   \calA
\ar[r]^-{F}
&  \calB
}
\]
commutes.
\end{bfhpg}

The following results are taken from Sun \cite{2019Sun}.

\begin{lem} \label{G-fun comp}
Let $(F,\sigma):(\calA,\rho^\calA,\theta^\calA)\to(\calB,\rho^\calB,\theta^\calB)$ and $(H,\tau):(\calB,\rho^\calB,\theta^\calB)\to(\calC,\rho^\calC,\theta^\calC)$ be $G$-functors.
Then the composition
\[
(H\circ F,\{\tau_g F \circ H\sigma_g\}_{g\in G}):
(\calA,\rho^\calA,\theta^\calA)\to(\calC,\rho^\calC,\theta^\calC)
\]
is again a $G$-functor.
Moreover, $(H\circ F)^G = H^G\circ F^G$.
\end{lem}

\begin{lem} \label{G-frobenius pair}
Let $F : \calA \rightleftarrows \calB : H$ be an adjoint pair.
Then $F$ is a $G$-functor if and only if $H$ is a $G$-functor. In this case, the induced equivariant functors form an adjoint pair $F^G : \calA^G \rightleftarrows \calB^G : H^G$.
\end{lem}

\begin{lem} \label{G-equivalence}
Let $(F,\sigma):(\calA,\rho^\calA,\theta^\calA)\to(\calB,\rho^\calB,\theta^\calB)$ be a $G$-functor.
If $F:\calA\to\calB$ is an equivalence of categories, then so is $F^G:\calA^G\to\calB^G$.
\end{lem}

\section{Transfer of Hovey triples along Frobenius functors}\label{section Frobenius pair}
\noindent
In this section, we show that, under suitable assumptions, along a Frobenius functor $F : \calA \to \calB$ between Grothendieck categories with enough projectives, cofibrantly generated Hovey triples in $\calA$ induce cofibrantly generated Hovey triples in $\calB$ (Theorem \ref{HOVEY case1}), and that a Quillen adjunction $(F, H)$ can be promoted to a Quillen equivalence (Theorem \ref{Ho(M_A)=Ho(M_B)}).
These results will be used later to obtain cofibrantly generated Hovey triples in the equivariant category from those in the base category, and to prove Theorem \ref{retract (Q^G,W^G,R^G)}.

We record the following transfer criterion, which will be used in this section. As mentioned in the introduction, it determines exactly when model structures are lifted along adjunctions.

\begin{bfhpg}[\bf Right-induced model structures]\label{ex for rim}
Let $(F,H)$ be an adjoint pair between Grothendieck categories $\calA$ and $\calB$. For a class $\mathcal K$ of morphisms in $\mathcal A$, let $H^{-1}(\mathcal K)$ denote the class of morphisms $f$ in $\mathcal B$ such that $H(f)\in\mathcal K$. If $(\Cofib(\mathcal A),\Weq(\mathcal A),\Fib(\mathcal A))$ is a cofibrantly generated model structure on $\mathcal A$, then, by \cite[Corollary 2.7]{20GKR} or \cite[Corollary 3.3.4]{KMEB}, there is a model structure
\[
\bigl({}^{\boxslash}H^{-1}(\Weq(\mathcal A)\cap\Fib(\mathcal A)),\,
H^{-1}(\Weq(\mathcal A)),\,
H^{-1}(\Fib(\mathcal A))\bigr)
\]
on $\mathcal B$, called the \emph{right-induced} model structure, if and only if the \emph{right acyclicity condition} holds:
\[
{}^{\boxslash}H^{-1}(\Fib(\mathcal A)) \subseteq H^{-1}(\Weq(\mathcal A)). \tag{RAC}
\]
\end{bfhpg}

The next standard result shows that the condition (RAC) can be checked directly on the generating set of acyclic cofibrations. For the reader's convenience, we provide a proof.

\begin{lem}\label{acyclicity condition}
Let $F:\calA\rightleftarrows\calB:H$ be an adjoint pair between Grothendieck categories $\mathcal A$ and $\mathcal B$, and let $(\Cofib(\mathcal A),\Weq(\mathcal A),\Fib(\mathcal A))$ be a cofibrantly generated model structure on $\mathcal A$, with generating sets of cofibrations and trivial cofibrations denoted by $\mathcal I_{\mathcal A}$ and $\mathcal J_{\mathcal A}$, respectively.
Then the following are equivalent.
\begin{eqc}
  \item The right acyclicity condition \text{(RAC)} holds.
  \item Every relative $F(\mathcal J_{\mathcal A})$-cell complex lies in $H^{-1}(\Weq(\mathcal A))$.
\end{eqc}
\end{lem}

\begin{prf*}
By assumption, we have $\Fib(\mathcal A)=\mathcal J_{\mathcal A}^{\boxslash}$. Since $(F,H)$ is an adjunction, a morphism $g$ in $\mathcal B$ has the right lifting property with respect to $F(j)$ for $j\in\mathcal J_{\mathcal A}$ if and only if $H(g)$ has the right lifting property with respect to $j$ for $j\in\mathcal J_{\mathcal A}$. Hence,
\[
H^{-1}(\Fib(\mathcal A)) = H^{-1}(\mathcal J_{\mathcal A}^{\boxslash}) = F(\mathcal J_{\mathcal A})^{\boxslash}.
\]
By the small object argument (see \cite[Theorem 2.1.14]{Hovey1999}), the left lifting complement of $F(\mathcal J_{\mathcal A})^{\boxslash}$ is precisely the class of retracts of relative $F(\mathcal J_{\mathcal A})$-cell complexes, that is,
$
{}^{\boxslash}(F(\mathcal J_{\mathcal A})^{\boxslash}) = {^{\oplus}\mathrm{Cell}(F(\mathcal J_{\mathcal A}))}.
$
Now, condition (RAC) states that
$
{}^{\boxslash}H^{-1}(\Fib(\mathcal A)) \subseteq H^{-1}(\Weq(\mathcal A)).
$
By the above equalities, this is equivalent to
$
^{\oplus}\mathrm{Cell}(F(\mathcal J_{\mathcal A})) \subseteq H^{-1}(\Weq(\mathcal A)).
$
Since $\Weq(\mathcal A)$ is closed under retracts and $H$ preserves retracts, the above inclusion is equivalent to
$
\mathrm{Cell}(F(\mathcal J_{\mathcal A})) \subseteq H^{-1}(\Weq(\mathcal A)).
$
Thus, \eqclbl{i} and \eqclbl{ii} are equivalent.
\end{prf*}

\begin{rmk}\label{cofibrantly generated}
Under the equivalent conditions of Lemma \ref{acyclicity condition}, the right-induced model structure
$(^{\boxslash}H^{-1}(\Weq(\mathcal A)\cap\Fib(\mathcal A)),H^{-1}(\Weq(\mathcal A)),H^{-1}(\Fib(\mathcal A)))$
on $\mathcal B$ is cofibrantly generated, with $F(\mathcal I_{\mathcal A})$ and $F(\mathcal J_{\mathcal A})$ as the sets of generating cofibrations and generating acyclic cofibrations, respectively; see \cite[Theorem 11.3.2]{HPS2003}.
\end{rmk}

\begin{bfhpg}[\bf Frobenius pairs] \label{def of Frobenius pairs}
Let $\calA$ and $\calB$ be Grothendieck categories.
Recall from Morita \cite{Morita65} that an adjoint pair $(F,H)$ between $\calA$ and $\calB$ is called a (classical) \emph{Frobenius pair} if $(H,F)$ is also an adjoint pair; such pairs are also known as ambidextrous adjunctions in the literature (see, e.g., \cite{2020BD}).
A functor $F$ is called a \emph{Frobenius functor} if it is part of a Frobenius pair $(F,H)$; in that case, $H$ is also a Frobenius functor. Thus, Frobenius functors always occur in pairs.
\end{bfhpg}

It is clear that both $F$ and $H$ are exact whenever $(F,H)$ is a Frobenius pair.
Hence, the following is a consequence of \lemcite[5.1]{HolmQUIVER}.

\begin{lem} \label{Ext(F,H)}
Let $(F,H)$ be a Frobenius pair between Grothendieck categories $\calA$ and $\calB$.
Then for all $X\in\calA$, $Y\in\calB$ and $i\ge 1$, there are isomorphisms
\[
\Ext_{\calB}^{i}(F(X),Y) \cong \Ext_{\calA}^{i}(X,H(Y))
\quad\text{and}\quad
\Ext_{\calA}^{i}(H(Y),X) \cong \Ext_{\calB}^{i}(Y,F(X)).
\]
\end{lem}

We mention that checking the right acyclicity condition in Lemma \ref{acyclicity condition} directly is often difficult, as it requires verifying all relative cell complexes. In the following, we show that, when $(F,H)$ is a Frobenius pair, this condition admits a much simpler reformulation.

\begin{setup}
In the rest of this section, unless otherwise stated, we always assume the following:
\begin{rqm}
\item[$\bullet$]
$\calA$ and $\calB$ are Grothendieck categories with enough projectives;

\item[$\bullet$]
$F : \calA \rightleftarrows \calB : H$ is a Frobenius pair with unit $\eta:\Id_{\calA}\to HF$ and counit $\varepsilon:FH\to\Id_{\calB}$;

\item[$\bullet$] for a class $\calE$ of objects in $\calA$, set $H^{-1}(\calE):=\{E\in\calB\mid H(E)\in\calE\}$.
\end{rqm}
\end{setup}

\begin{prp} \label{HF condition}
Let $\calM_{\calA}=(\Cofib(\mathcal A),\Weq(\mathcal A),\Fib(\mathcal A))$ be a cofibrantly generated abelian model structure on $\mathcal A$. If $HF$ preserves cofibrant objects in $\calM_{\calA}$, then the following are equivalent:
\begin{eqc}
\item There is a right-induced model structure on $\mathcal B$
\[
\calM_{\calB}=\bigl({}^{\boxslash}H^{-1}(\Weq(\mathcal A)\cap\Fib(\mathcal A)),\,
H^{-1}(\Weq(\mathcal A)),\,
H^{-1}(\Fib(\mathcal A))\bigr).
\]
\item $HF$ preserves trivial cofibrant objects in $\calM_{\calA}$.
\end{eqc}
In this case, every cofibrant object in $\calM_{\calB}$ lies in $H^{-1}(\calC_{\calA})$, where $\calC_{\calA}$ is the class of cofibrant objects in $\calM_{\calA}$.
\end{prp}

\begin{prf*}
Let $\mathcal I_{\calA}$ and $\mathcal J_{\calA}$ be generating sets of cofibrations and trivial cofibrations of $\mathcal M_{\mathcal A}$, respectively.

\proofofimp{i}{ii}
Note that the right acyclicity condition (RAC) holds; see \ref{ex for rim}. Equivalently, by Lemma \ref{acyclicity condition}, every relative $F(\mathcal J_{\calA})$-cell complex lies in $\Weq(\mathcal B)$.

Let $X$ be a trivial cofibrant object in $\calM_{\calA}$. Then $0\to X$ is a trivial cofibration in $\mathcal A$. Since $\mathcal M_{\mathcal A}$ is cofibrantly generated, $0\to X$ is a retract of a relative $\mathcal J_{\calA}$-cell complex. Applying the exact functor $F$, we obtain that $0\to F(X)$ is a retract of a relative $F(\mathcal J_{\calA})$-cell complex. Hence, $0\to F(X)$ is a weak equivalence in $\calM_{\calB}$. Therefore,
$
H(0\to F(X))=(0\to HF(X))
$
is a weak equivalence in $\calM_{\calA}$. Thus, $HF(X)$ is a trivial object. By the assumption that $HF$ preserves cofibrant objects, $HF(X)$ is cofibrant. Hence, it is a trivial cofibrant object in $\calM_{\calA}$.

\proofofimp{ii}{i}
It suffices to prove that the right acyclicity condition (RAC) holds; see \ref{ex for rim}. By Lemma \ref{acyclicity condition}, it is enough to show that every relative
$F(\mathcal J_\calA)$-cell complex is sent by $H$ to a weak equivalence in
$\calM_\calA$.

Let $j: X\to Y$ be in $\mathcal J_{\calA}$. Since $j\in \Cofib(\mathcal A)\cap\Weq(\mathcal A)$, it is a monomorphism with
$\operatorname{coker}(j)$ a trivial cofibrant object in $\calM_\calA$.
Consider a pushout of $F(j)$:
\[
\xymatrix{
  0
\ar[r]^{}
& F(X)
\ar[d]_{}
\ar[r]^-{F(j)}
& F(Y)
\ar[d]_-{}
\ar[r]^-{}
& F(\coker(j))
\ar@{=}[d]_{}
\ar[r]^{}
& 0                               \\
  0
\ar[r]^{}
& U
\ar[r]^-{i}
& V
\ar[r]^{}
& F(\coker(j))
\ar[r]^{}
& 0
}
\]
Applying $H$ yields
\[
\xymatrix{
  0
\ar[r]^{}
& HF(X)
\ar[d]_{}
\ar[r]^-{HF(j)}
& HF(Y)
\ar[d]_{}
\ar[r]^{}
& HF(\coker(j))
\ar@{=}[d]_{}
\ar[r]^{}
& 0 \\
  0
\ar[r]^{}
& H(U)
\ar[r]^={H(i)}
& H(V)
\ar[r]^{}
& HF(\coker(j))
\ar[r]^{}
& 0
}
\]
By assumption, $HF(\operatorname{coker}(j))$ is a trivial cofibrant object in $\calM_\calA$. Since $H(i)$ is a monomorphism, it follows that
$H(i)\in \Cofib(\mathcal A)\cap\Weq(\mathcal A)$.

Now take a relative $F(\mathcal J_\calA)$-cell complex
\[
U_0\to U_1\to \cdots \to U_\alpha\to \cdots,
\]
where each successor map $i_\alpha: U_\alpha\to U_{\alpha+1}$ is a pushout of $F(j_\alpha)$ for some $j_\alpha\in \mathcal J_{\calA}$. Repeating the argument above, we obtain that
$H(i_\alpha)\in \Cofib(\mathcal A)\cap\Weq(\mathcal A)$ for each $\alpha$.

Since $\Cofib(\mathcal A)\cap\Weq(\mathcal A)$ is closed under composition and transfinite composition, we conclude that $H(i)\in \Weq(\mathcal A)$ for each relative $F(\mathcal J_{\mathcal A})$-cell complex $i: U\to V$; see \cite[Proposition 10.3.4]{HPS2003}. Thus, by Lemma \ref{acyclicity condition}, the right acyclicity condition holds.
Hence, the right-induced model structure $\calM_{\calB}$ exists on $\mathcal B$.

Let $\calC_{\calA}$ (resp., $\calW_{\calA}$, $\calF_{\calA}$) denote the class of cofibrant objects (resp., trivial objects, fibrant objects) in $\calM_{\calA}$. It is immediate that the class of fibrant objects in $\calM_{\calB}$ is precisely $H^{-1}(\calF_{\calA})$, and the class of trivial objects in $\calM_{\calB}$ is precisely $H^{-1}(\calW_{\calA})$. We show that every cofibrant object in $\calM_{\calB}$ lies in $H^{-1}(\calC_{\calA})$. To this end, fix a cofibrant object $X$ in $\calM_{\calB}$. Then $0\to X$ is a cofibration, and so it is a retract of a relative $F(\mathcal I_{\calA})$-cell complex, since $\calM_{\calB}$ is a cofibrantly generated model structure with $F(\mathcal I_{\mathcal A})$ the set of generating cofibrations; see Remark \ref{cofibrantly generated}.

We claim that $0\to H(X)$ is a retract of a cofibration in $\calM_\calA$. Let $i\in\mathcal I_{\calA}$. Then $i$ is a monomorphism with
$\operatorname{coker}(i)\in\mathcal C_{\calA}$. As shown above, each pushout of $F(i)$ is sent by $H$ to a monomorphism with cokernel
$HF(\operatorname{coker}i)$, which lies in $\mathcal C_{\calA}$ since
$HF(\mathcal C_{\calA})\subseteq\mathcal C_{\calA}$ by assumption. Therefore, $H$ sends relative
$F(\mathcal I_{\calA})$-cell complexes to cofibrations in $\calM_\calA$, as $\Cofib(\mathcal A)$ is closed under composition and transfinite composition. The claim follows.

Thus, $0\to H(X)$ is in $\Cofib(\mathcal A)$, and so $H(X)\in\calC_{\calA}$.
\end{prf*}

The preceding result ensures the existence of the right-induced model structure on $\mathcal B$. However, it does not guarantee that this induced structure is abelian. To determine when it is abelian, we next study how cotorsion pairs are lifted along Frobenius functors. We begin with the following lemma.

\begin{lem}\label{cot pair (C,D)}
Let $(\calC,\calD)$ be a cotorsion pair in $\calA$. Then we have the equalities
\[
H^{-1}(\calD)=F(\calC)^\perp
\quad\text{and}\quad
H^{-1}(\calC)={^\perp F(\calD)}.
\]
Moreover, both
\[
(^\perp H^{-1}(\calD),\, H^{-1}(\calD))
\quad\text{and}\quad
(H^{-1}(\calC),\, H^{-1}(\calC)^\perp)
\]
are cotorsion pairs in $\calB$.
Furthermore, they are hereditary whenever $(\calC,\calD)$ is hereditary.
\end{lem}

\begin{prf*}
We first show that $H^{-1}(\mathcal D)=F(\calC)^\bot$.
Once this equality is established, it follows immediately that
$
(^\bot H^{-1}(\mathcal D), H^{-1}(\mathcal D)) =
({^\bot (F(\calC)^\bot)}, F(\calC)^\bot)
$
is a cotorsion pair in $\calB$ (see \ref{def of Cot pairs}).
Let $X\in H^{-1}(\mathcal D)$. Then $H(X)\in\calD=\calC^\bot$, and hence,
\[\Ext_{\calB}^1(F(C),X)\cong \Ext_{\calA}^1(C,H(X))=0\]
for all $C\in\calC$ by Lemma \ref{Ext(F,H)}. Thus, $X\in F(\calC)^\bot$. Conversely, take $X\in F(\calC)^\bot$.
Then, again by Lemma \ref{Ext(F,H)},
\[\Ext_{\calA}^1(C,H(X))\cong \Ext_{\calB}^1(F(C),X)=0\]
for all $C\in\calC$. Hence, $H(X)\in\calC^\bot=\calD$, so $X\in H^{-1}(\mathcal D)$.
Thus, $H^{-1}(\mathcal D)=F(\calC)^\bot$, as claimed.

The equality $H^{-1}(\mathcal C)={^\bot F(\calD)}$ is proved similarly.
Consequently,
\[
(H^{-1}(\mathcal C), H^{-1}(\mathcal C)^\bot) =
({^\bot F(\calD)}, ({^\bot F(\calD)})^\bot)
\]
is also a cotorsion pair in $\calB$.

Finally, since $F$ and $H$ are exact, it is readily checked that the resulting cotorsion pairs
$(^\bot H^{-1}(\mathcal D), H^{-1}(\mathcal D))$ and
$(H^{-1}(\mathcal C), H^{-1}(\mathcal C)^\bot)$
are hereditary if $(\calC,\calD)$ is hereditary.
\end{prf*}

\begin{rmk} \label{comolete for (C,C)}
If the cotorsion pair $(\calC,\calD)$ is cogenerated by a set $\calS$, then the same argument as in Lemma \ref{cot pair (C,D)} applies, using $\calS^\bot=\calD$ instead of $\calC^\bot=\calD$. One obtains $H^{-1}(\mathcal D)=F(\calS)^\bot$.
Consequently, the cotorsion pair $(^\bot H^{-1}(\mathcal D), H^{-1}(\mathcal D))$ is cogenerated by the set $F(\calS)$. In this case, \corcite[2.15(3)]{2013Sto} ensures that $(^\bot H^{-1}(\mathcal D), H^{-1}(\mathcal D))$ is complete.
\end{rmk}

\begin{ipg} \label{a key seq}
Suppose that $F : \calA \rightleftarrows \calB : H$ is an adjoint pair.
If $F$ is exact, then it is known that $F$ is faithful if and only if the unit $\eta:\Id_{\calA}\to HF$ is a monomorphism; see, for example, \lemcite[2.1(3)]{2022CR}.
Consequently, for every $X\in\calA$, we have a short exact sequence
\begin{equation*} \label{key ses mono}
\begin{gathered}
0 \to X \xrightarrow{\eta_X} HF(X) \to \coker(\eta_X) \to 0.
\end{gathered}
\end{equation*}
Dually, if $H$ is exact, then $H$ is faithful if and only if the counit $\varepsilon:FH\to\Id_{\calB}$ is an epimorphism.
Hence, for any $Y\in\calB$, there is a short exact sequence
\begin{equation*} \label{key ses} \tag{\ref{a key seq}.1}
\begin{gathered}
0 \to \ker(\varepsilon_Y) \xrightarrow{\lambda_Y} FH(Y)
      \xrightarrow{\varepsilon_Y} Y \to 0.
\end{gathered}
\end{equation*}
These two short exact sequences will be used frequently below.
\end{ipg}

\begin{prp} \label{C_H=D_H}
Let $(\calC,\calD)$ be a cotorsion pair in $\calA$ cogenerated by a set $\calS$.
Suppose that $H$ is faithful.
Then $(H^{-1}(\mathcal C), H^{-1}(\mathcal D))$ is a complete cotorsion pair in $\calB$ cogenerated by a set if and only if the following two conditions are satisfied:
\begin{prt}
\item $HF(\calC)\subseteq\calC$;
\item $\Hom_\calB(\lambda_X,Y)$ is an epimorphism for every $X\in H^{-1}(\mathcal C)$ and every $Y\in H^{-1}(\mathcal D)$.
\end{prt}
\end{prp}

\begin{prf*}
Suppose that $(\calC,\calD)$ is a cotorsion pair in $\calA$ cogenerated by a set $\calS$. By Lemma \ref{cot pair (C,D)}, we have
$
(H^{-1}(\calC), H^{-1}(\calC)^\bot) = ({^\bot F(\calD)}, ({^\bot F(\calD)})^\bot)
$
as a cotorsion pair in $\calB$.
Moreover, using the same lemma along with Remark \ref{comolete for (C,C)}, we obtain that
$
({^\bot H^{-1}(\calD)}, H^{-1}(\calD)) = ({^\bot(F(\calS)^\bot)}, F(\calS)^\bot)
$
is a complete cotorsion pair in $\calB$.

$(\Leftarrow)$
We first prove the inclusion $F(\calS) \subseteq {^\bot F(\calD)}$.
Once this is established, we get
\[
{^\bot H^{-1}(\calD)} = {^\bot(F(\calS)^\bot)}
\subseteq {^\bot(({^\bot F(\calD)})^\bot)}
= {^\bot F(\calD)}
= H^{-1}(\calC).
\]
Indeed, let $S\in\calS$. Since $\calS\subseteq {^\bot(\calS^\bot)}=\calC$, we have $S\in\calC$.
Condition (a) gives $HF(\calC)\subseteq\calC$, hence $HF(S)\in\calC$.
For any $D\in\calD$, Lemma \ref{Ext(F,H)} yields
\[
\Ext_{\calB}^1(F(S),F(D)) \cong \Ext_{\calA}^1(HF(S),D)=0.
\]
Thus, $F(S)\in {^\bot F(\calD)}$, and therefore $F(\calS)\subseteq {^\bot F(\calD)}$.

Next, we prove the inclusion ${^\bot F(\calD)} \subseteq {^\bot H^{-1}(\calD)}$.
Once this is shown, we obtain
$H^{-1}(\calC)={^\bot F(\calD)}\subseteq {^\bot H^{-1}(\calD)}$,
and together with the previous inclusion we have
$H^{-1}(\calC)={^\bot H^{-1}(\calD)}$.
Consequently, $(H^{-1}(\calC),H^{-1}(\calD))$ is a complete cotorsion pair in $\calB$, as desired.

To see this, take $X\in {^\bot F(\calD)}=H^{-1}(\calC)$ and $Y\in H^{-1}(\calD)$.
Note that $FH(Y)\in F(\calD)$ because $H(Y)\in\calD$.
Applying $\Hom_\calB(-,Y)$ to the short exact sequence (\ref{key ses}) gives the exact sequence
\begin{equation*} \label{key ses 2} \tag{$\dagger$}
\begin{gathered}
\xymatrix@C=0.75cm{
   \cdots
\ar[rr]^-{\Hom_\calB(\lambda_X, Y)}
&& \Hom_\calB(\ker(\varepsilon_X), Y)
\ar[r]
& \Ext_\calB^1(X, Y)
\ar[r]
& \Ext_\calB^1(FH(X), Y)
\ar[r]
& \cdots.}
\end{gathered}
\end{equation*}
By condition (b), the map $\Hom_\calB(\lambda_X,Y)$ is an epimorphism.
Moreover, using Lemma \ref{Ext(F,H)}, we have
\[
\Ext_\calB^1(FH(X),Y) \cong \Ext_\calB^1(X,FH(Y))=0.
\]
It then follows from the exact sequence that $\Ext_\calB^1(X,Y)=0$.
Hence, $X\in {^\bot H^{-1}(\calD)}$, which proves the inclusion ${^\bot F(\calD)}\subseteq {^\bot H^{-1}(\calD)}$.

$(\Rightarrow)$
Now suppose that $(H^{-1}(\calC),H^{-1}(\calD))$ is a cotorsion pair in $\calB$.
Then by definition, ${^\bot H^{-1}(\calD)}=H^{-1}(\calC)$.
Consider the exact sequence (\ref{key ses 2}) for any $X\in H^{-1}(\calC)$ and $Y\in H^{-1}(\calD)$.
Since $\Ext_\calB^1(X,Y)=0$, the connecting homomorphism in that sequence implies that $\Hom_\calB(\lambda_X,Y)$ is an epimorphism.
Thus, condition (b) is satisfied.

Moreover, by Lemma \ref{cot pair (C,D)}, we have the equalities
$H^{-1}(\calD)=F(\calC)^\perp$ and $H^{-1}(\calC)={^\bot F(\calD)}$.
Then we get
\[
F(\calC) \subseteq {^\bot(F(\calC)^\bot)} = {^\bot H^{-1}(\calD)} = H^{-1}(\calC) = {^\bot F(\calD)}.
\]
Consequently, for every $C\in\calC$ and $D\in\calD$,
\[
\Ext_\calA^1(HF(C),D) \cong \Ext_\calB^1(F(C),F(D))=0.
\]
Since ${^\bot\calD}=\calC$, this gives the inclusion $HF(\calC)\subseteq\calC$.
\end{prf*}

\begin{rmk}\label{here for cp}
If the given cotorsion pair $(\calC,\calD)$ in $\calA$ is hereditary, then by Lemma \ref{cot pair (C,D)} the induced cotorsion pair $(H^{-1}(\calC),H^{-1}(\calD))$ in $\calB$ is also hereditary.
\end{rmk}

Together with the discussion above, Proposition \ref{C_H=D_H} yields the main lifting theorem of this section. It gives a criterion for lifting a cofibrantly generated Hovey triple from $\mathcal A$ to $\mathcal B$. This result will be applied in the next section to obtain cofibrantly generated Hovey triples in the equivariant category $\calA^G$.

\begin{prp} \label{HOVEY case0}
Let $(\calC,\calW,\calF)$ be a cofibrantly generated Hovey triple in $\calA$.
Suppose that $H$ is faithful.
Then $(H^{-1}(\calC),H^{-1}(\calW),H^{-1}(\calF))$ is a cofibrantly generated Hovey triple in $\calB$ if and only if the following conditions hold:
\begin{prt}
\item $HF(\calC)\subseteq\calC$ and $HF(\calC\cap\calW)\subseteq\calC\cap\calW$;
\item $\Hom_\calB(\lambda_X,Y)$ is an epimorphism for every $X\in H^{-1}(\calC)$ and every $Y\in H^{-1}(\calW\cap\calF)$;
\item $\Hom_\calB(\lambda_X,Y)$ is an epimorphism for every $X\in H^{-1}(\calC\cap\calW)$ and every $Y\in H^{-1}(\calF)$.
\end{prt}
\end{prp}

\begin{prf*}
Since $(\calC,\calW,\calF)$ is a cofibrantly generated Hovey triple in $\calA$, both $(\calC,\calW\cap\calF)$ and $(\calC\cap\calW,\calF)$ are complete cotorsion pairs in $\calA$ cogenerated by sets.

$(\Leftarrow)$
Applying Proposition \ref{C_H=D_H} gives that
\[
(H^{-1}(\calC),\, H^{-1}(\calW)\cap H^{-1}(\calF))
\quad\text{and}\quad
(H^{-1}(\calC)\cap H^{-1}(\calW),\, H^{-1}(\calF))
\]
are complete cotorsion pairs in $\calB$ cogenerated by sets.
Moreover, since $\calW$ is thick and $H$ is exact, it follows that $H^{-1}(\calW)$ is also thick.
Thus, $(H^{-1}(\calC),H^{-1}(\calW),H^{-1}(\calF))$ is a cofibrantly generated Hovey triple in $\calB$.

$(\Rightarrow)$
Suppose that $(H^{-1}(\calC),H^{-1}(\calW),H^{-1}(\calF))$ is a cofibrantly generated Hovey triple in $\calB$.
Then
$(H^{-1}(\calC), H^{-1}(\calW)\cap H^{-1}(\calF))$
and
$(H^{-1}(\calC)\cap H^{-1}(\calW), H^{-1}(\calF))$
are complete cotorsion pairs in $\calB$ cogenerated by sets.
Observe that
$H^{-1}(\calW)\cap H^{-1}(\calF)=H^{-1}(\calW\cap\calF)$
and
$H^{-1}(\calC)\cap H^{-1}(\calW)=H^{-1}(\calC\cap\calW)$.
Consequently,
\[
(H^{-1}(\calC), H^{-1}(\calW\cap\calF))
\quad\text{and}\quad
(H^{-1}(\calC\cap\calW), H^{-1}(\calF))
\]
are complete cotorsion pairs in $\calB$ cogenerated by sets.
Applying Proposition \ref{C_H=D_H} yields conditions (a), (b), and (c) in the statement.
\end{prf*}

In the following, we show that condition (c) in Proposition \ref{HOVEY case0} can be removed.

\begin{rmk}\label{rmk2.13}
Suppose that $(\calC,\calW,\calF)$ is a cofibrantly generated Hovey triple in $\calA$. Then the corresponding abelian model structure
\[
\calM_\calA=(\Mon(\calC),\sfW,\Epi(\calF))
\]
is cofibrantly generated, where
\[
\sfW=\{ w \mid w=fc \text{ with } c\in\Mon(\calC\cap\calW),\ f\in\Epi(\calW\cap\calF)\}.
\]
Furthermore, if $HF(\calC)\subseteq\calC$ and $HF(\calC\cap\calW)\subseteq\calC\cap\calW$, then $HF$ preserves cofibrant and trivial cofibrant objects in $\calM_\calA$. It follows from Proposition \ref{HF condition} that there is a right-induced model structure
\[
\calM_{\calB}=(^{\boxslash}H^{-1}(\sfW\cap\Epi(\calF)),\, H^{-1}(\sfW),\, H^{-1}(\Epi(\calF)))
\]
on $\mathcal B$. Let $\calC_{\calB}$, $\calW_{\calB}$, and $\calF_{\calB}$ denote the classes of cofibrant, trivial, and fibrant objects in $\calM_{\calB}$, respectively. Then $\calW_{\calB}=H^{-1}(\calW)$ and $\calF_{\calB}=H^{-1}(\calF)$, while $\calC_{\calB}\subseteq H^{-1}(\calC)$ by Proposition \ref{HF condition}.
\end{rmk}

\begin{lem}\label{bc condition}
Adopt the setup and notation from Proposition \ref{HOVEY case0} and Remark \ref{rmk2.13}. Suppose that condition (a) holds. Then condition (b) is equivalent to $\calC_{\calB}=H^{-1}(\calC)$, which in turn implies condition (c).
\end{lem}

\begin{prf*}
We know that $\calC_{\calB}\subseteq H^{-1}(\calC)$; see Remark \ref{rmk2.13}. We first give two claims.

\vspace{1ex}\textbf{Claim 1}. Trivial fibrations $p$ in $\calM_\calB$ are precisely the epimorphisms such that $\ker(p)\in H^{-1}(\mathcal W\cap\mathcal F)$.

Indeed, $p$ is a trivial fibration if and only if $H(p)$ is a trivial fibration in $\calM_\calA$. Equivalently, $H(p)$ is an epimorphism whose kernel lies in $\calW\cap\calF$. Since $H$ is faithful and exact, it reflects epimorphisms and satisfies $H(\ker p)\cong \ker H(p)$. This proves the claim.

\vspace{1ex}\textbf{Claim 2}.
For every $X\in\mathcal B$, $X\in\calC_{\calB}$ if and only if $\Ext^1_{\calB}(X,Y)=0$ for all $Y\in H^{-1}(\mathcal W\cap\mathcal F)$.

Suppose first that $X\in\calC_{\calB}$. Consider the exact sequence
$0\to Y\to E\xrightarrow{p} X\to 0$
with $Y\in H^{-1}(\mathcal W\cap\mathcal F)$.
By Claim 1, $p$ is a trivial fibration in $\calM_\calB$. Since $0\to X$ is a cofibration, the diagram
\[
\xymatrix@C=1.25cm{
  0
\ar[r]
\ar[d]
& E
\ar[d] \\
  X
\ar[r]^-{\id_X}
\ar@{-->}[ur]
& X
}
\]
admits a lift, which gives a section of $E\to X$. Hence, the above short exact sequence splits, and so $\Ext^1_{\mathcal B}(X,Y)=0$.
Conversely, suppose that
$\Ext^1_{\mathcal B}(X,Y)=0$
for every $Y\in H^{-1}(\mathcal W\cap\mathcal F)$. Let $p:E\to Z$ be a trivial fibration in $\calM_\calB$. By Claim 1, $p$ is an epimorphism with
$\ker(p)\in H^{-1}(\mathcal W\cap\mathcal F)$. For any morphism $f: X\to Z$ in $\calB$, consider the pullback
\[
\xymatrix@C=1.25cm{
P
\ar[r]
\ar[d]
& E
\ar[d]^p \\
X
\ar[r]^-{f}
& Z .
}
\]
Then we have a short exact sequence
\[
0\to \ker(p)\to P\to X\to 0.
\]
By the assumed Ext-vanishing, this sequence splits. A section $X\to P$, followed by the morphism $P\to E$, gives a lift $X\to E$ in the diagram
\[
\xymatrix@C=1.25cm{
  0
\ar[r]
\ar[d]
& E
\ar[d]^-{p} \\
  X
\ar[r]^-{f}
\ar@{-->}[ur]
& Z .
}
\]
Hence, $0\to X$ is a cofibration in $\calM_\calB$, that is, $X\in\calC_{\calB}$. This proves Claim 2.
\vspace{1ex}

Now take $X\in H^{-1}(\mathcal C)$ and $Y\in H^{-1}(\mathcal W\cap\mathcal F)$.
Since $H$ is faithful, the counit gives a short exact sequence
$0\to \ker(\varepsilon_X) \xrightarrow{\lambda_X} FH(X) \xrightarrow{\varepsilon_X} X\to 0$.
Applying $\Hom_{\mathcal B}(-,Y)$, we obtain an exact sequence
\[
\Hom_{\mathcal B}(FH(X),Y)
\xrightarrow{\Hom_{\mathcal B}(\lambda_X,Y)}
\Hom_{\mathcal B}(\ker(\varepsilon_X),Y)
\to
\Ext^1_{\mathcal B}(X,Y)
\to
\Ext^1_{\mathcal B}(FH(X),Y).
\]
By Lemma \ref{Ext(F,H)}, we have
$\Ext^1_{\mathcal B}(FH(X),Y) \cong \Ext^1_{\mathcal A}(H(X),H(Y))=0$,
since $H(X)\in\mathcal C$ and $H(Y)\in\mathcal W\cap\mathcal F$.
It follows that $\Hom_{\mathcal B}(\lambda_X,Y)$ is an epimorphism if and only if
$\Ext^1_{\mathcal B}(X,Y)=0$.
Thus, by Claim 2, condition (b) is equivalent to $\calC_{\calB}=H^{-1}(\mathcal C)$.

Finally, we prove that the equality $\calC_{\calB}=H^{-1}(\mathcal C)$ implies condition (c). To this end, let
$X\in H^{-1}(\mathcal C\cap\mathcal W)$ and $Y\in H^{-1}(\mathcal F)$.
Then $X\in H^{-1}(\mathcal C)=\calC_{\calB}$.
Moreover, since $H(X)\in\mathcal W$, the morphism $0\to X$ is a weak equivalence in $\calM_\calB$. Hence, $0\to X$ is a trivial cofibration in $\calM_\calB$.

Let $0\to Y\to E\to X\to 0$ be a short exact sequence in $\calB$.
Since $Y\in H^{-1}(\mathcal F)$, the epimorphism $E\to X$ is a fibration in $\calM_\calB$. Thus, there is a morphism $s:X\to E$ making the diagram
\[
\xymatrix@C=1.25cm{
  0
\ar[r]
\ar[d]
& E
\ar[d]^{q} \\
  X
\ar[r]^-{\id_X}
\ar@{-->}[ur]^{s}
& X
}
\]
commute. Hence, $s$ is a section of $q:E\to X$, so the above short exact sequence splits. This yields $\Ext^1_{\mathcal B}(X,Y)=0$.
Applying $\Hom_{\mathcal B}(-,Y)$ to the exact sequence
$$0\to \ker(\varepsilon_X) \xrightarrow{\lambda_X} FH(X) \xrightarrow{\varepsilon_X} X\to 0,$$
we get that $\Hom_{\mathcal B}(\lambda_X,Y)$ is an epimorphism. Thus, condition (c) holds.
\end{prf*}

Now the following result is immediate from Proposition \ref{HOVEY case0} and Lemma \ref{bc condition}.

\begin{thm} \label{HOVEY case1}
Let $(\calC,\calW,\calF)$ be a cofibrantly generated Hovey triple in $\calA$. Suppose that $H$ is faithful. Then the following are equivalent.
\begin{eqc}
\item
$\bigl(H^{-1}(\calC),H^{-1}(\calW),H^{-1}(\calF)\bigr)$
is a cofibrantly generated Hovey triple in $\calB$.

\item The following two conditions hold:
\begin{prt}
\item
$HF(\calC)\subseteq \calC$
and
$HF(\calC\cap\calW)\subseteq \calC\cap\calW;$
\item
$\Hom_{\calB}(\lambda_X,Y)$
is an epimorphism for every
$X\in H^{-1}(\calC)$
and
$Y\in H^{-1}(\calW\cap\calF).$
\end{prt}
\end{eqc}
Moreover, in this case, the corresponding abelian model structure $\calM_{\calB}$ on $\calB$ is the right-induced model structure along $H$.
\end{thm}

\begin{rmk} \label{compare with DP}
Theorem \ref{HOVEY case1} is closely related to the right-lifting theorem for hereditary abelian model structures established by Dalezios and Psaroudakis \thmcite[5.12]{2023DP}.
Their result is formulated for a general adjunction and gives sufficient conditions for the existence of a right-lifted abelian model structure, one of the essential hypotheses being a right acyclicity condition.
In contrast, under the additional assumption that $(F, H)$ is a Frobenius pair and that $H$ is faithful, Theorem \ref{HOVEY case1} gives an if-and-only-if criterion for the right-induced model structure to be represented by the Hovey triple
\[
(H^{-1}(\calC), H^{-1}(\calW), H^{-1}(\calF)).
\]
More precisely, the Frobenius structure allows the right acyclicity condition to be reformulated in terms of the preservation properties of $H F$ together with the epimorphism condition in Theorem \ref{HOVEY case1}(ii)(b).
Thus, the two results are complementary: \thmcite[5.12]{2023DP} applies to general adjunctions, whereas Theorem \ref{HOVEY case1} provides a more explicit characterization in the Frobenius setting.
\end{rmk}

\begin{rmk}\label{here for ht}
If the given Hovey triple $(\calC,\calW,\calF)$ in $\calA$ is hereditary, then by Remark \ref{here for cp} the induced Hovey triple $(H^{-1}(\calC),H^{-1}(\calW),H^{-1}(\calF))$ in $\calB$ is hereditary as well.
\end{rmk}

Let $\calM_\calA = (\calC_\calA,\calW_\calA,\calF_\calA)$ and $\calM_\calB = (\calC_\calB,\calW_\calB,\calF_\calB)$ be Hovey triples in $\calA$ and $\calB$, respectively.
The following result shows that, under suitable conditions, a Quillen adjunction between $\calM_\calA$ and $\calM_\calB$ can be promoted to a Quillen equivalence.
Note that the following result does not require $(F,H)$ to be a Frobenius pair.

\begin{thm} \label{Ho(M_A)=Ho(M_B)}
Let $(F, H)$ be a Quillen adjunction between $\calM_\calA$ and $\calM_\calB$.
Suppose that both $F$ and $H$ are exact and faithful and that $H(\calC_\calB \cap \calW_\calB) \subseteq \calW_\calA$ and $F(\calW_\calA \cap \calF_\calA) \subseteq \calW_\calB$.
Then the following are equivalent.
\begin{eqc}
\item
$(F,H)$ is a Quillen equivalence between $\calM_\calA$ and $\calM_\calB$.
\item
$\coker(\eta_X) \in \calW_\calA$ for every $X\in\calC_\calA$, and
$\ker(\varepsilon_Y) \in \calW_\calB$ for every $Y\in\calF_\calB$.
\end{eqc}
\end{thm}

\begin{prf*}
It is known that a Quillen adjunction $(F,H)$ is a Quillen equivalence if and only if for every $X\in\calC_\calA$ and $Y\in\calF_\calB$, the compositions
\[
X \xrightarrow{\eta_X} HF(X) \xrightarrow{H(p)} HRF(X)
\quad\text{and}\quad
FQH(Y) \xrightarrow{F(q)} FH(Y) \xrightarrow{\varepsilon_Y} Y
\]
are weak equivalences; see \thmcite[7.17]{2025GillBook} or \prpcite[1.3.13]{Hovey1999}.
Here, $RF(X)$ is a fibrant replacement of $F(X)$ in $\calM_\calB$ with $p: F(X) \to RF(X)$ a trivial cofibration, and $QH(Y)$ is a cofibrant replacement of $H(Y)$ in $\calM_\calA$ with $q: QH(Y) \to H(Y)$ a trivial fibration.

We first show that both $H(p)$ and $F(q)$ are weak equivalences.
We prove the statement for $H(p)$; the argument for $F(q)$ is analogous.
Since $p$ is a trivial cofibration, there is a short exact sequence
\[
0 \to F(X) \xrightarrow{p} RF(X) \to \coker(p) \to 0,
\]
with $\coker(p) \in \calC_\calB \cap \calW_\calB$ (see \prpcite[4.2]{HOVEY2002}).
Applying the exact functor $H$ yields
\[
0 \to HF(X) \xrightarrow{H(p)} HRF(X) \to H(\coker(p)) \to 0.
\]
By hypothesis, $H(\coker(p)) \in \calW_\calA$, and \lemcite[5.8]{HOVEY2002} then implies that $H(p)$ is a weak equivalence.
Consequently, to finish the proof, it remains to establish the following two equivalences:
\begin{rqm}
\item[$\bullet$]
$\eta_X$ is a weak equivalence if and only if $\coker(\eta_X) \in \calW_\calA$;

\item[$\bullet$]
$\varepsilon_Y$ is a weak equivalence if and only if $\ker(\varepsilon_Y) \in \calW_\calB$.
\end{rqm}
Both follow directly from \lemcite[5.8]{HOVEY2002}.
Indeed, $\eta_X$ is a monomorphism by faithfulness of $F$, and $\varepsilon_Y$ is an epimorphism by faithfulness of $H$.
\end{prf*}

\section{Proof of Theorem \ref{retract (Q^G,W^G,R^G)}}
\label{Proofs of Theorems A}
\noindent
In this section, we prove Theorem \ref{retract (Q^G,W^G,R^G)} stated in the introduction. {\sf Throughout this section}, let $\calA$ be a Grothendieck category with enough projectives and $G$ a finite group acting on $\calA$.

\subsection{Proof of Theorem \ref{retract (Q^G,W^G,R^G)}(1)} \label{Thm A(1)}
This subsection is devoted to giving a proof of Theorem \ref{retract (Q^G,W^G,R^G)}(1).
From \ref{frobenius pair (Ind,w)}, one sees that
\[
\Ind : \calA \rightleftarrows \calA^G : U
\]
forms a Frobenius pair with $U$ faithful. We denote by $\eta: \Id_{\calA} \to U \circ \Ind$ the unit and by $\varepsilon: \Ind \circ U \to \Id_{\calA^G}$ the counit.

Recall that an integer $n$ is said to be \emph{invertible} in an additive category if for every morphism $f:X\to Y$ in it, there exists a unique morphism $g:X\to Y$ such that $f=ng$.
For a finite group $G$, denote by $|G|$ the order of $G$. The following result may be found in \cite[Lemma 2.14]{2019Sun}. Note that every Grothendieck category is automatically idempotent complete.

\begin{lem} \label{split cond}
Suppose that $|G|$ is invertible in $\calA$.
Then for each object $X$ in $\calA^G$, the short exact sequence
$
0  \to  \kernel(\varepsilon_X)  \overset{\lambda_X} \longrightarrow      (\Ind \circ U)(X)
                                \overset{\varepsilon_X} \longrightarrow  X  \to  0
$
is split.
\end{lem}

In what follows, for any class $\calE$ of objects in $\calA$, we set
\[
\calE^G = \{(X, \phi) \in \calA^G \,|\, U(X, \phi) = X \in \calE \}.
\]

\begin{bfhpg}[\bf Proof of Theorem \ref{retract (Q^G,W^G,R^G)}(1)]
We will use Theorem \ref{HOVEY case1} to show that $\calM^G = (\calC^G, \calW^G, \calF^G)$ is a cofibrantly generated Hovey triple in $\calA^G$.
It suffices to verify the following two conditions:
\begin{rqm}
\item
$(U \circ \Ind)(\calC) \subseteq \calC$ and $(U \circ \Ind)(\calC \cap \calW) \subseteq \calC \cap \calW$;

\item
$\Hom_{\calA^G}(\lambda_X, Y)$ is an epimorphism for every $X \in \calC^G$ and every $Y \in (\calW \cap \calF)^G$,
\end{rqm}
where $\lambda_X$ appears in the short exact sequence
in Lemma \ref{split cond}.
Indeed, condition (2) is guaranteed by Lemma \ref{split cond}.
For condition (1), we only show that $(U \circ \Ind)(\calC) \subseteq \calC$;
the other inclusion is proved similarly.
To this end, take $X \in \calC$.
Write $(\rho, \theta)$ for the action of $G$ on $\calA$.
Then $(U \circ \Ind)(X) = \bigoplus_{g \in G} \rho_g X$;
see \ref{frobenius pair (Ind,w)}.
By Lemma \ref{(Q,W,R) G-invariant}, the class $\calC$ is $G$-invariant.
Thus, $\rho_g X \in \calC$ for every $g \in G$.
Since $\calC$ is the left part of the cotorsion pair $(\calC, \calW \cap \calF)$, it is closed under direct sums.
It follows that $(U \circ \Ind)(X) = \bigoplus_{g \in G} \rho_g X \in \calC$.
Hence, $(U \circ \Ind)(\calC) \subseteq \calC$, as desired.

By Remark \ref{here for ht}, if the given Hovey triple
$(\calC, \calW, \calF)$ in $\calA$ is hereditary, then the induced Hovey triple $(\calC^G, \calW^G, \calF^G)$ in $\calA^G$ is also hereditary.
\qed
\end{bfhpg}

\subsection{Proof of Theorem \ref{retract (Q^G,W^G,R^G)}(2)} \label{Thm Ho(M)G}
\noindent
In this subsection, we provide a proof of Theorem \ref{retract (Q^G,W^G,R^G)}(2).

\begin{setup}
Throughout this subsection, let $\calM = (\calC, \calW, \calF)$ be a hereditary Hovey triple in $\calA$ with $\gamma_{\calA}: \calA \to \Ho(\calM)$ the localization functor.
\end{setup}

We first extend the given $G$-action on $\calA$ to the homotopy category $\Ho(\calM)$ and to the stable category $\St_\omega(\calC\cap\calF)$, where $\omega = \calC\cap\calW\cap\calF$.
Using the triangle equivalence
\[
\iota : \St_\omega(\calC \cap \calF) \longrightarrow \Ho(\calM)
\]
from Theorem \ref{hered Hovey correspondence}, we then obtain an induced triangle equivalence between the corresponding equivariant categories.
As a consequence, the canonical triangulated structure of the stable category transfers to the homotopy category.

We will need the following result, taken from \corcite[5.16]{2025GillBook}. It provides a useful tool for constructing natural isomorphisms between functors defined on the homotopy category, and will be used repeatedly in the remainder of the paper.

\begin{lem} \label{natural iso}
Suppose that $K : \Ho(\calM) \to \calB$ and $L : \Ho(\calM) \to \calB$ are functors and $\eta = \{ \eta_D : K(D) \to L(D) \}_{D \in \calA}$ is a family of morphisms in $\calB$.
Then $\eta : K \to L$ is a natural transformation if and only if $\eta : K \circ \gamma_{\calA} \to L \circ \gamma_{\calA}$ is a natural transformation.
\end{lem}

\begin{bfhpg}[\bf Admissible $G$-actions on categories]\label{admiss GrAct on Cat}
Let $\calT$ be a triangulated category with the shift functor $\mathsf{\Sigma}$, and let $G$ be a finite group.
Recall from \dfncite[3.1]{2019Sun} that a $G$-action $(\rho,\theta)$ on $\calT$ is called \emph{admissible} if it is equipped with a family of natural isomorphisms
\[
\varphi = \{\varphi_g : \mathsf{\Sigma} \circ \rho_g \overset{\cong}{\longrightarrow} \rho_g \circ \mathsf{\Sigma}\}_{g \in G}
\]
such that each $(\rho_g,\varphi_g)$ is a triangle functor and the pair $(\mathsf{\Sigma},\varphi)$ is a $G$-functor.
\end{bfhpg}

As noted in Remark \ref{tri for Ho(M)}, the homotopy category $\Ho(\calM)$ is triangulated. The following lemma shows that, under mild assumptions, a $G$-action on $\calA$ lifts to an admissible $G$-action on $\Ho(\calM)$; its proof is deferred to the appendix (see \ref{detailed proof}).

\begin{lem} \label{admiss for Ho(M)}
If any two of the classes $\calC$, $\calW$, $\calF$ are $G$-invariant, then the $G$-action $(\rho,\theta)$ on $\calA$ induces an admissible $G$-action $(\widehat{\rho},\widehat{\theta})$ on $\Ho(\calM)$.
\end{lem}

By Happel's classical construction, the stable category $\St_\omega(\calC\cap\calF)$ of the Frobenius category $\calC\cap\calF$ modulo $\omega=\calC\cap\calW\cap\calF$ is triangulated. The next result gives an analogous lifting statement for the stable category; again, its proof appears in the appendix (see \ref{pf of lem 3.7}).

\begin{lem} \label{admiss for St(C cap F)}
If any two of the classes $\calC$, $\calW$, $\calF$ are $G$-invariant, then the $G$-action $(\rho,\theta)$ on $\calA$ induces an admissible $G$-action $(\widetilde{\rho},\widetilde{\theta})$ on $\St_\omega(\calC\cap\calF)$.
\end{lem}

Recall from Karakikes, Kontogeorgis, and Psaroudakis \dfncite[2.21]{2025Equ} that given an admissible $G$-action on a triangulated category $\calT$, a \emph{canonical triangulated structure} on $\calT^G$ is such that the forgetful functor is a triangulated functor. It is well known that $\calT^G$ is not always a triangulated category.

Based on Chen's work \cite{2015Chen}, Sun gave in \prpcite[3.3]{2019Sun} the following result, which shows that under certain conditions, $\calT^G$ always has a unique canonical pre-triangulated structure.

\begin{lem} \label{G-pre-tri stru}
Let $\calT$ be a triangulated category and $G$ a finite group admissibly acting on $\calT$ such that $|G|$ is invertible in $\calT$.
Then $\calT^G$ admits a unique canonical pre-triangulated structure.
Moreover, a triangle in $\calT^G$ is exact precisely when its image under the forgetful functor is exact in $\calT$.
\end{lem}

\begin{dfn} \cite{2025Equ} \label{tri G functor}
Let $\calT$ and $\widehat{\calT}$ be triangulated categories with shift functors $\mathsf{\Sigma}$ and $\widehat{\mathsf{\Sigma}}$, respectively, endowed with admissible
$G$-actions $(\rho,\theta,\varphi)$ and $(\widehat{\rho},\widehat{\theta},\widehat{\varphi})$.
A $G$-functor
\[
(F,\sigma):(\calT,\rho,\theta)
\longrightarrow
(\widehat{\calT},\widehat{\rho},\widehat{\theta})
\]
is called a \emph{triangulated $G$-functor} if $F$ is a triangle functor with structure isomorphism $\eta : F \mathsf{\Sigma} \overset{\cong}{\longrightarrow} \widehat{\mathsf{\Sigma}} F$, and $\eta$ is $G$-natural. Equivalently, for every $g\in G$, the following diagram is commutative:
\[
\xymatrix{
F\mathsf{\Sigma}\rho_g
\ar[r]^-{\eta\rho_g}
\ar[d]_-{(\sigma_g\mathsf{\Sigma})\circ(F\varphi_g)}
&
\widehat{\mathsf{\Sigma}}F\rho_g
\ar[d]^-{(\widehat{\varphi}_gF)\circ
(\widehat{\mathsf{\Sigma}}\sigma_g)}
\\
\widehat{\rho}_gF\mathsf{\Sigma}
\ar[r]_-{\widehat{\rho}_g\eta}
&
\widehat{\rho}_g\widehat{\mathsf{\Sigma}}F
}
\]
\end{dfn}

The following result, taken from \lemcite[2.23]{2025Equ}, may be viewed as the triangulated analogue of \ref{equi fun}.

\begin{lem} \label{G-tri fun}
Let $\calT$ and $\widehat{\calT}$ be triangulated categories, and let $G$ be a finite group acting admissibly on both categories such that $|G|$ is invertible in each.
Suppose that $F : \calT \to \widehat{\calT}$ is a triangulated $G$-functor.
If $\calT^G$ and $\widehat{\calT}^G$ are equipped with the canonical triangulated structure, then the induced equivariant functor $F^G : \calT^G \to \widehat{\calT}^G$ is triangulated.
\end{lem}

Note that the proof of Lemma \ref{G-tri fun} does not depend on the octahedral axiom. Hence, the same conclusion remains valid when $\calT^G$ and $\widehat{\calT}^G$ are only equipped with their canonical pre-triangulated structures. Combining with Lemma \ref{G-pre-tri stru} yields the following result.

\begin{lem} \label{G-pretri fun}
Let $\calT$ and $\widehat{\calT}$ be triangulated categories, and let $G$ be a finite group acting admissibly on both categories such that $|G|$ is invertible in each.
Suppose that $F : \calT \to \widehat{\calT}$ is a triangulated $G$-functor.
Then the induced equivariant functor $F^G : \calT^G \to \widehat{\calT}^G$ is triangulated, where both $\calT^G$ and $\widehat{\calT}^G$ are equipped with their unique canonical pre-triangulated structure.
\end{lem}

We are now ready to prove Theorem \ref{retract (Q^G,W^G,R^G)}(2).

\begin{bfhpg}[\bf Proof of Theorem \ref{retract (Q^G,W^G,R^G)}(2)]\label{3.11}
By Theorem \ref{hered Hovey correspondence}, we have the triangle equivalence
\[
\iota : \St_\omega(\calC \cap \calF) \to \Ho(\calM).
\]
Moreover, $\iota$ is a $G$-functor with the natural $G$-structure $\sigma^\iota_g = \Id$.
By Lemmas \ref{admiss for Ho(M)} and \ref{admiss for St(C cap F)}, the $G$-actions on $\St_{\omega}(\calC \cap \calF)$ and $\Ho(\calM)$ are admissible.
Choosing the same defining short exact sequences for the shift functors, their admissible structures satisfy
\[
\iota \circ \widetilde{\varphi}_g = \widehat{\varphi}_g \circ \iota
\]
for every $g\in G$.
Hence the structure isomorphism $\iota \circ \mathsf{\Sigma} \cong \mathsf{\Sigma} \circ \iota$ is $G$-natural, and therefore $\iota$ is a triangulated $G$-functor.

Note that $\calC \cap \calF$ is a full subcategory of $\calA$. The order $|G|$ is invertible in that category, and hence, also in $\St_\omega(\calC \cap \calF)$.
Furthermore, every morphism in $\Ho(\calM)$ is represented by a morphism in $\calA$ via the bifibrant replacement functor, so $|G|$ is invertible in $\Ho(\calM)$ as well.
Thus, by Lemma \ref{G-pretri fun}, the induced equivariant functor
\begin{equation*}\tag{\ref{3.11}.1}
\iota^G : \St_\omega(\calC \cap \calF)^G \to \Ho(\calM)^G
\end{equation*}
is triangulated, where both are equipped with their unique canonical pre-triangulated structures. Lemma \ref{G-equivalence} tells us that $\iota^G$ is an equivalence, hence, a triangle equivalence.

By Theorem \ref{retract (Q^G,W^G,R^G)}(1), $\calM^G = (\calC^G, \calW^G, \calF^G)$ is a hereditary Hovey triple in $\calA^G$.
Hence, by Proposition \ref{C cap F Frobenius}, $(\calC \cap \calF)^G = \calC^G \cap \calF^G$ is a Frobenius category.
Together with the admissible $G$-action on $\calC \cap \calF$ and the invertibility of $|G|$ noted above, the assumptions of \thmcite[3.14]{2019Sun} are satisfied.
Therefore, the canonical pre-triangulated structure on $\St_\omega(\calC\cap\calF)^G$ is in fact triangulated.
Consequently, the canonical pre-triangulated structure on $\Ho(\calM)^G$ is also triangulated.
\qed
\end{bfhpg}

\subsection{Proof of Theorem \ref{retract (Q^G,W^G,R^G)}(3)} \label{Thm A(2)}
In this subsection, we prove Theorem \ref{retract (Q^G,W^G,R^G)}(3).

\begin{setup}
In this subsection, suppose that $|G|$ is invertible in $\calA$, let $\calM = (\calC, \calW, \calF)$ be a cofibrantly generated hereditary Hovey triple in $\calA$ such that any two of the classes $\calC$, $\calW$, $\calF$ are $G$-invariant, and denote by $\gamma_{\calA} : \calA \to \Ho(\calM)$ the localization functor.
\end{setup}

\begin{bfhpg}[\bf The induced functor $\Ho(\gamma_{\calA}^G)$] \label{exis for Ho^G}
By Theorem \ref{retract (Q^G,W^G,R^G)}(1), $\calM^G = (\calC^G, \calW^G, \calF^G)$ is a cofibrantly generated Hovey triple in $\calA^G$.
Hence, we may consider the localization functor
\[\gamma_{\calA^G} : \calA^G \to \Ho(\calM^G).\]
On the other hand, Remark \ref{gamma is G fun} gives the equivariant functor \[\gamma_{\calA}^G : \calA^G \to \Ho(\calM)^G\] between the corresponding equivariant categories.

Note that a morphism $f$ in $\calA^G$ is a weak equivalence if and only if its image $U_\calA(f)$ under the forgetful functor $U_\calA : \calA^G \to \calA$ is a weak equivalence in $\mathcal{A}$.
Moreover, we have $U_{\Ho(\calM)} \circ \gamma_{\calA}^G = \gamma_{\calA} \circ U_\calA$.
It follows that $\gamma_{\calA}^G$ sends weak equivalences in $\calA^G$ to isomorphisms.
Hence, by the Localization Theorem \ref{loc thm}, there is a unique functor
\[
\Ho(\gamma_{\calA}^G) : \Ho(\calM^G) \to \Ho(\calM)^G
\]
making the following diagram commute:
\[
\xymatrix@R=0.75CM{
   \calA^G
\ar[d]_-{\gamma_{\calA}^G}
\ar[r]^-{\gamma_{\calA^G} \quad}
&  \Ho(\calM^G)
\ar@{-->}[dl]^(0.4){\Ho(\gamma_{\calA}^G)}    \\
   \Ho(\calM)^G                      }
\]
\end{bfhpg}

Recall that a functor $F:\calC\to\calD$ is called \emph{dense up to retracts} if every object $D \in \calD$ is a retract of $F(C)$ for some object $C\in\calC$; that is, there exist morphisms $u:D\to F(C)$ and $v:F(C)\to D$
such that $v\circ u=\Id_D$. The functor $F$ is called an \emph{equivalence up to retracts} if it is fully faithful and dense up to retracts;
see \cite[Lemma 3.4]{2015Chen} for an equivalent formulation.

Our next goal is to show that $\Ho(\gamma_{\calA}^G)$ is a triangle equivalence up to retracts. To this end, we first construct a functor
\[
\Theta:\Ho(\calM^G)\to\Ho(\calM)^G
\]
which is a triangle equivalence up to retracts, and then prove that $\Theta$ is naturally isomorphic to $\Ho(\gamma_{\calA}^G)$.

\begin{bfhpg}[\bf Construction of the functor $\Theta$] \label{con for Theta}
The $G$-action on $\calA$ restricts to a $G$-action on the full subcategory $\calC \cap \calF$.
By Lemma \ref{admiss for St(C cap F)}, this induces an admissible $G$-action on the stable category $\St_\omega(\calC \cap \calF)$.
The canonical functor $K: \calC \cap \calF \to \St_\omega(\calC \cap \calF)$ is a $G$-functor, and therefore, gives rise to an equivariant functor
$K^G: (\calC \cap \calF)^G \to \St_\omega(\calC \cap \calF)^G$.
We then obtain a triangulated functor
\[
\St(K^G) : \St_{\omega^G}(\calC^G \cap \calF^G) \to \St_\omega(\calC \cap \calF)^G,
\]
which serves as the comparison functor between the stable categories, as described in \thmcite[3.14]{2019Sun}.
It is defined as follows:
\begin{rqm}
\item[$\bullet$]
for each $(X, \{\phi_g\}) \in \St_{\omega^G}(\calC^G \cap \calF^G)$, set
$
\St(K^G)(X, \{\phi_g\}) = K^G(X, \{\phi_g\}) = (X, \{ [\phi_g] \});
$
\item[$\bullet$]
for each
$
[f] \in \Hom_{\St_{\omega^G}(\calC^G \cap \calF^G)}((X, \{\phi_g\}), (Y, \{\psi_g\})),
$
the representative $f : (X, \{\phi_g\}) \to (Y, \{\psi_g\})$ is a morphism in $\calC^G \cap \calF^G$.
Then define
\[
\St(K^G)([f]) = K^G(f) = [f] \in \Hom_{\St_\omega(\calC \cap \calF)^G}((X, \{[\phi_g]\}), (Y, \{[\psi_g]\})).
\]
\end{rqm}

By Theorem \ref{retract (Q^G,W^G,R^G)}(1), $\calM^G=(\calC^G,\calW^G,\calF^G)$ is a cofibrantly generated hereditary Hovey triple in $\calA^G$. Hence, Theorem \ref{hered Hovey correspondence} yields the triangle equivalence
\[
RQ:\Ho(\calM^G) \longrightarrow \St_{\omega^G}(\calC^G\cap\calF^G),
\]
given by the bifibrant replacement functor. In particular, for an object $(X,\phi)\in\Ho(\calM^G)$, we have $RQ(X,\phi)=(RQX,\lambda^{QX})$, where $(RQX,\lambda^{QX})$ is taken from
\[
(RQX, \lambda^{QX}) \leftarrowtail  (QX, \tau^X) \twoheadrightarrow  (X, \phi).
\]
On the other hand, Theorem \ref{retract (Q^G,W^G,R^G)}(2) provides the triangle equivalence
\[
\iota^G:\St_\omega(\calC\cap\calF)^G \longrightarrow \Ho(\calM)^G;
\]
see (\ref{3.11}.1).
We now define the comparison functor $\Theta$ as the composition
\[
\Theta:
\Ho(\calM^G)
\xra{RQ}
\St_{\omega^G}(\calC^G\cap\calF^G)
\xra{\St(K^G)}
\St_\omega(\calC\cap\calF)^G
\xra{\iota^G}
\Ho(\calM)^G.
\]
According to \thmcite[3.14]{2019Sun}, the comparison functor $\St(K^G)$ is a triangle equivalence up to retracts, and becomes a genuine triangle equivalence whenever $\St_{\omega^G}(\calC^G\cap\calF^G)$ is idempotent complete. Consequently, $\Theta$ is also a triangle equivalence up to retracts, and is a triangle equivalence under the same idempotent completeness condition.
\end{bfhpg}

\begin{rmk}
In the derived category setting, Chen \cite[Proposition 4.5]{2015Chen} and Elagin \thmcite[7.1]{2014Elag} independently proved that the comparison functor is a triangle equivalence.
Their approaches differ: Elagin established the result in a general framework without giving an explicit description of the functor, whereas Chen provided a complete explicit characterization.
Later, Sun \thmcite[3.14]{2019Sun} showed that, in the stable category setting, the comparison functor is a triangle equivalence up to retracts.
Notably, all these comparison functors are special instances of Beck's theorem \cite[Chapter VI 7.1]{1998Mac}.
\end{rmk}

\begin{lem} \label{Ho iso theta}
There is a natural isomorphism $\Theta \cong \Ho(\gamma_{\calA}^G)$.
\end{lem}

\begin{prf*}
By Lemma \ref{natural iso}, it suffices to show that
$
\Theta \circ \gamma_{\calA^G} \cong \Ho(\gamma_{\calA}^G) \circ \gamma_{\calA^G}.
$
Since $\Ho(\gamma_{\calA}^G) \circ \gamma_{\calA^G} = \gamma_{\calA}^G$,
we need to prove
$
\Theta \circ \gamma_{\calA^G} \cong \gamma_{\calA}^G.
$

For each object $(X, \phi) \in \calA^G$, take a bifibrant replacement in $\calA^G$, giving a sequence
\[
\xymatrix@C=0.75CM{
(RQX, \lambda^{QX})  &  (QX, \tau^X) \ar@{>->}[l]_{\quad h_X}  \ar@{->>}[r]^{\,\,\,\,\,\, p_X}   &  (X, \phi)
}
\]
in $\calA^G$ with $h_X$ and $p_X$ weak equivalences.
For each $g \in G$, we have a commutative diagram
\[
\xymatrix@C=1.5CM@R=0.75CM{
   \rho_g RQX
\ar[d]_-{\lambda^{QX}_g}
&  \rho_g QX
\ar[d]^-{\tau^X_g}
\ar@{->>}[r]^-{\rho_g p_X}
\ar@{>->}[l]_-{\rho_g h_X}
&  \rho_g X
\ar[d]^-{\phi_g}                             \\
   RQX
&  QX
\ar@{->>}[r]^-{p_X}
\ar@{>->}[l]_-{h_X}
&  X                                                        }
\]
in $\calA$. Applying the localization functor $\gamma_{\calA} : \calA \to \Ho(\calM)$ yields the commutative diagram
\[
\xymatrix@C=1.5CM@R=0.75CM{
   \rho_g RQX
\ar[d]_-{\gamma_{\calA}(\lambda^{QX}_g)}
&  \rho_g QX
\ar[d]^-{\gamma_{\calA}(\tau^X_g)}
\ar@{->>}[r]^-{\gamma_{\calA}(\rho_g p_X)}
\ar@{>->}[l]_-{\gamma_{\calA}(\rho_g h_X)}
&  \rho_g X
\ar[d]^-{\gamma_{\calA}(\phi_g)}                             \\
   RQX
&  QX
\ar@{->>}[r]^-{\gamma_{\calA}(p_X)}
\ar@{>->}[l]_-{\gamma_{\calA}(h_X)}
&  X                                                        }
\]
in $\Ho(\calM)$.
Using the identities $\gamma_{\calA}(\rho_g h_X) = \widehat{\rho}_g \gamma_{\calA}(h_X)$ and $\gamma_{\calA}(\rho_g p_X) = \widehat{\rho}_g \gamma_{\calA}(p_X)$, we obtain the following sequence in $\Ho(\calM)^G$:
\[
\xymatrix@C=1.2CM{
   (RQX, \{\gamma_{\calA}(\lambda^{QX}_g)\})
&  (QX, \{\gamma_{\calA}(\tau^X_g)\})
\ar@{>->}[l]_-{\gamma_{\calA}(h_X)}
\ar@{->>}[r]^-{\gamma_{\calA}(p_X)}
&  (X, \{\gamma_{\calA}(\phi_g)\})
}
\]
Here $\gamma_{\calA}(h_X)$ and $\gamma_{\calA}(p_X)$ are isomorphisms.
Recall that $\gamma_{\calA}^G(X, \phi) =  (X, \{\gamma_{\calA}(\phi_g)\})$.
Moreover,
\begin{align*}
     \Theta \circ \gamma_{\calA^G}(X, \phi)
& =  \iota^G \circ \St(K^G) \circ RQ \circ \gamma_{\calA^G}(X, \phi)  \\
& =  \iota^G \circ \St(K^G)(RQX, \lambda^{QX})\\
&  =  (RQX, \{[\lambda^{QX}_g]\})\\
&=  (RQX, \{\gamma_{\calA}(\lambda^{QX}_g)\}).
\end{align*}
Thus, the composition $\gamma_{\calA}(p_X) \circ \gamma_{\calA}(h_X)^{-1}$ gives an isomorphism $\Theta \circ \gamma_{\calA^G}(X, \phi) \to \gamma_{\calA}^G(X, \phi)$.

We claim that these isomorphisms assemble into a natural isomorphism
\[
\{\gamma_{\calA}(p_X) \circ \gamma_{\calA}(h_X)^{-1}\}_{(X, \phi) \in \calA^G} :
\Theta \circ \gamma_{\calA^G} \to \gamma_{\calA}^G.
\]
To verify naturality, take any morphism $f:(X,\phi)\to(Y,\psi)$ in $\calA^G$. The following commutative diagram in $\Ho(\calM)^G$ is obtained by the usual bifibrant replacement functoriality:
\[
\xymatrix@C=1.5CM@R=0.75CM{
   (RQX,\{\gamma_{\calA}(\lambda^{QX}_g)\})
\ar[d]_-{\gamma_{\calA}(RQf)}
&  (QX,\{\gamma_{\calA}(\tau^X_g)\})
\ar@{>->}[l]_-{\gamma_{\calA}(h_X)}
\ar@{->>}[r]^-{\gamma_{\calA}(p_X)}
\ar[d]^-{\gamma_{\calA}(Qf)}
&  (X,\{\gamma_{\calA}(\phi_g)\})
\ar[d]^-{\gamma_{\calA}(f)}           \\
   (RQY,\{\gamma_{\calA}(\lambda^{QY}_g)\})
&  (QY,\{\gamma_{\calA}(\tau^Y_g)\})
\ar@{>->}[l]_-{\gamma_{\calA}(h_Y)}
\ar@{->>}[r]^{\gamma_{\calA}(p_Y)}
&  (Y,\{\gamma_{\calA}(\psi_g)\})
}
\]
Lemma \ref{fundamental Lemma} ensures that $\gamma_{\calA}(Qf)$ is a morphism of $G$-equivariant objects; its dual ensures the same for $\gamma_{\calA}(RQf)$.
Since $\gamma_{\calA}^G(f)=\gamma_{\calA}(f)$ and $\Theta\circ\gamma_{\calA^G}(f)=\gamma_{\calA}(RQf)$, the commutativity of the outer rectangle gives the required naturality.
\end{prf*}

\begin{bfhpg}[\bf Proof of Theorem \ref{retract (Q^G,W^G,R^G)}(3)]
The functor $\Theta$ is a triangle equivalence up to retracts (see \ref{con for Theta}), and becomes a genuine triangle equivalence if $\St_{\omega^G}(\calC^G\cap\calF^G)$ is idempotent complete. Hence, by Lemma \ref{Ho iso theta}, the induced functor $\Ho(\gamma_{\calA}^G)$ is also a triangle equivalence up to retracts, and is a triangle equivalence under the same idempotent completeness condition.
\qed
\end{bfhpg}

\section{Squares for homotopy categories and stable categories}
\label{Proofs of Theorem B}
\noindent
This section is devoted to the proof of Theorem \ref{thmB homotopy square}. In Subsection \ref{ThmB(1)}, we show that a $G$-functor can be lifted to homotopy and stable categories. Subsection \ref{ThmB(2)} then completes the proof of Theorem \ref{thmB homotopy square}.
{\sf Throughout this section}, let $\calA$ and $\calB$ be Grothendieck categories with enough projectives, and let $G$ be a finite group acting on both categories via $(\rho^\calA,\theta^\calA)$ and $(\rho^\calB,\theta^\calB)$, respectively.

\subsection{Induced $G$-functors on homotopy and stable categories} \label{ThmB(1)}
This subsection is devoted to proving that, under suitable conditions, a $G$-functor between Grothendieck categories equipped with Hovey triples induces $G$-functors on the corresponding homotopy and stable categories.

\begin{setup}
In this subsection, let $\calM_\calA=(\calC_\calA,\calW_\calA,\calF_\calA)$ and $\calM_\calB=(\calC_\calB,\calW_\calB,\calF_\calB)$ be hereditary Hovey triples in $\calA$ and $\calB$, respectively, such that in each triple at least two of the three classes are $G$-invariant, and denote by $\gamma_{\calA}: \calA \to \Ho(\calM_\calA)$ and $\gamma_{\calB}: \calB \to \Ho(\calM_\calB)$ the corresponding localization functor.
\end{setup}



\begin{lem} \label{G-act for St(A)}
The $G$-action $(\rho^\calA,\theta^\calA)$ on $\calA$ induces a $G$-action $(\overline{\rho^\calA},\overline{\theta^\calA})$ on the stable category $\St_{\calC_\calA \cap \calW_\calA}(\calA)$.
\end{lem}

\begin{prf*}
The construction is entirely analogous to the first part of the proof of Lemma \ref{admiss for St(C cap F)}. Specifically, for each $g\in G$, the induced functor
\[
(\overline{\rho^\calA})_g : \St_{\calC_\calA \cap \calW_\calA}(\calA) \longrightarrow \St_{\calC_\calA \cap \calW_\calA}(\calA)
\]
is defined as follows:
\begin{rqm}
\item[$\bullet$]
for an object $X \in \St_{\calC_\calA \cap \calW_\calA}(\calA)$, set $(\overline{\rho^\calA})_g X = (\rho^\calA)_g X$;
\item[$\bullet$]
for a morphism $[f] : X \to Y \in \St_{\calC_\calA \cap \calW_\calA}(\calA)$, set $(\overline{\rho^\calA})_g[f] = [(\rho^\calA)_g f] : (\rho^\calA)_g X \to (\rho^\calA)_g Y$.
\end{rqm}
For any $g,h\in G$, the natural isomorphism
\[
(\overline{\theta^\calA})_{g,h}:(\overline{\rho^\calA})_g\circ(\overline{\rho^\calA})_h \longrightarrow (\overline{\rho^\calA})_{hg}
\]
is given by
\[
(\overline{\theta^\calA})_{g,h,X} = [(\theta^\calA)_{g,h,X}] : (\rho^\calA)_g(\rho^\calA)_h X \longrightarrow (\rho^\calA)_{hg} X
\]
for every object $X\in\St_{\calC_\calA\cap\calW_\calA}(\calA)$.
\end{prf*}

\begin{rmk} \label{gamma F is G}
Let $(F,\sigma):(\calA,\rho^\calA,\theta^\calA)\to(\calB,\rho^\calB,\theta^\calB)$ be a $G$-functor.
By Lemma \ref{G-act for St(A)} we have the induced $G$-action $(\overline{\rho^\calA},\overline{\theta^\calA})$ on $\St_{\calC_\calA \cap \calW_\calA}(\calA)$, and the induced $G$-action $(\widehat{\rho^\calB},\widehat{\theta^\calB})$ on $\Ho(\calM_\calB)$ by Lemma \ref{admiss for Ho(M)}. A direct check shows that
\[
(\gamma_{\calB}F,\{\gamma_{\calB}\sigma_g\}) :
(\St_{\calC_\calA \cap \calW_\calA}(\calA),\overline{\rho^\calA},\overline{\theta^\calA})
\longrightarrow
(\Ho(\calM_\calB),\widehat{\rho^\calB},\widehat{\theta^\calB})
\]
is a $G$-functor.
Consequently, we obtain the equivariant functor
\[
(\gamma_{\calB}F)^G :
\St_{\calC_\calA \cap \calW_\calA}(\calA)^G
\longrightarrow
\Ho(\calM_\calB)^G
\]
defined as follows:
\begin{rqm}
\item[$\bullet$]
for an object $(X, \{[\phi_g]\}) \in \St_{\calC_\calA \cap \calW_\calA}(\calA)^G$, set
\[
(\gamma_\calB F)^G(X, \{[\phi_g]\}) = (\gamma_\calB(FX), \{\gamma_\calB(F \phi_g) \circ \gamma_\calB(\sigma_{g, X})^{-1}\});
\]
\item[$\bullet$]
for a morphism $[f] : (X, \{[\phi_g]\}) \to (Y, \{[\psi_g]\}) \in
\St_{\calC_\calA \cap \calW_\calA}(\calA)^G$, set
\[(\gamma_\calB F)^G([f]) = \gamma_\calB(F f).\]
\end{rqm}
\end{rmk}

\begin{lem} \label{G-functor Q}
For the cofibrant replacement functor $Q$, there exists a family $\{\sigma_g^Q\}_{g\in G}$ of natural isomorphisms such that $(Q,\sigma^Q)$ is a $G$-functor. Consequently, the equivariant functor
$
Q^G:\calA^G \rightarrow \St_{\calC_\calA \cap \calW_\calA}(\calA)^G
$
is well-defined.
\end{lem}

\begin{prf*}
For each $g\in G$, we construct a natural isomorphism
$
\sigma^Q_g : Q\circ(\rho^\calA)_g \rightarrow (\overline{\rho^\calA})_g\circ Q.
$
Diagrammatically, this is a natural transformation filling the square
\[
\xymatrix@C=1.5CM@R=0.75CM{
   \calA
\ar[d]_-{Q}
\ar[r]^-{(\rho^\calA)_g}
&  \calA
\ar[d]^-{Q}             \\
   \St_{\calC_\calA \cap \calW_\calA}(\calA)
\ar[r]^-{(\overline{\rho^\calA})_g}
&  \St_{\calC_\calA \cap \calW_\calA}(\calA)                         }
\]
Fix an object $X \in \calA$.
Since $(\calC_\calA, \calW_\calA \cap \calF_\calA)$ is a complete cotorsion pair in $\calA$, we have short exact sequences of rows
\[
\xymatrix@C=1.5CM@R=0.75CM{
   W_{(\rho^\calA)_g X}
\ar@{-->}[d]
\ar@{>->}[r]
&  Q(\rho^\calA)_g X
\ar@{-->}[d]_-{\sigma^Q_{g, X}}
\ar@{->>}[r]
&  (\rho^\calA)_g X
\ar@{=}[d]                       \\
   (\rho^\calA)_g W_X
\ar@{>->}[r]
&  (\rho^\calA)_g QX
\ar@{->>}[r]
&  (\rho^\calA)_g X                               }
\]
in $\calA$, where $Q(\rho^\calA)_g X, QX \in \calC_\calA$ and $W_{(\rho^\calA)_g X}, W_X\in \calW_\calA \cap \calF_\calA$.
By the $G$-invariance assumptions, we also have $(\rho^\calA)_g QX\in\calC_\calA$ and $(\rho^\calA)_g W_X\in\calW_\calA\cap\calF_\calA$.
Hence, the dashed arrows exist and make the diagram commute.
We then define
\[
\sigma^Q_g(X)=[\sigma^Q_{g,X}] :
Q(\rho^\calA)_g X \longrightarrow (\rho^\calA)_g QX
\]
in the stable category $\St_{\calC_\calA\cap\calW_\calA}(\calA)$. The uniqueness part of Lemma \ref{fundamental Lemma} ensures that $\sigma^Q_g(X)$ is an isomorphism.

The naturality of $\sigma^Q_g$ follows by an argument analogous to the second step of the proof of Lemma \ref{admiss for St(C cap F)}; hence each $\sigma^Q_g$ is a natural isomorphism. Moreover, a direct verification shows that the family $\{\sigma^Q_g\}_{g\in G}$ satisfies the associativity and unit conditions of Definition \ref{G fun}. Therefore, $(Q,\sigma^Q)$ is a $G$-functor.
\end{prf*}

\begin{ipg} \label{total LF}
Let $(F,H)$ be a Quillen adjunction between $\calM_\calA$ and $\calM_\calB$ with $(F,\sigma^F):(\calA,\rho^\calA,\theta^\calA)\to (\calB,\rho^\calB,\theta^\calB)$ a $G$-functor. Since $F$ preserves trivial cofibrations, the composition
\[\gamma_{\calB}\circ F : \calA \longrightarrow \Ho(\calM_{\calB})\]
sends trivial cofibrations between cofibrant objects to isomorphisms. In particular, $\gamma_{\calB}\circ F$ sends trivial cofibrant objects to zero objects in $\Ho(\calM_{\calB})$. Hence, it descends to a functor
\[\gamma_{\calB} F : \St_{\calC_\calA \cap \calW_\calA}(\calA) \longrightarrow \Ho(\calM_{\calB}),\]
which acts on objects by $X\mapsto F(X)$ and on morphisms by $[f]\mapsto \gamma_{\calB}(Ff)$.
Therefore, the Localization Theorem \ref{loc thm} applies to the composition $\gamma_{\calB} F\circ Q$, where $Q$ is a cofibrant replacement functor. We thus obtain a unique functor
\[
\mathbb{L}F : \Ho(\calM_\calA) \longrightarrow \Ho(\calM_\calB),
\]
called the \emph{total left derived functor} of $F$, making the following diagram commute:
\[
\xymatrix@R=0.75cm{
   \calA
\ar[r]^-{\gamma_{\calA}}
\ar[d]_{Q}
&  \Ho(\calM_\calA)
\ar@{-->}[d]^-{\mathbb{L}F} \\
   \St_{\calC_\calA \cap \calW_\calA}(\calA)
\ar[r]^-{\gamma_{\calB} F}
&  \Ho(\calM_\calB)
}
\]
For more details, we refer the reader to \thmcite[7.11]{2025GillBook}.
\end{ipg}

\begin{ipg} \label{equ for F}
We recall from the proof of Lemma \ref{admiss for Ho(M)} (see \ref{detailed proof}) the commutative diagram:
\[
\xymatrix{
   \calA
\ar[r]^-{F}
&  \calB
\ar[d]_-{(\rho^\calB)_g}
\ar[r]^-{\gamma_\calB}
&  \Ho(\calM_\calB)
\ar[d]^-{(\widehat{\rho^\calB})_g}            \\
&  \calB
\ar[r]^-{\gamma_\calB}
&  \Ho(\calM_\calB)                             }
\]
As shown in \ref{total LF}, the composite $\gamma_\calB\circ F$ yields a functor $\gamma_\calB F : \St_{\calC_\calA \cap \calW_\calA}(\calA) \longrightarrow \Ho(\calM_\calB)$.
We will use the above diagram to define an additional descent functor and to relate it to $\gamma_\calB F$.

Since any two of the classes $\calC_\calB$, $\calW_\calB$ and $\calF_\calB$ are $G$-invariant by assumption, each $(\rho^\calB)_g$ preserves trivial cofibrations.
Consequently, so does the composite $(\rho^\calB)_g \circ F$.
It follows that the composition
\[
\gamma_{\calB} \circ (\rho^\calB)_g \circ F : \calA \to \Ho(\calM_{\calB})
\]
sends trivial cofibrations to isomorphisms.
In particular, this functor sends every trivially cofibrant object to a zero object in $\Ho(\calM_{\calB})$.
Hence it descends to a well-defined functor
\[
\gamma_{\calB} (\rho^\calB)_g F : \St_{\calC_\calA \cap \calW_\calA}(\calA) \to \Ho(\calM_{\calB})
\]
given on objects by $X\mapsto (\rho^\calB)_g F(X)$ and on morphisms by $[f]\mapsto \gamma_\calB((\rho^\calB)_g Ff)$.

A direct verification then yields the identity
\[
(\widehat{\rho^\calB})_g \circ \gamma_\calB F = \gamma_\calB(\rho^\calB)_g F,
\]
which will be used in the proof of Proposition \ref{G-functor for LF}.
\end{ipg}

The following result shows that, under suitable conditions, a $G$-functor between $\calA$ and $\calB$ induces $G$-functors at the level of the associated homotopy categories.

\begin{prp} \label{G-functor for LF}
Let $(F, H)$ be a Quillen adjunction between $\calM_\calA$ and $\calM_\calB$ with $(F, \sigma^F): (\calA, \rho^\calA, \theta^\calA) \to (\calB, \rho^\calB, \theta^\calB)$ a $G$-functor.
Then there exists a family $\sigma^{\mathbb{L}F} = \{\sigma^{\mathbb{L}F}_g\}_{g\in G}$ of natural isomorphisms such that $(\mathbb{L}F, \sigma^{\mathbb{L}F}): \Ho(\calM_\calA) \to \Ho(\calM_\calB)$ is a triangulated $G$-functor.
In particular, the equivariant functor
\[
\mathbb{L}(F)^G:
\Ho(\calM_\calA)^G
\to
\Ho(\calM_\calB)^G
\]
is well-defined.

If, in addition, $|G|$ is invertible in both $\calA$ and $\calB$, then $\mathbb{L}(F)^G$ is triangulated, where both equivariant categories are equipped with their unique canonical pre-triangulated structures.
Moreover, $\mathbb{L}(F)^G$ is a triangle equivalence under the following conditions:
\begin{prt}
\item Both $F$ and $H$ are exact and faithful;
\item $H(\calC_\calB\cap\calW_\calB)\subseteq\calW_\calA$ and
      $F(\calW_\calA\cap\calF_\calA)\subseteq\calW_\calB$;
\item $\coker(\eta_X)\in\calW_\calA$ for every $X\in\calC_\calA$, and
      $\ker(\varepsilon_Y)\in\calW_\calB$ for every $Y\in\calF_\calB$.
\end{prt}
\end{prp}

\begin{prf*}
The total left derived functor $\mathbb{L}F:\Ho(\calM_\calA)\to\Ho(\calM_\calB)$ exists and satisfies
$
\mathbb{L}F\circ\gamma_{\calA}=\gamma_{\calB}F\circ Q;
$
see \ref{total LF}.
For each $g\in G$, we need to construct a natural isomorphism
\[
\sigma_g^{\mathbb{L}F}:
\mathbb{L}F\circ(\widehat{\rho^\calA})_g
\longrightarrow
(\widehat{\rho^\calB})_g\circ\mathbb{L}F,
\]
i.e., the square commutes up to natural isomorphism
\[
\xymatrix@C=1.25cm{
   \Ho(\calM_\calA) \ar[d]_{\mathbb{L}F} \ar[r]^{(\widehat{\rho^\calA})_g}
&  \Ho(\calM_\calA) \ar[d]^{\mathbb{L}F} \\
   \Ho(\calM_\calB) \ar[r]_{(\widehat{\rho^\calB})_g}
&  \Ho(\calM_\calB).
}
\]
By Lemma \ref{natural iso}, it suffices to produce a natural isomorphism
\[
\mathbb{L}F\circ(\widehat{\rho^\calA})_g\circ\gamma_{\calA}
\cong
(\widehat{\rho^\calB})_g\circ\mathbb{L}F\circ\gamma_{\calA}.
\]
Using the identities
$
(\widehat{\rho^\calA})_g\circ\gamma_{\calA}=\gamma_{\calA}\circ(\rho^\calA)_g
$
and
$(\widehat{\rho^\calB})_g\circ\gamma_{\calB}=\gamma_{\calB}\circ(\rho^\calB)_g,
$
together with \ref{equ for F}, we obtain
\begin{align*}
\mathbb{L}F\circ(\widehat{\rho^\calA})_g\circ\gamma_{\calA}
&= \mathbb{L}F\circ\gamma_{\calA}\circ(\rho^\calA)_g
   = \gamma_\calB F\circ Q\circ(\rho^\calA)_g, \\
(\widehat{\rho^\calB})_g\circ\mathbb{L}F\circ\gamma_{\calA}
&= (\widehat{\rho^\calB})_g\circ\gamma_\calB F\circ Q
   = \gamma_\calB(\rho^\calB)_g F\circ Q.
\end{align*}
Thus it remains to construct an isomorphism
\[
\gamma_\calB F\circ Q\circ(\rho^\calA)_g
\cong
\gamma_\calB(\rho^\calB)_g F\circ Q.
\]

Consider the diagram
\[
\xymatrix@C=1.25CM@R=1.25CM{
   \calA \ar[d]_{Q} \ar[r]^{(\rho^\calA)_g}
&  \calA \ar[r]^{Q\ \ \ \ \ \ }
&  \St_{\calC_\calA \cap \calW_\calA}(\calA) \ar[d]^{\gamma_\calB F} \\
   \St_{\calC_\calA \cap \calW_\calA}(\calA)
   \ar[rr]_{\gamma_\calB(\rho^\calB)_g F}
   \ar@{-->}[urr]^{(\rho^\calA)_g}
&& \Ho(\calM_\calB)
}
\]
Since $(\rho^\calA)_g$ preserves trivial cofibrant objects, it induces the dashed arrow
\[
(\rho^\calA)_g:\St_{\calC_\calA\cap\calW_\calA}(\calA)
\to
\St_{\calC_\calA\cap\calW_\calA}(\calA).
\]
Because $\sigma^F_g:F\circ(\rho^\calA)_g\to(\rho^\calB)_g\circ F$ is a natural isomorphism, one checks directly that
\[
\gamma_\calB\sigma_g:
\gamma_\calB F\circ(\rho^\calA)_g
\cong
\gamma_\calB(\rho^\calB)_g F.
\]

We next construct a natural isomorphism
\[
Q\circ(\rho^\calA)_g \cong (\rho^\calA)_g\circ Q.
\]
Once this is done, the desired $\sigma_g^{\mathbb{L}F}$ follows.

For an object $X\in\calA$, recall that $(\calC_\calA,\calW_\calA\cap\calF_\calA)$ is a complete cotorsion pair.
We have the commutative diagram
\[
\xymatrix{
   W_{(\rho^\calA)_g X} \ar@{-->}[d] \ar@{>->}[r]
&  Q(\rho^\calA)_g X \ar@{-->}[d]_{h_X} \ar@{->>}[r]
&  (\rho^\calA)_g X \ar@{=}[d] \\
   (\rho^\calA)_g W_X \ar@{>->}[r]
&  (\rho^\calA)_g QX \ar@{->>}[r]
&  (\rho^\calA)_g X
}
\]
where $Q(\rho^\calA)_g X,(\rho^\calA)_g QX\in\calC_\calA$ and
$W_{(\rho^\calA)_g X},(\rho^\calA)_g W_X\in\calW_\calA\cap\calF_\calA$.
By the uniqueness part of Lemma \ref{fundamental Lemma}, $[h_X]$ is an isomorphism in the stable category $\St_{\calC_\calA\cap\calW_\calA}(\calA)$. We claim that these assemble into a natural isomorphism
\[
\{[h_X]:Q(\rho^\calA)_g X\to(\rho^\calA)_g QX\}_{X\in\calA}
:
Q\circ(\rho^\calA)_g
\longrightarrow
(\rho^\calA)_g\circ Q.
\]
To verify naturality, let $f:X\to Y$ be a morphism in $\calA$. We must show that the square
\[
\xymatrix@C=1.75CM{
   Q(\rho^\calA)_g X \ar[d]_{[h_X]} \ar[r]^{[Q(\rho^\calA)_g f]}
&  Q(\rho^\calA)_g Y \ar[d]^{[h_Y]} \\
   (\rho^\calA)_g QX \ar[r]^{[(\rho^\calA)_g Qf]}
&  (\rho^\calA)_g QY
}
\]
commutes in $\St_{\calC_\calA\cap\calW_\calA}(\calA)$. Consider the diagram with exact rows:
\[
\xymatrix@C=0.7CM @R=0.6CM{
&     W_{(\rho^\calA)_g X}   \ar@{>->}[rr]     \ar'[d][dd]              \ar[dl]
&  &  Q(\rho^\calA)_g X     \ar@{->>}[rr]     \ar'[d][dd]_(0.3){h_X}   \ar[dl]_{Q(\rho^\calA)_g f}
&  &  (\rho^\calA)_g X       \ar@{=}[dd]                                \ar[dl]_{(\rho^\calA)_g f}          \\
      W_{(\rho^\calA)_g Y}   \ar@{>->}[rr]     \ar[dd]
&  &  Q(\rho^\calA)_g Y     \ar@{->>}[rr]     \ar[dd]_(0.3){h_Y}
&  &  (\rho^\calA)_g Y       \ar@{=}[dd]                                                                  \\
   &  (\rho^\calA)_g W_X     \ar@{>->}'[r][rr]                          \ar[dl]
&  &  (\rho^\calA)_g QX     \ar@{->>}'[r][rr]                          \ar[dl]_{(\rho^\calA)_g Qf}
&  &  (\rho^\calA)_g X                                                  \ar[dl]_{(\rho^\calA)_g f}          \\
      (\rho^\calA)_g W_Y     \ar@{>->}[rr]
&  &  (\rho^\calA)_g QY     \ar@{->>}[rr]
&  &  (\rho^\calA)_g Y                     }
\]
The upper and lower squares are commutative. The uniqueness part of Lemma \ref{fundamental Lemma} then yields
\[
[h_Y]\circ[Q(\rho^\calA)_g f]
=
[(\rho^\calA)_g Qf]\circ[h_X],
\]
as required.

Finally, a routine verification shows that the family $\{\sigma^{\mathbb{L}F}_g\}_{g\in G}$ satisfies the associativity and unit conditions of Definition \ref{G fun}.
Hence $(\mathbb{L}F,\sigma^{\mathbb{L}F})$ is a $G$-functor.
By \thmcite[7.12]{2025GillBook}, $\mathbb{L}F$ is a triangle functor.
Let $\xi : \mathbb{L}F\mathsf{\Sigma}_{\calA} \overset{\cong}{\rightarrow} \mathsf{\Sigma}_{\calB}\mathbb{L}F$ be its structure isomorphism, and denote by $\widehat{\varphi}^{\calA}_g$ and $\widehat{\varphi}^{\calB}_g$ the
admissible structures on $\Ho(\calM_{\calA})$ and $\Ho(\calM_{\calB})$, respectively.

We claim that $\xi$ is $G$-natural.
Indeed, the triangle structure on $\mathbb{L}F$ is induced by the suspension construction on cofibrant
replacements.
By the construction above, $\sigma^{\mathbb{L}F}_g$ is induced by the natural isomorphisms
\[
\sigma^Q_g:
Q\circ(\rho^{\calA})_g
\overset{\cong}{\rightarrow}
(\rho^{\calA})_g\circ Q
\quad  \text{and} \quad
\sigma^F_g:
F\circ(\rho^{\calA})_g
\overset{\cong}{\rightarrow}
(\rho^{\calB})_g\circ F.
\]
These natural isomorphisms induce morphisms between the corresponding suspension constructions.
Hence, by the naturality of the connecting morphism and the construction of $\widehat{\varphi}^{\calA}_g$ and
$\widehat{\varphi}^{\calB}_g$ in Lemma \ref{admiss for Ho(M)}, the following diagram is commutative for every
$g\in G$:
\[
\xymatrix@C=1.8cm{
\mathbb{L}F\mathsf{\Sigma}_{\calA}(\widehat{\rho^{\calA}})_g
\ar[r]^-{\xi(\widehat{\rho^{\calA}})_g}
\ar[d]_-{(\sigma^{\mathbb{L}F}_g\mathsf{\Sigma}_{\calA})
\circ(\mathbb{L}F\widehat{\varphi}^{\calA}_g)}
&
\mathsf{\Sigma}_{\calB}\mathbb{L}F(\widehat{\rho^{\calA}})_g
\ar[d]^-{(\widehat{\varphi}^{\calB}_g\mathbb{L}F)
\circ(\mathsf{\Sigma}_{\calB}\sigma^{\mathbb{L}F}_g)}
\\
(\widehat{\rho^{\calB}})_g\mathbb{L}F\mathsf{\Sigma}_{\calA}
\ar[r]^-{(\widehat{\rho^{\calB}})_g\xi}
&
(\widehat{\rho^{\calB}})_g\mathsf{\Sigma}_{\calB}\mathbb{L}F .
}
\]
Thus $\xi$ is $G$-natural.
Therefore, $(\mathbb{L}F,\sigma^{\mathbb{L}F})$ is a triangulated $G$-functor.
In particular, it induces the equivariant functor $\mathbb{L}(F)^G : \Ho(\calM_{\calA})^G \to \Ho(\calM_{\calB})^G$.

If $|G|$ is invertible in both $\calA$ and $\calB$, then it is also invertible in $\Ho(\calM_{\calA})$ and $\Ho(\calM_{\calB})$.
Hence Lemma \ref{G-pretri fun} implies that $\mathbb{L}(F)^G$ is triangulated with respect to the unique canonical pre-triangulated structures.

Under the additional assumptions {\rm (a)--(c)}, Theorem \ref{Ho(M_A)=Ho(M_B)} implies that $(F,H)$ is a Quillen
equivalence.
Hence $\mathbb{L}F$ is an equivalence.
By Lemma \ref{G-equivalence}, $\mathbb{L}(F)^G$ is also an equivalence.
Consequently, $\mathbb{L}(F)^G$ is a triangle equivalence.
\end{prf*}

\subsection{Proof of Theorem \ref{thmB homotopy square}} \label{ThmB(2)}
In this subsection we prove Theorem \ref{thmB homotopy square}. As a first step, Proposition \ref{Quillen adj L(F^G)} records that a Quillen adjunction lifts to one between the equivariant model structures.

\begin{setup}
Throughout this subsection, assume that $|G|$ is invertible in both $\calA$ and $\calB$. Let
$\calM_\calA=(\calC_\calA,\calW_\calA,\calF_\calA)$
and
$\calM_\calB=(\calC_\calB,\calW_\calB,\calF_\calB)$
be cofibrantly generated hereditary Hovey triples in $\calA$ and $\calB$, respectively, such that in each triple at least two of the three classes are $G$-invariant. Denote by
$
\gamma_{\calA}:\calA\to\Ho(\calM_\calA)
$ and
$\gamma_{\calB}:\calB\to\Ho(\calM_\calB)
$
the corresponding localization functors.
\end{setup}

By Theorem \ref{retract (Q^G,W^G,R^G)}(1), one gets that both
\[
\calM_\calA^G = (\calC_\calA^G, \calW_\calA^G, \calF_\calA^G)
\quad \text{and} \quad
\calM_\calB^G = (\calC_\calB^G, \calW_\calB^G, \calF_\calB^G)
\]
are cofibrantly generated hereditary Hovey triples in $\calA^G$ and $\calB^G$, respectively.

The following result shows that, under suitable conditions, a Quillen adjunction between $\calM_\calA$ and $\calM_\calB$ induces a Quillen adjunction between the corresponding equivariant abelian model structures.

\begin{prp} \label{Quillen adj L(F^G)}
Let $(F, H)$ be a Quillen adjunction between $\calM_\calA$ and $\calM_\calB$.
If either $F$ or $H$ is a $G$-functor, then $(F^G, H^G)$ is a Quillen adjunction between $\calM_\calA^G$ and $\calM_\calB^G$.
Consequently, the total left derived functor $\mathbb{L}(F^G) : \Ho(\calM_\calA^G) \to \Ho(\calM_\calB^G)$ is well-defined.
Moreover, $\mathbb{L}(F^G)$ is a triangle equivalence under the following additional conditions:
\begin{prt}
  \item $H :\calB \rightleftarrows \calA: F$ is an adjoint pair with both $F$ and $H$ faithful;
  \item $H(\calC_\calB \cap \calW_\calB) \subseteq \calW_\calA$ and $F(\calW_\calA \cap \calF_\calA) \subseteq \calW_\calB$;
  \item $\cok(\eta_X) \in \calW_\calA$ for every $X \in \calC_\calA$, and $\kernel(\varepsilon_Y) \in \calW_\calB$ for every $Y \in \calF_\calB$.
 \end{prt}
\end{prp}

\begin{prf*}
By Lemma \ref{G-frobenius pair}, the pair $(F^G, H^G)$ forms an adjunction between $\calA^G$ and $\calB^G$.
We show that $F^G$ preserves cofibrations; the preservation of trivial cofibrations follows by the same argument.
Once these facts are established, it follows that $(F^G, H^G)$ is a Quillen adjunction between $\calM_\calA^G$ and $\calM_\calB^G$, and hence the total left derived functor $\mathbb{L}(F^G)$ exists.

Let $f: (X, \phi) \to (Y, \psi)$ be a cofibration.
Then $f$ is a monomorphism whose cokernel is cofibrant.
Thus we obtain a short exact sequence in $\calA^G$:
\[
0 \to (X, \phi) \overset{f}{\rightarrow} (Y, \psi) \to (\cok(f), \mu) \to 0,
\]
where $(\cok(f), \mu) \in \calC_\calA^G$ and the linearization $\mu$ is induced by the universal property of the cokernel; see \rmkcite[2.6]{2025Equ}.
To prove that $F^G$ preserves cofibrations, it remains to show that $F^G(\cok(f), \mu) \in \calC_\calB^G$ and that $F^G(f)$ is a monomorphism.
Since $F^G(\cok(f), \mu) = (F(\cok(f)), \mu')$ with $\mu'$ as in \ref{equi fun}, and $F^G(f)=F(f)$, it suffices to verify that $F(\cok(f)) \in \calC_\calB$ and that $F(f)$ is a monomorphism.

Applying the forgetful functor $U_\calA$ to the above short exact sequence yields a short exact sequence in $\calA$:
\[
0 \to X \overset{f}\rightarrow Y \to \cok(f) \to 0,
\]
where $\cok(f) \in \calC_\calA$.
Since $F$ preserves cofibrations, we obtain a short exact sequence
\[
\xymatrix@C=0.8cm{
0  \to  F(X) \ar[r]^{F(f) \qquad \quad\,\,\,}  &  F(Y)  \to  F(\cok(f))  \to  0 }
\]
in $\calB$, where $F(\cok(f)) \in \calC_\calB$, as desired.

We now prove that $\mathbb{L}(F^G)$ is a triangle equivalence.
By \thmcite[7.12]{2025GillBook}, $\mathbb{L}(F^G)$ is a triangle functor.
We apply Theorem \ref{Ho(M_A)=Ho(M_B)} to show that $(F^G,H^G)$ is a Quillen equivalence.
If this is established, then $\mathbb{L}(F^G)$ is an equivalence, hence a triangle equivalence.
It suffices to verify the following conditions:
\begin{rqm}
\item
Both $F^G$ and $H^G$ are exact and faithful.
\item
$H^G(\calC_\calB^G \cap \calW_\calB^G) \subseteq \calW_\calA^G$ and $F^G(\calW_\calA^G \cap \calF_\calA^G) \subseteq \calW_\calB^G$.
\item
$\cok(\eta^G_{(X, \phi)}) \in \calW_\calA^G$ for each $(X, \phi) \in \calC_\calA^G$, and $\ker(\varepsilon^G_{(Y, \psi)}) \in \calW_\calB^G$ for each $(Y, \psi) \in \calF_\calB^G$.
\end{rqm}

It is clear that $F^G$ and $H^G$ are exact.
Since $F^G(f)=F(f)$ for each $f$ in $\calA^G$, faithfulness of $F$ implies that of $F^G$; the same holds for $H^G$.
Condition (2) follows directly from assumption (b) by a routine check.
For (3), we only prove the cokernel condition for the unit; the kernel condition for the counit is analogous.

For any object $(X,\phi)\in\calC_\calA^G$, consider the short exact sequence
\[
\xymatrix@C=0.9cm{
0 \ar[r] & (X,\phi) \ar[r]^{\eta^G_{(X,\phi)} \quad \,\,} & H^G F^G(X,\phi) \ar[r] & \cok(\eta^G_{(X,\phi)}) \ar[r] & 0.
}
\]
We need to show $\cok(\eta^G_{(X,\phi)})\in\calW_\calA^G$.
Recall from \ref{equi cat} that $U_\calB \circ F^G = F \circ U_\calA$ and $U_\calA \circ H^G = H \circ U_\calB$.
Applying $U_\calA$ to the above sequence yields
\[
U_\calA(H^G F^G(X,\phi)) = H F U_\calA(X,\phi) = H F(X),
\]
and by \lemcite[2.9]{2019Sun}, $U_\calA(\eta^G_{(X,\phi)}) = \eta_X$.
Since $U_\calA$ is exact, the resulting short exact sequence in $\calA$ is
\[
\xymatrix@C=0.5cm{
0 \ar[r] & X \ar[r]^{\eta_X \quad} & H F(X) \ar[r] & U_\calA(\cok(\eta^G_{(X,\phi)})) \ar[r] & 0.
}
\]
Thus $U_\calA(\cok(\eta^G_{(X,\phi)})) \cong \cok(\eta_X)$.
Because $(X,\phi)\in\calC_\calA^G$, we have $X = U_\calA(X,\phi)\in\calC_\calA$.
By the assumption (c), $\cok(\eta_X)\in\calW_\calA$; hence $U_\calA(\cok(\eta^G_{(X,\phi)}))\in\calW_\calA$, which implies $\cok(\eta^G_{(X,\phi)})\in\calW_\calA^G$, as desired.
\end{prf*}

\begin{rmk} \label{descend and com}
Recall the following commutative diagram from \ref{exis for Ho^G}:
\[
\xymatrix{
   \calA^G \ar[r]^-{F^G} & \calB^G \ar[d]_(0.45){\gamma_\calB^G} \ar[r]^-{\gamma_{\calB^G}} & \Ho(\calM_\calB^G) \ar[dl]^(0.4){\Ho(\gamma_\calB^G)} \\
   & \Ho(\calM_\calB)^G
}
\]
We will use this diagram to define two descent functors and to establish their compatibility, which will be used in the proof of Theorem \ref{thmB homotopy square}.

Consider first the composite $\gamma_\calB^G \circ F^G: \calA^G \to \Ho(\calM_\calB)^G$.
By Proposition \ref{Quillen adj L(F^G)}, $F^G$ is left Quillen, hence preserves all trivial cofibrations.
Moreover, by \ref{exis for Ho^G}, $\gamma_\calB^G$ sends weak equivalences in $\calB^G$ to isomorphisms.
Consequently, the composite sends trivial cofibrations to isomorphisms; in particular, it sends every trivially cofibrant object to a zero object in $\Ho(\calM_\calB)^G$.
Hence, it descends to a functor
\[
\gamma_\calB^G \circ F^G: \St_{\calC_\calA^G \cap \calW_\calA^G}(\calA^G) \to \Ho(\calM_{\calB})^G,
\]
which is defined as follows:
\begin{rqm}
\item[$\bullet$]
for an object $(X, \phi) \in \St_{\calC_\calA^G \cap \calW_\calA^G}(\calA^G)$, $\gamma_\calB^G \circ F^G(X, \phi) = (FX, \{\gamma_\calB(F \phi_g) \circ \gamma_\calB((\sigma_{g, X})^{-1})\})$;
\item[$\bullet$]
for a morphism $[f] : (X, \phi) \to (Y, \psi)$ in $\St_{\calC_\calA^G \cap \calW_\calA^G}(\calA^G)$, $\gamma_\calB^G \circ F^G([f]) = \gamma_\calB(F f)$.
\end{rqm}

Next consider the composite $\gamma_{\calB^G} \circ F^G : \calA^G \to \Ho(\calM_\calB^G)$.
By the same argument as above, it descends to a functor
\[
\gamma_{\calB^G} \circ F^G: \St_{\calC_\calA^G \cap \calW_\calA^G}(\calA^G) \to \Ho(\calM_{\calB}^G),
\]
defined as follows:
\begin{rqm}
\item[$\bullet$]
for an object $(X, \phi) \in \St_{\calC_\calA^G \cap \calW_\calA^G}(\calA^G)$, $\gamma_{\calB^G} \circ F^G(X, \phi) = (FX, \{F \phi_g \circ (\sigma_{g, X})^{-1}\})$;
\item[$\bullet$]
for a morphism $[f] : (X, \phi) \to (Y, \psi)$ in $\St_{\calC_\calA^G \cap \calW_\calA^G}(\calA^G)$, $\gamma_{\calB^G} \circ F^G([f]) = \gamma_{\calB^G}(F(f))$.
\end{rqm}

We claim that
\[
\Ho(\gamma_\calB^G) \circ \gamma_{\calB^G} \circ F^G = \gamma_\calB^G\circ F^G.
\]
Indeed, for each $(X,\phi) \in \St_{\calC_\calA^G \cap \calW_\calA^G}(\calA^G)$, one has
\begin{align*}
\Ho(\gamma_\calB^G) \circ \gamma_{\calB^G}\circ F^G(X,\phi)
&= \Ho(\gamma_\calB^G)\bigl(FX,\{F\phi_g \circ (\sigma_{g,X})^{-1}\}\bigr) \\
&= \Ho(\gamma_\calB^G)\Bigl(\gamma_{\calB^G}\bigl(FX,\{F\phi_g \circ (\sigma_{g,X})^{-1}\}\bigr)\Bigr) \\
&= \gamma_\calB^G\bigl(FX,\{F\phi_g \circ (\sigma_{g,X})^{-1}\}\bigr) \\
&= \bigl(FX,\;\{\gamma_\calB(F\phi_g) \circ \gamma_\calB((\sigma_{g,X})^{-1})\}\bigr) \\
&= \gamma_\calB^G\circ F^G(X,\phi),
\end{align*}
where the second equality holds as $\calB^G$ and $\Ho(\calM_{\calB}^G)$ share the same objects, and $\gamma_{\calB^G}$ is the identity on objects.

For a morphism $[f] : (X, \phi) \to (Y, \psi)$ in $\St_{\calC_\calA^G \cap \calW_\calA^G}(\calA^G)$, we have
\begin{align*}
   \Ho(\gamma_\calB^G) \circ \gamma_{\calB^G}\circ F^G([f])
=  \Ho(\gamma_\calB^G)(\gamma_{\calB^G}(F(f)))
=  \gamma_\calB^G(F(f))
=  \gamma_\calB(F(f))
=  \gamma_\calB^G  F^G([f]).
\end{align*}
Thus the claimed equality holds.
\end{rmk}

We now give the proof of Theorem \ref{thmB homotopy square}.

\begin{bfhpg}[\bf Proof of Theorem \ref{thmB homotopy square}]  \label{thmB origin}
By Theorem \ref{retract (Q^G,W^G,R^G)}, there exist comparison functors
\[
\Ho(\gamma_{\calA}^G): \Ho(\calM_\calA^G) \to \Ho(\calM_\calA)^G
\quad\text{and}\quad
\Ho(\gamma_\calB^G): \Ho(\calM_\calB^G) \to \Ho(\calM_\calB)^G,
\]
both of which are triangle equivalences up to retracts.
Proposition \ref{Quillen adj L(F^G)} yields the total left derived functor
\[
\mathbb{L}(F^G) : \Ho(\calM_\calA^G) \to \Ho(\calM_\calB^G),
\]
and Proposition \ref{G-functor for LF} shows that $\mathbb{L}F$ induces an equivariant functor
\[
\mathbb{L}(F)^G : \Ho(\calM_\calA)^G \to \Ho(\calM_\calB)^G.
\]

We now prove that there is a natural isomorphism
\[
\Ho(\gamma_\calB^G) \circ \mathbb{L}(F^G) \cong \mathbb{L}(F)^G \circ \Ho(\gamma_{\calA}^G).
\]
By Lemma \ref{natural iso}, it suffices to establish
\[
      \Ho(\gamma_\calB^G) \circ \mathbb{L}(F^G) \circ \gamma_{\calA^G}
\cong \mathbb{L}(F)^G \circ \Ho(\gamma_{\calA}^G) \circ \gamma_{\calA^G}.
\]
For the left-hand side,
\begin{align*}
& \phantom{\quad \,\,}  \Ho(\gamma_\calB^G) \circ \mathbb{L}(F^G) \circ \gamma_{\calA^G}   \\
& =  \Ho(\gamma_\calB^G) \circ \gamma_{\calB^G} F^G \circ Q_{\calA^G}
& \text{by Remark } \ref{total LF}                                                         \\
& =  \gamma_\calB^G F^G \circ Q_{\calA^G}.
& \text{by Remark } \ref{descend and com}
\end{align*}
For the right-hand side,
\begin{align*}
& \phantom{\quad \,\,}  \mathbb{L}(F)^G \circ \Ho(\gamma_{\calA}^G) \circ \gamma_{\calA^G}      \\
& =  \mathbb{L}(F)^G \circ \gamma_{\calA}^G
& \text{by } \ref{exis for Ho^G}                                                            \\
& =  (\mathbb{L}(F) \circ \gamma_{\calA})^G
& \text{by Lemma } \ref{G-fun comp}                                                        \\
& =  (\gamma_\calB F \circ Q_\calA)^G
& \text{by Remark } \ref{total LF}                                                            \\
& =  (\gamma_\calB F)^G \circ Q_\calA^G.
& \text{by Lemma } \ref{G-functor Q} \text{ and Remark } \ref{gamma F is G}
\end{align*}

Therefore, it remains to prove the natural isomorphism
\[
\gamma_\calB^G F^G \circ Q_{\calA^G} \cong (\gamma_\calB F)^G \circ Q_\calA^G,
\]
i.e., the commutativity of the following diagram:
\[
\xymatrix@C=1.5CM{
   \calA^G
\ar[d]_-{Q_\calA^G}
\ar[r]^-{Q_{\calA^G}}
&  \St_{\calC_\calA^G \cap \calW_\calA^G}(\calA^G)
\ar[d]^-{\gamma_\calB^G F^G}               \\
   \St_{\calC_\calA \cap \calW_\calA}(\calA)^G
\ar[r]^-{(\gamma_\calB F)^G}
&  \Ho(\calM_\calB)^G                                                                         }
\]

\vspace{1ex}\textbf{Step 1. Establish the object-wise isomorphism.}

Fix $(X,\phi)\in\calA^G$.
We first compute the value of $\gamma_\calB^G F^G \circ Q_{\calA^G}$ at $(X, \phi)$.
Since $(\calC_\calA^G,\calW_\calA^G\cap\calF_\calA^G)$ is a complete cotorsion pair in $\calA^G$, there exists a trivial fibration
\[
p_X:(Q_{\calA^G}X,\tau^X)\twoheadrightarrow(X,\phi)
\]
in $\calA^G$ with $(Q_{\calA^G}X,\tau^X)\in\calC_\calA^G$.
Hence $Q_{\calA^G}(X, \phi) = (Q_{\calA^G} X, \tau^X)$.
Note that
\[
[p_X] : (Q_{\calA^G}X, \tau^X) \twoheadrightarrow (X, \phi)
\]
is a morphism in $\St_{\calC_\calA^G \cap \calW_\calA^G}(\calA^G)$.
Applying the functor $\gamma_\calB^G F^G$ to this morphism yields
\[
\gamma_\calB^G F^G([p_X]) = \gamma_\calB(F p_X),
\]
which is an epimorphism in $\Ho(\calM_\calB)^G$:
\[
(FQ_{\calA^G}X,\{\gamma_\calB(F\tau^X_g)\circ\gamma_\calB((\sigma_{g,Q_{\calA^G}X})^{-1})\})
\twoheadrightarrow
(FX,\{\gamma_\calB(F\phi_g)\circ\gamma_\calB((\sigma_{g,X})^{-1})\}).
\]
Indeed, this epimorphism follows from the commutative diagram in $\Ho(\calM_\calB)$:
\[
\xymatrix@C=2.5CM{
   \rho^\calB_g F Q_{\calA^G} X   \ar@{->>}[d]_-{\gamma_\calB(\rho^\calB_g F p_X)}
                                   \ar[r]^-{\gamma_\calB((\sigma_{g, Q_{\calA^G} X})^{-1})}
&  F \rho^\calA_g Q_{\calA^G} X   \ar@{->>}[d]^-{\gamma_\calB(F \rho^\calA_g p_X)}
                                      \ar[r]^-{\gamma_\calB(F \tau^X_g)}
&  F Q_{\calA^G} X                \ar@{->>}[d]^-{\gamma_\calB(F p_X)}             \\
   \rho^\calB_g F X               \ar[r]^-{\gamma_\calB((\sigma_{g, X})^{-1})}
&  F \rho^\calA_g X               \ar[r]^-{\gamma_\calB(F \phi_g)}
&  F X                                                        }
\]
Since $F$ is exact and $F(\calF_\calA \cap \calW_\calA) \subseteq \calW_\calB$ by assumption, the morphism $Fp_X$ is a weak equivalence in $\calB^G$.
Hence its image $\gamma_\calB(F p_X)$ is an isomorphism in $\Ho(\calM_\calB)^G$.
Thus the explicit value of the left-hand composite at $(X,\phi)$ is
\[
\gamma_\calB^G F^G \circ Q_{\calA^G}(X,\phi)
= (FQ_{\calA^G}X,\{\gamma_\calB(F\tau^X_g)\circ\gamma_\calB((\sigma_{g,Q_{\calA^G}X})^{-1})\}).
\]

We next compute the value of $(\gamma_\calB F)^G \circ Q_\calA^G$ at $(X,\phi)$.
By Lemma \ref{G-functor Q},
\[
Q_\calA^G(X,\phi) = \bigl(Q_\calA X,\; \{[Q_\calA \phi_g] \circ [\sigma^{Q_\calA}_{g,X}]^{-1}\}\bigr)
\]
in $\St_{\calC_\calA \cap \calW_\calA}(\calA)^G$.
Moreover, for each $g\in G$, we have the commutative diagram in $\St_{\calC_\calA\cap\calW_\calA}(\calA)$:
\[
\xymatrix@C=1.2CM{
   (\rho^\calA)_g Q_\calA X
\ar@{->>}[d]_-{[(\rho^\calA)_g q_X]}
\ar[r]^-{[\sigma^{Q_\calA}_{g, X}]^{-1}}
&  Q_\calA (\rho^\calA)_g X
\ar@{->>}[d]^-{[q_{(\rho^\calA)_g X}]}
\ar[r]^-{[Q_\calA \phi_g]}
&  Q_\calA X
\ar@{->>}[d]^-{[q_X]}             \\
   (\rho^\calA)_g X
\ar@{=}[r]
&  (\rho^\calA)_g X
\ar[r]^-{[\phi_g]}
&  X                                                        }
\]
This diagram yields an epimorphism in $\St_{\calC_\calA\cap\calW_\calA}(\calA)^G$:
\[
[q_X] : \bigl(Q_\calA X,\; \{[Q_\calA \phi_g] \circ [\sigma^{Q_\calA}_{g,X}]^{-1}\}\bigr)
        \twoheadrightarrow
        \bigl(X,\; \{[\phi_g]\}\bigr),
\]
where the underlying morphism $q_X : Q_\calA X \twoheadrightarrow X$ is a trivial fibration in $\calA$.

Applying $(\gamma_\calB F)^G$ to the epimorphism $[q_X]$ gives
\[
(\gamma_\calB F)^G([q_X]) = \gamma_\calB(F q_X),
\]
which is an epimorphism in $\Ho(\calM_\calB)^G$:
\[
\bigl(F Q_\calA X,\; \{\gamma_\calB(F Q_\calA \phi_g) \circ \gamma_\calB(F \sigma^{Q_\calA}_{g,X})^{-1} \circ \gamma_\calB(\sigma_{g,Q_\calA X})^{-1}\}\bigr)
\twoheadrightarrow
\bigl(FX,\; \{\gamma_\calB(F \phi_g) \circ \gamma_\calB(\sigma_{g,X})^{-1}\}\bigr).
\]
Indeed, the following commutative diagram in $\Ho(\calM_\calB)$ establishes the claim:
\[
\xymatrix@C=1.75CM{
   (\rho^\calB)_g F Q_\calA X
\ar@{->>}[d]_-{\gamma_\calB((\rho^\calB)_g F q_X)}
\ar[r]^-{\gamma_\calB(\sigma_{g, Q_\calA X})^{-1}}
&  F (\rho^\calA)_g Q_\calA X
\ar@{->>}[d]^-{\gamma_\calB(F (\rho^\calA)_g q_X)}
\ar[r]^-{\gamma_\calB(F \sigma^{Q_\calA}_{g, X})^{-1}}
&  F Q_\calA (\rho^\calA)_g X
\ar@{->>}[d]^-{\gamma_\calB(F q_{(\rho^\calA)_g X})}
\ar[r]^-{\gamma_\calB(F Q_\calA \phi_g)}
&  F Q_\calA X
\ar@{->>}[d]^-{\gamma_\calB(F q_X)}             \\
   (\rho^\calB)_g F X
\ar[r]^-{\gamma_\calB(\sigma_{g, X})^{-1}}
&  F (\rho^\calA)_g X
\ar@{=}[r]
&  F (\rho^\calA)_g X                                               \ar[r]^-{\gamma_\calB(F \phi_g)}
&  F X                                                        }
\]
Since $F$ is exact and $F(\calF_\calA \cap \calW_\calA) \subseteq \calW_\calB$ by assumption, the morphism $Fq_X$ is a weak equivalence.
Hence $\gamma_\calB(Fq_X)$ is an isomorphism in $\Ho(\calM_\calB)^G$.
Furthermore, the composite $(\gamma_\calB F)^G \circ Q_\calA^G$ has the explicit value
\[
(\gamma_\calB F)^G \circ Q_\calA^G(X,\phi)
= \bigl(F Q_\calA X,\; \{\gamma_\calB(F Q_\calA \phi_g)\circ\gamma_\calB(F \sigma^{Q_\calA}_{g,X})^{-1}\circ\gamma_\calB(\sigma_{g,Q_\calA X})^{-1}\}\bigr).
\]

Observe that
\[
(FX, \{\gamma_\calB(F \phi_g) \circ \gamma_\calB((\sigma_{g, X})^{-1})\})
=
(FX, \{\gamma_\calB(F \phi_g) \circ \gamma_\calB(\sigma_{g, X})^{-1}\}).
\]
Thus the composite
\[
\gamma_\calB(F q_X)^{-1} \circ \gamma_\calB(F p_X) :
\gamma_\calB^G F^G \circ Q_{\calA^G}(X, \phi) \to (\gamma_\calB F)^G \circ Q_\calA^G(X, \phi)
\]
is an isomorphism in $\Ho(\calM_\calB)^G$.
We claim that these isomorphisms, depending naturally on $(X,\phi)\in\calA^G$, assemble into a natural isomorphism
\[
\{\gamma_\calB(F q_X)^{-1} \circ \gamma_\calB(F p_X)\}_{(X, \phi) \in \calA^G} :
\gamma_\calB^G F^G \circ Q_{\calA^G} \to (\gamma_\calB F)^G \circ Q_\calA^G.
\]

\vspace{1ex}\textbf{Step 2. Prove the naturality of the isomorphism.}

For any morphism $f : (X, \phi) \to (Y, \psi)$ in $\calA^G$, we have the following commutative square in the stable category $\St_{\calC_\calA^G \cap \calW_\calA^G}(\calA^G)$:
\[
\xymatrix{
   (Q_{\calA^G} X, \tau^X)
\ar[d]_-{[Q_\calA^G f]}
\ar@{->>}[r]^-{[p_X]}
&  (X, \phi)
\ar[d]^{[f]}               \\
   (Q_{\calA^G} Y, \tau^Y)
\ar@{->>}[r]^-{[p_Y]}
&  (Y, \psi).                                            }
\]
Applying the functor $\gamma_\calB^G F^G$ to this square yields the commutative square in $\Ho(\calM_\calB)^G$:
\begin{equation*} \label{key square 1} \tag{\ref{thmB origin}.1}
\begin{gathered}
\xymatrix@C=1.5CM{
   \gamma_\calB^G F^G \circ Q_{\calA^G}(X, \phi)
\ar@{->>}[r]^-{\gamma_\calB(F p_X)}
\ar[d]_-{\gamma_\calB(F Q_{\calA^G} f)}
&  (FX, \{\gamma_\calB(F \phi_g) \circ \gamma_\calB((\sigma_{g, X})^{-1})\})
\ar[d]^-{\gamma_\calB(F f)}                   \\
   \gamma_\calB^G F^G \circ Q_{\calA^G}(Y, \psi)
\ar@{->>}[r]^-{\gamma_\calB(F p_Y)}
&  (FY, \{\gamma_\calB(F \psi_g) \circ \gamma_\calB((\sigma_{g, Y})^{-1})\})
}
\end{gathered}
\end{equation*}

Moreover, we have the following commutative square in $\St_{\calC_\calA \cap \calW_\calA}(\calA)^G$:
\[
\xymatrix{
   (Q_\calA X, \{[Q_\calA \phi_g] \circ [\sigma^{Q_\calA}_{g, X}]^{-1}\})
\ar[d]_{[Q_\calA f]}  \ar@{->>}[r]^-{[q_X]}
&  (X, \{[\phi_g]\})
\ar[d]^-{[f]}               \\
   (Q_\calA Y, \{[Q_\calA \psi_g] \circ [\sigma^{Q_\calA}_{g, Y}]^{-1}\})
\ar@{->>}[r]^-{[q_Y]}
&  (Y, \{[\psi_g]\})                                            }
\]
By Lemma \ref{fundamental Lemma}, $[Q_\calA f]$ is a morphism of $G$-equivariant objects.
Applying $(\gamma_\calB F)^G$ to this square yields the commutative square in $\Ho(\calM_\calB)^G$:
\begin{equation*} \label{key square 2} \tag{\ref{thmB origin}.2}
\begin{gathered}
\xymatrix@C=1.5CM{
   (\gamma_\calB F)^G \circ Q_\calA^G(X, \phi)
\ar@{->>}[r]^-{\gamma_\calB(F q_X)}
\ar[d]_-{\gamma_\calB(F Q_\calA f)}
&  (FX, \{\gamma_\calB(F \phi_g) \circ \gamma_\calB(\sigma_{g, X})^{-1}\})
\ar[d]^-{\gamma_\calB(F f)}                   \\
   (\gamma_\calB F)^G \circ Q_\calA^G(Y, \psi)
\ar@{->>}[r]^-{\gamma_\calB(F q_Y)}
&  (FY, \{\gamma_\calB(F \psi_g) \circ \gamma_\calB(\sigma_{g, Y})^{-1}\})
}
\end{gathered}
\end{equation*}

Finally, combining squares (\ref{key square 1}) and (\ref{key square 2}) yields the desired commutative diagram, which establishes the naturality of the isomorphisms constructed above.
\qed
\end{bfhpg}

\section{Illustrating the main theorems via PGF Hovey triples}
\label{Abelian model structures on Mod(RG) int}
\noindent
In this section, we illustrate the main results of the paper---Theorems \ref{retract (Q^G,W^G,R^G)} and \ref{thmB homotopy square}---in a concrete module-theoretic setting. More precisely, for any associative ring $R$, we focus on the cofibrantly generated hereditary Hovey triple
\[
\mathfrak{PGF}(R):=(\PGF(R),\PGF(R)^\bot,\Mod(R))
\]
associated with projectively coresolved Gorenstein flat (or PGF for short) $R$-modules, as established by \v{S}aroch and \v{S}\v{t}ov\'{\i}\v{c}ek \cite{2020Sto}.
In Subsection \ref{Abelian model structures on Mod(RG)}, we apply Theorem \ref{retract (Q^G,W^G,R^G)} to describe the trivial class $\PGF(RG)^\bot$ of the PGF Hovey triple
\[
\mathfrak{PGF}(RG)=(\PGF(RG),\PGF(RG)^\bot,\Mod(RG))
\]
in $\Mod(RG)$, the category of modules over the skew group algebra $RG$, and simultaneously obtain the comparison functor between the homotopy categories
\[
\Ho(\mathfrak{PGF}(RG)) \quad\text{and}\quad \Ho(\mathfrak{PGF}(R))^G,
\]
which has an explicit description as in \ref{con for Theta}; see Corollary \ref{cor of PGF}.
As illustrations of Theorem \ref{thmB homotopy square}, we turn in Subsections \ref{Homotopy squares via bimodule} and \ref{Homotopy squares via stable equivalence} to a two-ring setting involving $A$ and $B$.
In Subsection \ref{Homotopy squares via bimodule}, we use a Frobenius bimodule ${}_BM_A$ to induce a homotopy square involving PGF modules (see Corollary \ref{homotopy square for F-bimodules}).
To remove the restrictive assumptions of Corollary \ref{homotopy square for F-bimodules}, we then employ a stable equivalence of adjoint type, defined by a pair of bimodules, to recover the same homotopy square (see Corollary \ref{homotopy square for adjoint type}).

\begin{rmk} \label{rmk for G homo modules}
We emphasise that PGF modules are not the only way to illustrate Theorems \ref{retract (Q^G,W^G,R^G)} and \ref{thmB homotopy square}.
The same machinery works equally well for the other standard Gorenstein model structures, mutatis mutandis, as we now briefly explain.

For Gorenstein projective modules, the verification is entirely analogous to the PGF case, and in fact simpler, as it does not involve tensor products.
The only additional assumption needed is that the pair consisting of the class of Gorenstein projective modules and its right orthogonal forms a complete hereditary cotorsion pair cogenerated by a set---a hypothesis that is automatically satisfied in many settings of interest but remains unknown for arbitrary rings.
The Gorenstein injective case is formally dual to the Gorenstein projective one.
The Gorenstein flat case is completely parallel to the PGF case.

In summary, among these four Gorenstein model structures, we have chosen the PGF case as our primary example precisely because its verification process is the most representative, subsuming the key features and technical challenges that also appear in the other cases.
\end{rmk}

We continue with the definition of PGF modules.

\begin{ipg}
Following \cite{2020Sto}, an $R$-module $M$ is called \emph{projectively coresolved Gorenstein flat} (or PGF for short) if there exists an acyclic complex of projective $R$-modules
\[
\cdots \to P^{-1} \to P^0 \to P^1 \to \cdots
\]
with $M\cong \ker(P^0\to P^1)$, and such that the complex remains acyclic after applying $E\otimes_R-$ for every injective $\Rop$-module $E$. We denote by $\PGF(R)$ the class of PGF $R$-modules.
\end{ipg}

\begin{ipg}\label{exa of Goren}
For any ring $R$, \thmcite[4.9]{2020Sto} shows that $(\PGF(R),\PGF(R)^\bot)$ is a complete hereditary cotorsion pair in $\Mod(R)$, cogenerated by a set, with
\[
\PGF(R)\cap\PGF(R)^\bot=\Proj(R),
\]
the class of projective $R$-modules. Hence, by \thmcite[1.1]{HOWTO},
\[
\mathfrak{PGF}(R):=(\PGF(R),\PGF(R)^\bot,\Mod(R))
\]
is a cofibrantly generated hereditary Hovey triple in $\Mod(R)$.
\end{ipg}

\begin{setup}
In Subsection \ref{Abelian model structures on Mod(RG)}, we work with a ring $R$ and a finite group $G$ acting on $R$ via a group homomorphism $\rho : G \to \Aut(R)^{\op}$.
In Subsections \ref{Homotopy squares via bimodule} and \ref{Homotopy squares via stable equivalence}, we consider two rings $A$ and $B$ such that there exist group homomorphisms $\rho_1 : G \to \Aut(A)^{\op}$ and $\rho_2 : G \to \Aut(B)^{\op}$, and assume that $|G|$ is invertible in both $\Mod(A)$ and $\Mod(B)$.
Throughout, we identify $\Mod(RG)$ with the equivariant category $\Mod(R)^G$ via the canonical equivalence described in Example \ref{exa skew}.
\end{setup}

\subsection{PGF Hovey triple over skew group rings}
\label{Abelian model structures on Mod(RG)}
\noindent
In this subsection, to illustrate Theorem \ref{retract (Q^G,W^G,R^G)}, we show that, under certain conditions, $\PGF(RG)^\bot = (\PGF(R)^\bot)^G$, and simultaneously there is a comparison functor
from $\Ho(\mathfrak{PGF}(RG))$ to $\Ho(\mathfrak{PGF}(R))^G$,
which is a triangle equivalence up to retracts.

The following result will be used in the proof of Lemma \ref{PGF G-invariant}.

\begin{lem} \label{tensor product}
Let $M$ be an $\Rop$-module and $N$ an $R$-module.
Then there is an isomorphism $M^g\otimes_RN \cong M\otimes_R N^{g^{-1}}$, which is natural in both variables.
\end{lem}

\begin{prf*}
Consider the abelian group $M^g\otimes_R N$ together with the canonical $R$-biadditive map $f:M^g\times N\to M^g\otimes_R N$. Define a map
\[
h:M^g\times N \longrightarrow M\otimes_R N^{g^{-1}},\qquad (m,n)\longmapsto m\otimes n.
\]
It is straightforward to verify that $h$ is $R$-biadditive. By the universal property of the tensor product, there is a unique homomorphism
\[
\alpha:M^g\otimes_R N \longrightarrow M\otimes_R N^{g^{-1}}
\]
such that $\alpha(m\otimes n)=h(m,n)=m\otimes n$. This map is clearly an isomorphism, with inverse given by $m\otimes n\mapsto m\otimes n$ (using the natural identification of the two tensor products). Naturality in both $M$ and $N$ follows by a direct verification.
\end{prf*}

\begin{lem} \label{PGF G-invariant}
The class $\PGF(R)$ is $G$-invariant.
\end{lem}

\begin{prf*}
Let $M$ be a PGF $R$-module. We show that $M^g$ is again PGF.
By definition, there is an acyclic complex of projective $R$-modules
\[
P^\bullet:\cdots\to P^{-1}\to P^0\to P^1\to\cdots
\]
with $M\cong \ker(P^0\to P^1)$, which remains acyclic after applying $E\otimes_R-$ for every injective $\Rop$-module $E$.
Since $\Proj(R)$ is $G$-invariant, applying the exact functor $(-)^g$ yields an acyclic complex of projective $R$-modules
\[
(P^\bullet)^g:\cdots\to (P^{-1})^g\to (P^0)^g\to (P^1)^g\to\cdots
\]
with $M^g\cong \ker((P^0)^g\to(P^1)^g)$.
By Lemma \ref{tensor product}, we have an isomorphism of complexes
\[
E\otimes_R (P^\bullet)^g \cong E^{g^{-1}} \otimes_R P^\bullet.
\]
Since $E^{g^{-1}}$ is injective, the complex $E^{g^{-1}} \otimes_R P^\bullet$ is acyclic; hence so is $E\otimes_R(P^\bullet)^g$.
Therefore $M^g$ is PGF.
\end{prf*}

Next, we give a characterization of PGF $RG$-modules.

\begin{lem} \label{PGF equ}
Suppose that $|G|$ is invertible in $\Mod(R)$.
Then $\PGF(RG)=\PGF(R)^G$.
\end{lem}

\begin{prf*}
We first prove the inclusion $\PGF(R)^G\subseteq\PGF(RG)$.
Let $(M,\phi)\in\PGF(R)^G$. Then $M\in\PGF(R)$. By Lemma \ref{split cond}, $(M,\phi)$ is a direct summand of $(\Ind\circ U)(M,\phi)=\Ind(M)$. Since $\PGF(RG)$ is closed under direct summands, it suffices to show that $\Ind(M)\in\PGF(RG)$.
To this end, by definition, there is an acyclic complex of projective $R$-modules
\[
P^\bullet:\cdots\to P^{-1}\to P^0\to P^1\to\cdots
\]
with $M\cong\ker(P^0\to P^1)$, which remains acyclic after applying $E\otimes_R-$ for every injective $\Rop$-module $E$. Applying the exact functor $\Ind$, which preserves projectives (see \ref{frobenius pair (Ind,w)}), yields an acyclic complex of projective $RG$-modules
\[
\Ind(P^\bullet):\cdots\to \Ind(P^{-1})\to \Ind(P^0)\to \Ind(P^1)\to\cdots
\]
with $\Ind(M)\cong\ker(\Ind(P^0)\to\Ind(P^1))$.
Let $(E,\lambda)$ be an injective $(RG)^{\op}$-module. By Remark \ref{watts thm}, we have
\[
(E,\lambda)\otimes_{RG}\Ind(P^\bullet)
\cong (E,\lambda)\otimes_{RG}(RG\otimes_R P^\bullet)
\cong E\otimes_R P^\bullet.
\]
Since $E\otimes_R P^\bullet$ is acyclic, so is $(E,\lambda)\otimes_{RG}\Ind(P^\bullet)$. Hence, $\Ind(M)\in\PGF(RG)$, as desired.

For the reverse inclusion, take $(M,\phi)\in\PGF(RG)$. We show that $M\in\PGF(R)$.
Indeed, there exists an acyclic complex of projective $RG$-modules
\[
(P^\bullet,\mu^\bullet):
\cdots\to (P^{-1},\mu^{-1})\to (P^0,\mu^0)\to (P^1,\mu^1)\to\cdots
\]
with $(M,\phi)\cong\ker((P^0,\mu^0)\to(P^1,\mu^1))$, which remains acyclic after applying $(I,\lambda)\otimes_{RG}-$ for every injective $(RG)^{\op}$-module $(I,\lambda)$.
Applying the forgetful functor $U$ gives an acyclic complex of projective $R$-modules
\[
P^\bullet:\cdots\to P^{-1}\to P^0\to P^1\to\cdots
\]
with $M\cong\ker(P^0\to P^1)$.
Let $E$ be an injective $R^{\op}$-module. Since $(\Ind,U)$ is a Frobenius pair, $\Ind(E)$ is an injective $(RG)^{\op}$-module. Thus, $\Ind(E)\otimes_{RG}(P^\bullet,\mu^\bullet)$ is acyclic. But
\[
\Ind(E)\otimes_{RG}(P^\bullet,\mu^\bullet)
\cong (E\otimes_R RG)\otimes_{RG}(P^\bullet,\mu^\bullet)
\cong E\otimes_R P^\bullet.
\]
Therefore, $E\otimes_R P^\bullet$ is acyclic for every injective $R^{\op}$-module $E$, so $M\in\PGF(R)$.
\end{prf*}

By \ref{exa of Goren}, we know that
\[
\mathfrak{PGF}(RG)=(\PGF(RG),\PGF(RG)^\bot,\Mod(RG))
\]
is a cofibrantly generated hereditary Hovey triple in $\Mod(RG)$.
Suppose that $|G|$ is invertible in $\Mod(R)$.
Then, by Lemma \ref{PGF G-invariant}, Theorem \ref{retract (Q^G,W^G,R^G)} yields the cofibrantly generated hereditary Hovey triple
\[
(\PGF(R)^G,(\PGF(R)^\bot)^G,\Mod(R)^G)
\]
in $\Mod(RG)$.
Consequently, by Lemma \ref{PGF equ} and the identification $\Mod(RG)=\Mod(R)^G$, we obtain the following result.

\begin{cor} \label{cor of PGF}
Suppose that $|G|$ is invertible in $\Mod(R)$. Then
\[
\PGF(RG)^\bot = (\PGF(R)^\bot)^G.
\]
Moreover, there is a comparison functor
\[
\Theta:\Ho(\mathfrak{PGF}(RG)) \longrightarrow \Ho(\mathfrak{PGF}(R))^G,
\]
which is a triangle equivalence up to retracts.
\end{cor}

\subsection{Homotopy squares induced by a Frobenius bimodule}
\label{Homotopy squares via bimodule}
\noindent
In this subsection, as an illustration of Theorem \ref{thmB homotopy square}, we use a Frobenius bimodule with respect to rings $A$ and $B$ to induce homotopy squares involving PGF modules.

Recall from \cite{1992AF} that a $B$-$A$-bimodule $M$ is called a \emph{Frobenius bimodule} if both ${}_B M$ and $M_A$ are finitely generated projective and there is an $A$-$B$-bimodule isomorphism
\[
{}^*M:=\Hom_B(M,B)\cong \Hom_{A^{\mathrm{op}}}(M,A)=:M^*.
\]
We denote both ${}^*M$ and $M^*$ by $N$. Then, by \thmcite[2.1]{1999CJG},
\[
T_M:=M\otimes_A-:\Mod(A)\rightleftarrows\Mod(B):N\otimes_B-=:T_N
\]
is a Frobenius pair. Indeed, \thmcite[2.1]{1999CJG} shows that every Frobenius pair between module categories arises essentially from tensor functors associated with Frobenius bimodules. The functors $T_M$ and $T_N$ are exact and preserve projective, injective, and flat modules (see, for example, \cite[2.1]{2020HLGZ}).

\begin{exa}
Frobenius bimodules arise in many contexts; for instance, every Frobenius extension $A\to B$ gives a Frobenius bimodule ${}_B B_A$ (see \cite{KL99}). In particular, the skew group ring $RG$ is a Frobenius extension of $R$, and hence, $_{RG}(RG)_R$ is a Frobenius bimodule.
\end{exa}

Recall from Example \ref{exa of Goren} that $\Mod(A)$ and $\Mod(B)$ admit cofibrantly generated hereditary Hovey triples $\mathfrak{PGF}(A)$ and $\mathfrak{PGF}(B)$, respectively.
The following result shows that the adjunctions $(T_M,T_N)$ and $(T_N,T_M)$ are Quillen adjunctions between the corresponding abelian model structures.

\begin{lem} \label{TM Quillen adj}
Let $_BM_A$ be a Frobenius bimodule.
Then $T_M(\PGF(A)) \subseteq \PGF(B)$ and $T_M(\PGF(A)^\bot) \subseteq \PGF(B)^\bot$.
\end{lem}

\begin{prf*}
We first prove that $T_M(\PGF(A))\subseteq\PGF(B)$. Let $X\in\PGF(A)$. By definition, there is an acyclic complex of projective $A$-modules
\[
P^\bullet:\cdots\to P^{-1}\to P^0\to P^1\to\cdots
\]
with $X\cong \ker(P^0\to P^1)$, which remains acyclic after applying $E\otimes_A-$ for every injective $A^{\mathrm{op}}$-module $E$. Applying the exact functor $T_M$ yields an acyclic complex of projective $B$-modules
\[
T_M(P^\bullet):\cdots\to T_M(P^{-1})\to T_M(P^0)\to T_M(P^1)\to\cdots
\]
with $T_M(X)\cong \ker(T_M(P^0)\to T_M(P^1))$.
Now let $I$ be an injective $B^{\mathrm{op}}$-module. Since $(-\otimes_A N,-\otimes_B M)$ is an adjoint pair, $I\otimes_B M$ is injective as an $A^{\mathrm{op}}$-module. Hence,
\[
I\otimes_B T_M(P^\bullet)\cong (I\otimes_B M)\otimes_A P^\bullet
\]
is acyclic. Therefore, $T_M(X)\in\PGF(B)$.

Next, we prove that $T_M(\PGF(A)^\perp)\subseteq\PGF(B)^\perp$. Let $X\in\PGF(A)^\perp$. For any $Y\in\PGF(B)$, Lemma \ref{Ext(F,H)} gives an isomorphism
\[
\Ext_B^1(Y,T_M(X))\cong \Ext_A^1(T_N(Y),X).
\]
By the same argument as above, $T_N(\PGF(B))\subseteq\PGF(A)$. Thus, $\Ext_A^1(T_N(Y),X)=0$, and consequently, $\Ext_B^1(Y,T_M(X))=0$. Hence, $T_M(X)\in\PGF(B)^\perp$.
\end{prf*}

In the following, let ${}_BM_A$ be a Frobenius bimodule such that $M_A$ and $N_B$ are generators, and assume that either $T_M$ or $T_N$ is a $G$-functor.
Then by \thmcite[2.2]{2020HLGZ} one gets that both $T_M$ and $T_N$ are faithful.
Consequently, for every object $X\in\Mod(A)$, we have a short exact sequence
\[
0 \to X \xrightarrow{\eta_X} T_N T_M(X) \to \coker(\eta_X) \to 0
\]
in $\Mod(A)$. Similarly, for every $Y\in\Mod(B)$, there is a short exact sequence
\[
0 \to \ker(\varepsilon_Y) \to T_M T_N(Y) \xrightarrow{\varepsilon_Y} Y \to 0
\]
in $\Mod(B)$; see \ref{a key seq}.

Lemma \ref{TM Quillen adj} now ensures that the hypotheses of Theorem \ref{thmB homotopy square} are satisfied, yielding the following result.

\begin{cor} \label{homotopy square for F-bimodules}
Suppose that $\cok(\eta_X) \in \PGF(A)^\bot$ for every $X \in \PGF(A)$, and $\kernel(\varepsilon_Y) \in \PGF(B)^\bot$ for every $Y \in \Mod(B)$.
Then the diagram
\[
\xymatrix@C=1.2CM{
   \Ho(\mathfrak{PGF}(AG))   \ar[d]_{\Ho(\gamma_\calA^G)}^{ }  \ar[r]^{\mathbb{L}(T_M^G)}_{\simeq}
&  \Ho(\mathfrak{PGF}(BG))   \ar[d]^{\Ho(\gamma_\calB^G)}_{ }             \\
   \Ho(\mathfrak{PGF}(A))^G                              \ar[r]^{\mathbb{L}(T_M)^G}_{\simeq}
&  \Ho(\mathfrak{PGF}(B))^G                                          }
\]
commutes up to natural isomorphism, where both the total left derived functor $\mathbb{L}(T_M^G)$ and the equivariant functor $\mathbb{L}(T_M)^G$ are triangle equivalences.
\end{cor}

\subsection{Homotopy squares induced by a stable equivalence of adjoint type}
\label{Homotopy squares via stable equivalence}
\noindent
In this subsection, to remove the restrictive assumptions of Corollary \ref{homotopy square for F-bimodules}, we employ a stable equivalence of adjoint type, defined by a pair of bimodules, to recover the same homotopy square.

The notion of stable equivalences of Morita type was introduced by Brou\'{e} \cite{1994Bro} for finite dimensional algebras, particularly those over algebraically closed fields. For general rings, it is more natural to adopt the following generalization, taken from \dfncite[5.3.15]{2014Zim}.

\begin{dfn} \label{Stable equivalences of Morita type}
Two bimodules ${}_BM_A$ and ${}_AN_B$ are said to define a \emph{stable equivalence of Morita type in the generalized sense} between rings $A$ and $B$ if the following conditions hold:
\begin{prt}
\item both $M$ and $N$ are finitely generated projective as left and as right modules;
\item there are isomorphisms of bimodules
\[
N\otimes_B M \cong A\oplus P
\quad\text{and}\quad
M\otimes_A N \cong B\oplus Q,
\]
where $P$ is a finitely generated flat $A$-$A$-bimodule such that $P\otimes_A X$ is projective for every $A$-module $X$, and $Q$ is a finitely generated flat $B$-$B$-bimodule such that $Q\otimes_B Y$ is projective for every $B$-module $Y$.
\end{prt}
\end{dfn}

\begin{rmk} \label{rmk of Morita type}
In \dfncite[5.3.15]{2014Zim}, Zimmermann works over a commutative ring $K$, with $A$ and $B$ as $K$-algebras. Here we restrict to the special case $K=\mathbb Z$. In condition (a), we additionally require that $M$ and $N$ be finitely generated on both sides; this hypothesis will be needed in Lemma \ref{op adjoint pair}.

Note that both ${}_B M$ and ${}_A N$ are projective. Hence, $N\otimes_B M$ is a projective $A$-module, so its direct summand $P$ is also projective as an $A$-module. Similarly, one obtains that $P_A$, ${}_B Q$, and $Q_B$ are projective.
\end{rmk}

For finite dimensional algebras, the following results are well known; see, for example, \thmcite[4.1]{2002Xi} and \lemcite[2.2]{2003Liu}. The same statements remain valid for general rings. As in Subsection \ref{Homotopy squares via bimodule}, we write $T_M:=M\otimes_A-$ and $T_N:=N\otimes_B-$ and let $T_P:=P\otimes_A-$, and $T_Q:=Q\otimes_B-$.

\begin{lem} \label{Stable equivalences resules}
Let ${}_BM_A$ and ${}_AN_B$ define a stable equivalence of Morita type in the generalized sense between rings $A$ and $B$.
Then the following statements hold.
\begin{prt}
\item There exist natural isomorphisms
$T_N\circ T_M \cong \Id_{\Mod(A)}\oplus T_P$
and
$T_M\circ T_N \cong \Id_{\Mod(B)}\oplus T_Q$.
\item Each of the bimodules ${}_BM_A$ and ${}_AN_B$ is a projective generator on both sides.
\item The functors $T_M$ and $T_N$ are faithful.
\end{prt}
\end{lem}

\begin{prf*}
(a) It suffices to prove the first isomorphism; the second follows by the same argument. We compute
\begin{align*}
(T_N\circ T_M)(-)
&= N\otimes_B(M\otimes_A-) \\
&\cong (N\otimes_B M)\otimes_A- \\
&\cong (A\oplus P)\otimes_A- \\
&\cong (A\otimes_A-)\oplus(P\otimes_A-) \\
&\cong \Id_{\Mod(A)}\oplus T_P(-).
\end{align*}

(b) We prove that ${}_BM$ is a generator; the other cases are analogous. Since ${}_BM$ is projective, it suffices to show that for every nonzero $B$-module $X$, there is a nonzero $B$-homomorphism $M\to X$, i.e., $\Hom_B(M,X)\neq 0$. Indeed, by the isomorphism $M\otimes_A N\cong B\oplus Q$, where $Q$ is a finitely generated flat $B$-$B$-bimodule, there is a split epimorphism $p:M\otimes_A N\to B$. Applying $\Hom_B(-,X)$ gives a monomorphism
\[
p^*:\Hom_B(B,X)\hookrightarrow \Hom_B(M\otimes_A N,X).
\]
If $\Hom_B(M,X)=0$, then $\Hom_A(N,\Hom_B(M,X))=0$. Using the adjunction
\[
\Hom_B(M\otimes_A N,X)\cong \Hom_A(N,\Hom_B(M,X)),
\]
we obtain $\Hom_B(M\otimes_A N,X)=0$. Consequently, $\Hom_B(B,X)=0$, which forces $X\cong \Hom_B(B,X)=0$, a contradiction.

(c) To prove that $T_M$ is faithful, let $f:X\to Y$ be a morphism of $A$-modules such that $T_M(f)=0$. We show that $f=0$; the proof for $T_N$ is analogous.

Since $M_A$ is a projective generator by (b), there exist an index set $I$ and an epimorphism $\varphi:\bigoplus_{i\in I}M\to A$ of $A^{\mathrm{op}}$-modules. Tensoring with $X$ and $Y$ yields epimorphisms
\[
\varphi_X:\bigoplus_{i\in I}(M\otimes_A X)\to A\otimes_A X\cong X
\quad\text{and}\quad
\varphi_Y:\bigoplus_{i\in I}(M\otimes_A Y)\to A\otimes_A Y\cong Y.
\]
Consider the commutative diagram
\[
\xymatrix@C=1.75cm{
   \bigoplus_{i\in I}(M\otimes_A X)
\ar[d]_-{\varphi_X}
\ar[r]^-{\oplus(M\otimes f)}
&  \bigoplus_{i\in I}(M\otimes_A Y)
\ar[d]^-{\varphi_Y}  \\
   X
\ar[r]^-{f}
& Y
}
\]
Since $M\otimes f=0$, the top horizontal map is zero, so $f\circ\varphi_X=0$. As $\varphi_X$ is surjective, it follows that $f=0$. Therefore $T_M$ is faithful.
\end{prf*}

Suppose that ${}_BM_A$ and ${}_AN_B$ define a stable equivalence of Morita type between algebras $A$ and $B$. Following Xi \dfncite[3.2]{2008Xi}, we say that $M$ and $N$ define a \emph{stable equivalence of adjoint type} if both $(T_M,T_N)$ and $(T_N,T_M)$ are adjoint pairs. Examples of such equivalences can be found in \cite{2008Xi}, \cite{LX2007}, and \cite{DM2007}. Motivated by this, we introduce the following generalization for arbitrary rings.

\begin{dfn} \label{Stable equivalences of adjoints}
Let ${}_BM_A$ and ${}_AN_B$ define a stable equivalence of Morita type in the generalized sense between rings $A$ and $B$.
This equivalence is said to be of \emph{adjoint type in the generalized sense} if both
\[
(T_M,T_N) \quad\text{and}\quad (T_N,T_M)
\]
are adjoint pairs.
\end{dfn}

We now give an example of a stable equivalence of adjoint type in the generalized sense.

\begin{exa} \label{examples of adjoint type}
Suppose that the rings $A$ and $B$ are Morita equivalent, i.e., their module categories $\Mod(A)$ and $\Mod(B)$ are equivalent. By Morita theory (see \cite{1992AF}), there exist bimodules ${}_AN_B$ and ${}_BM_A$ such that:
\begin{eqc}
\item[$\bullet$] both $M$ and $N$ are finitely generated projective as left and as right modules;
\item[$\bullet$] the functors $M\otimes_A-$ and $N\otimes_B-$ are inverse equivalences; and
\item[$\bullet$] $M\otimes_A N\cong B$ and $N\otimes_B M\cong A$.
\end{eqc}
If we take $P=0=Q$ in Definition \ref{Stable equivalences of Morita type}, then ${}_AN_B$ and ${}_BM_A$ define a stable equivalence of Morita type in the generalized sense between $A$ and $B$. Moreover, both $(M\otimes_A-,N\otimes_B-)$ and $(N\otimes_B-,M\otimes_A-)$ are adjoint pairs. Thus, ${}_AN_B$ and ${}_BM_A$ define a stable equivalence of adjoint type in the generalized sense between $A$ and $B$.
\end{exa}

By Lemma \ref{Stable equivalences resules}(c), both $T_M$ and $T_N$ are faithful.
Hence, by \ref{a key seq}, for every $X\in\Mod(A)$ there is a short exact sequence
\[
0 \to X \xrightarrow{\eta_X} T_N \circ T_M(X) \to \cok(\eta_X) \to 0,
\]
and, for every $Y\in\Mod(B)$, there is a short exact sequence
\[
0 \to \kernel(\varepsilon_Y) \to T_M \circ T_N(Y) \xrightarrow{\varepsilon_Y} Y \to 0.
\]

We first record the following elementary observation on projective stable categories.

\begin{lem}\label{stable iso induces Ext iso}
Let $R$ be a ring and let $f : X \to Y$ be a morphism in $\Mod(R)$.
If $[f]$ is an isomorphism in the projective stable category $\St_{\Proj(R)}(\Mod(R)):=\Mod(R)/\Proj(R)$, then, for every $Z \in \Mod(R)$ and every integer $n \geq 1$, the induced homomorphism
\[
\Ext_R^n(f, Z) : \Ext_R^n(Y, Z) \to \Ext_R^n(X, Z)
\]
is an isomorphism.
\end{lem}

\begin{prf*}
Since $[f]$ is an isomorphism in $\St_{\Proj(R)}(\Mod(R))$, there exists a morphism $g : Y \to X$ such that $[g \circ f] = [\Id_X]$ and $[f \circ g] = [\Id_Y]$.
Hence both $g \circ f - \Id_X$ and $f \circ g - \Id_Y$ factor through projective $R$-modules.

We first note that if a morphism $h : U \to V$ factors through a projective $R$-module, then $\Ext_R^n(h,Z)=0$ for every $Z \in \Mod(R)$ and every $n \geq 1$.
Indeed, if $h = b \circ a$ for some factorization $U \xrightarrow{a} P \xrightarrow{b} V$ with $P \in \Proj(R)$, then, since $\Ext_R^n(-, Z)$ is contravariant in the first variable,
\[
\Ext_R^n(h, Z)
=
\Ext_R^n(a, Z) \circ \Ext_R^n(b, Z).
\]
This homomorphism factors through $\Ext_R^n(P, Z) = 0$, and hence $\Ext_R^n(h, Z) = 0$.

Applying this observation to $g \circ f - \Id_X$, we obtain $\Ext_R^n(g \circ f - \Id_X, Z) = 0$.
By the additivity and contravariance of $\Ext_R^n(-, Z)$, $\Ext_R^n(g \circ f, Z) = \Ext_R^n(f, Z) \circ \Ext_R^n(g, Z)$,
and therefore
\[
\Ext_R^n(f, Z) \circ \Ext_R^n(g, Z)
=
\Id_{\Ext_R^n(X, Z)}.
\]
Similarly, $\Ext_R^n(g, Z) \circ \Ext_R^n(f, Z) = \Id_{\Ext_R^n(Y, Z)}$.
Thus $\Ext_R^n(g, Z)$ is the inverse of $\Ext_R^n(f, Z)$, and the result follows.
\end{prf*}

We can now describe the homological properties of the unit and counit.

\begin{lem} \label{conditions}
Let $_BM_A$ and $_AN_B$ define a stable equivalence of adjoint type in the generalized sense between rings $A$ and $B$.
Then the following statements hold.
\begin{prt}
\item
For every $X\in\Mod(A)$, the $A$-module $\cok(\eta_X)$ has projective dimension at most one.
\item
For every $Y\in\Mod(B)$, the $B$-module $\kernel(\varepsilon_Y)$ is projective.
\end{prt}
\end{lem}

\begin{prf*}
(a) Let $X\in\Mod(A)$.
Consider the short exact sequence
\[
0  \to  X  \overset{\eta_X}  \longrightarrow  T_N\circ T_M(X)  \to  C_X  \to  0,
\]
where $C_X := \cok(\eta_X)$.
We show that $C_X$ has projective dimension at most one.
Take a short exact sequence $0 \to L\to P\to C_X \to 0$ with $P$ projective.
It suffices to show that $L$ is projective.
For this, let $Z$ be any $A$-module.
By dimension shifting, we have $\Ext_A^1(L, Z) \cong \Ext_A^2(C_X, Z)$.
Thus, it is enough to prove that $\Ext_A^2(C_X, Z) = 0$ for every $Z\in\Mod(A)$.

Applying the functor $\Hom_A(-,Z)$, we obtain the exact sequence
\[
\Ext_A^1(T_N \circ T_M(X),Z) \xrightarrow{\Ext_A^1(\eta_X,Z)} \Ext_A^1(X,Z) \longrightarrow \Ext_A^2(C_X,Z)
\]
\[
\longrightarrow \Ext_A^2(T_N \circ T_M(X),Z) \xrightarrow{\Ext_A^2(\eta_X,Z)} \Ext_A^2(X,Z).
\]
Therefore, to prove  $\Ext_A^2(C_X, Z) = 0$, it suffices to show that
\[
\Ext_A^n(\eta_X, Z) : \Ext_A^n(T_N \circ T_M(X), Z) \to \Ext_A^n(X, Z)
\]
is an isomorphism for $n = 1, 2$.
By Lemma~\ref{stable iso induces Ext iso}, it remains to show that $[\eta_X] : X \to T_N \circ T_M(X)$ is an isomorphism in $\St_{\Proj(A)}(\Mod(A))$.

Indeed, since $T_M$ and $T_N$ preserve projective modules, they induce functors
\[
T_M : \St_{\Proj(A)}(\Mod(A)) \rightleftarrows \St_{\Proj(B)}(\Mod(B)) : T_N.
\]
By Lemma \ref{Stable equivalences resules}(a), there are natural isomorphisms
\[
T_N \circ T_M \cong \Id_{\Mod(A)} \oplus T_P
\quad \text{and} \quad
T_M \circ T_N \cong \Id_{\Mod(B)} \oplus T_Q.
\]
Since $T_P(X')$ and $T_Q(Y')$ are projective for all $X'\in\Mod(A)$ and $Y'\in\Mod(B)$, these isomorphisms become
\[
T_N \circ T_M \cong \Id_{\St_{\Proj(A)}(\Mod(A))}
\quad \text{and} \quad
T_M \circ T_N \cong \Id_{\St_{\Proj(B)}(\Mod(B))}
\]
on the projective stable categories.
Hence the induced functors $T_M$ and $T_N$ are mutually quasi-inverse equivalences.
In particular, the induced functor $T_M$ is fully faithful.

It remains to relate this equivalence to the given unit $\eta$.
Since the stable equivalence is of adjoint type, $(T_M, T_N)$ is an adjoint pair on the module categories. Moreover, this adjunction descends to the projective stable categories.
To see this, consider the adjunction isomorphism
\[
\Hom_B(T_M(U), V) \cong \Hom_A(U, T_N(V)).
\]
If a morphism $T_M(U) \to V$ factors through a projective $B$-module, then its adjoint factors through a projective $A$-module, since $T_N$ preserves projective modules.
Conversely, if a morphism $U \to T_N(V)$ factors through a projective $A$-module, then its adjoint factors through a projective $B$-module, since $T_M$ preserves projective modules.
Thus the above adjunction isomorphism induces an adjunction
\[
T_M: \St_{\Proj(A)}(\Mod(A)) \rightleftarrows \St_{\Proj(B)}(\Mod(B)) : T_N,
\]
whose unit is precisely $[\eta]$.

Since the left adjoint $T_M$ is fully faithful, the unit of this adjunction is a natural isomorphism.
Therefore $[\eta_X] : X \to T_N \circ T_M(X)$ is an isomorphism in $\St_{\Proj(A)}(\Mod(A))$, as desired.

The proof of (b) is analogous to that of (a), and is even simpler, since it suffices to show that $\Ext_B^1(\ker(\varepsilon_Y), Z) = 0$ for every $Z \in \Mod(B)$.
\end{prf*}

For any ring $R$, every $R$-module of finite projective dimension belongs to $\PGF(R)^\bot$.
Therefore Lemma~ \ref{conditions} immediately yields the following consequence.

\begin{lem} \label{condition of adjoint type}
Let $_BM_A$ and $_AN_B$ define a stable equivalence of adjoint type in the generalized sense between rings $A$ and $B$.
Then $\cok(\eta_X) \in \PGF(A)^\bot$ for each $X \in \Mod(A)$ and $\kernel(\varepsilon_Y) \in \PGF(B)^\bot$ for each $Y \in \Mod(B)$.
\end{lem}

In the case of finite dimensional algebras, Xi \lemcite[3.3]{2008Xi} proved the following result.
Since both $M$ and $N$ are finitely generated projective on both sides, the same result remains valid in our context.

\begin{lem} \label{op adjoint pair}
Let $_BM_A$ and $_AN_B$ define a stable equivalence of adjoint type in the generalized sense between rings $A$ and $B$.
Then
\[
(- \otimes_A N, - \otimes_B M) \quad\text{and}\quad (- \otimes_B M, - \otimes_A N)
\]
are adjoint pairs between $\Mod(A^{\op})$ and $\Mod(B^{\op})$.
\end{lem}

\begin{prf*}
We only prove that $(- \otimes_A N, - \otimes_B M)$ is an adjoint pair.
The proof for the other adjoint pair is similar.
Recall that $(-\otimes_A N, \Hom_{\Bop}(N, -))$ is an adjoint pair.
Hence, it suffices to show that
\[
-\otimes_B M \cong \Hom_{\Bop}(N, -).
\]
For this, we need the following two facts:
\begin{rqm}
\item
$_BM_A \cong \Hom_{\Bop}(\Hom_B({_BM_A} , {_BB_B}), {_BB_B})$, and
\item
$- \otimes_B \Hom_{\Bop}({_AN_B} , {_BB_B}) \cong \Hom_{\Bop}({_AN_B} , -)$.
\end{rqm}
Now, since $(M \otimes_A -, N \otimes_B -)$ is an adjoint pair, we have
\[
\Hom_B({_BM_A}, B) \cong \Hom_B(M \otimes_A A, B) \cong \Hom_A(A, N \otimes_B B) \cong {_AN_B}.
\]
Hence, by (1) we obtain an isomorphism $_BM_A \cong \Hom_{\Bop}({_AN_B}, {_BB_B})$.
Then (2) yields the desired natural isomorphism.

It remains to verify (1) and (2).
Because $_BM$ is finitely generated projective, (1) holds by \prpcite[20.17]{1992AF}.
Since $N_B$ is finitely generated projective, (2) follows from \exercite[20.12]{1992AF}.
\end{prf*}

Recall again from Example \ref{exa of Goren} that $\Mod(A)$ and $\Mod(B)$ admit cofibrantly generated hereditary Hovey triples $\mathfrak{PGF}(A)$ and $\mathfrak{PGF}(B)$, respectively.
The following result implies that the adjunctions $(T_M, T_N)$ and $(T_N, T_M)$ are Quillen adjunctions between the corresponding abelian model structures on $\Mod(A)$ and $\Mod(B)$.

\begin{lem} \label{T_M preserve adjoint type}
Let $_BM_A$ and $_AN_B$ define a stable equivalence of adjoint type in the generalized sense between rings $A$ and $B$.
Then $T_M(\PGF(A)) \subseteq \PGF(B)$ and $T_M(\PGF(A)^\bot) \subseteq \PGF(B)^\bot$.
\end{lem}

\begin{prf*}
For the inclusion $T_M(\PGF(A)) \subseteq \PGF(B)$, take $X \in \PGF(A)$.
We need to prove that $T_M(X) \in \PGF(B)$.
There exists an exact sequence
\[
P^\bullet: \cdots \to P^{-1} \to P^0 \to P^1 \to \cdots
\]
of projective $A$-modules with $X \cong \kernel(P^0 \to P^1)$ such that it remains exact after applying $E\otimes_A-$ for every injective $A^{\op}$-module $E$.
Applying $T_M$ to $P^\bullet$ yields an exact sequence of projective $B$-modules
\[
T_M(P^\bullet): \cdots \to T_M(P^{-1}) \to T_M(P^0) \to T_M(P^1) \to \cdots
\]
with $T_M(X) \cong \kernel(T_M(P^0) \to T_M(P^1))$.
Let $I$ be an injective $B^{\op}$-module.
By Lemma \ref{op adjoint pair}, $(- \otimes_A N, - \otimes_B M)$ is an adjoint pair.
Then $I\otimes_BM$ is an injective $A^{\op}$-module as $\Hom_{A^{\op}}(-, I \otimes_B M) \cong \Hom_{B^{\op}}(- \otimes_A N, I)$.
Consequently, the complex
\[
I \otimes_B T_M(P^\bullet) \cong (I \otimes_B M) \otimes_A P^\bullet
\]
is acyclic.
Hence $T_M(X)$ is a PGF $B$-module.

Now take $X \in \PGF(A)^\bot$.
We show that $T_M(X) \in \PGF(B)^\bot$.
For any $Y \in \PGF(B)$, Lemma \ref{Ext(F,H)} provides an isomorphism
\[
\Ext_B^1(Y, T_M(X)) \cong \Ext_A^1(T_N(Y), X).
\]
As in the case of $T_M$, one has $T_N(\PGF(B)) \subseteq \PGF(A)$.
Therefore $\Ext_A^1(T_N(Y), X) = 0$, and consequently $\Ext_B^1(Y, T_M(X)) = 0$, as desired.
\end{prf*}

In the following, let $_BM_A$ and $_AN_B$ define a stable equivalence of adjoint type in the generalized sense between $A$ and $B$, and assume that either $T_M$ or $T_N$ is a $G$-functor. The conditions required by Theorem \ref{thmB homotopy square} are guaranteed by Lemmas \ref{T_M preserve adjoint type}, \ref{Stable equivalences resules}(c), and \ref{condition of adjoint type}. Consequently, we obtain the following result.

\begin{cor} \label{homotopy square for adjoint type}
The diagram
\[
\xymatrix@C=1.2CM{
   \Ho(\mathfrak{PGF}(AG))   \ar[d]_{\Ho(\gamma_\calA^G)}^{ }  \ar[r]^{\mathbb{L}(T_M^G)}_{\simeq}
&  \Ho(\mathfrak{PGF}(BG))   \ar[d]^{\Ho(\gamma_\calB^G)}_{ }             \\
   \Ho(\mathfrak{PGF}(A))^G                              \ar[r]^{\mathbb{L}(T_M)^G}_{\simeq}
&  \Ho(\mathfrak{PGF}(B))^G                                          }
\]
commutes up to natural isomorphism, where both the total left derived functor $\mathbb{L}(T_M^G)$ and the equivariant functor $\mathbb{L}(T_M)^G$ are triangle equivalences.
\end{cor}

\appendix
\section*{Appendix: Proofs of Lemmas \ref{admiss for Ho(M)} and \ref{admiss for St(C cap F)}}
\label{appendix}
\stepcounter{section}
\noindent
In this section, we give detailed proofs for Lemmas \ref{admiss for Ho(M)} and \ref{admiss for St(C cap F)}. The following result, taken from \lemcite[6.26]{2025GillBook} and also called the Connecting Morphism, will help us construct the natural isomorphism in Lemma \ref{admiss for Ho(M)}. {\sf Throughout the section}, let $\calM = (\calC, \calW, \calF)$ be a hereditary Hovey triple in $\calA$ with $\gamma_{\calA}: \calA \to \Ho(\calM)$ the localization functor.

\begin{lem} \label{The Connecting Mor}
Suppose that $0\to A \to B \to C\to 0$ and $0\to A \to W \to D\to 0$ are short exact sequences in $\calA$ with $W \in \calW$.
Construct a pushout diagram
\[
\xymatrix{
     A
\ar@{>->}[d]
\ar@{>->}[r]
&    B
\ar@{>->}[d]
\ar@{->>}[r]
&    C
\ar@{=}[d]   \\
     W
\ar@{->>}[d]
\ar@{>->}[r]
&    P
\ar@{->>}[d]^(0.45){v}
\ar@{->>}[r]^-{w}
&    C        \\
     D
\ar@{=}[r]
&  D
}
\]
Independently of the choices made when forming the pushout, the composition
\[
\gamma_{\calA}(v) \circ (\gamma_{\calA}(w))^{-1} \in \Hom_{\Ho(\calM)}(C, D)
\]
is well defined.
Moreover, the construction is natural in the following sense: for each pair of morphisms of short exact sequences, as shown below, with the same $f_A$ and $W' \in \calW$,
\[
\xymatrix{
   A
\ar[d]_-{f_A}
\ar@{>->}[r]
&  B
\ar[d]^-{f_B}
\ar@{->>}[r]
&  C \ar[d]^-{f_C}   \\
   A'
\ar@{>->}[r]
&  B'
\ar@{->>}[r]
&  C'               }
\quad \raisebox{-3.8ex}{\text{and}} \quad
\xymatrix{
   A
\ar[d]_-{f_A}
\ar@{>->}[r]
&  W
\ar[d]^-{f_W}
\ar@{->>}[r]
&  D
\ar[d]^-{f_D}   \\
A'
\ar@{>->}[r]
&  W'
\ar@{->>}[r]
&  D'               }
\]
the following square in $\Ho(\calM)$ commutes:
\[
\xymatrix@C=3CM{
  C
\ar[d]_-{\gamma_{\calA}(f_C)}
\ar[r]^-{\gamma_{\calA}(v) \circ (\gamma_{\calA}(w))^{-1}}
& D
\ar[d]^-{\gamma_{\calA}(f_D)}  \\
  C'
\ar[r]^-{\gamma_{\calA}(v') \circ (\gamma_{\calA}(w'))^{-1}}
& D'                      }
\]
\end{lem}

\begin{bfhpg}[Proof of Lemma \ref{admiss for Ho(M)}]\label{detailed proof}
We will give a $G$-action on the homotopy category $\Ho(\calM)$, and then show that it is admissible. The proof is divided into two steps.

\vspace{1ex}\textbf{Step 1}. Construct the $G$-action $(\widehat{\rho}, \widehat{\theta})$.

For each $g \in G$, we will use the Localization Theorem \ref{loc thm} to construct the functor $\widehat{\rho}_g : \Ho(\calM) \to \Ho(\calM)$.
Consider the diagram
\begin{equation*} \label{gamma diagram} \tag{\ref{detailed proof}.1}
\begin{gathered}
\xymatrix{
  \calA
\ar[d]_-{\rho_g}
\ar[r]^-{\gamma_{\calA}}
& \Ho(\calM)
\ar@{-->}[d]^-{\widehat{\rho}_g}  \\
  \calA
\ar[r]^-{\gamma_{\calA}}
& \Ho(\calM)    }
\end{gathered}
\end{equation*}
where $\gamma_{\calA}$ is the localization functor.
We have to show that the functor $\rho_g$ preserves weak equivalences.
Once this is done, there exists a unique functor $\widehat{\rho}_g$ such that the diagram commutes.
But this is obvious, because we already know that $\rho_g$ is exact, and by Lemma \ref{(Q,W,R) G-invariant}, the inclusions $\rho_g (\calC) \subseteq \calC$, $\rho_g (\calW) \subseteq \calW$, and $\rho_g (\calF) \subseteq \calF$ all hold.

Next we construct the natural isomorphism $\widehat{\theta}_{g,h}: \widehat{\rho}_g \circ \widehat{\rho}_h \to \widehat{\rho}_{hg}$.
Define $\widehat{\theta}_{g,h} \circ \gamma_{\calA}$ to be the composition
\[
\xymatrix@C=1CM{
   \widehat{\rho}_g \circ \widehat{\rho}_h \circ \gamma_{\calA}
=  \widehat{\rho}_g \circ \gamma_{\calA} \circ \rho_h
=  \gamma_{\calA} \circ \rho_g \circ \rho_h
\ar[r]^-{\gamma_{\calA} \circ \theta_{g,h}}_-{\cong}
&  \gamma_{\calA} \circ \rho_{hg}
=  \widehat{\rho}_{hg} \circ \gamma_{\calA},
}
\]
which is a natural isomorphism.
Then, by Lemma \ref{natural iso}, we see that $\widehat{\theta}_{g,h}$ is also a natural isomorphism.
In particular, for each $X \in \Ho(\calM)$, we have
\[
   \widehat{\theta}_{g,h,X}
=  \widehat{\theta}_{g,h,\gamma_{\calA}(X)}
=  \gamma_{\calA}(\theta_{g,h,X})
=  [RQ \theta_{g,h,X}] \in \Hom_{\Ho(\calM)}(\rho_g \rho_h X, \rho_{hg} X).
\]

Now, using the 2-cocycle condition for $(\rho, \theta)$, one can directly verify that the 2-cocycle condition for $(\widehat{\rho}, \widehat{\theta})$ holds.
This implies that $(\widehat{\rho}, \widehat{\theta})$ is a $G$-action on $\Ho(\calM)$.

\vspace{1ex}\textbf{Step 2}. Prove that $(\widehat{\rho}, \widehat{\theta})$ is admissible.

For each $g \in G$, we first show that there exists a natural isomorphism
\[
\varphi_g : \mathsf{\Sigma} \circ \widehat{\rho}_g \to \widehat{\rho}_g \circ \mathsf{\Sigma}.
\]
By Lemma \ref{natural iso}, it suffices to show that there exists a natural isomorphism
\[
\mathsf{\Sigma} \circ \widehat{\rho}_g \circ \gamma_{\calA} \cong \widehat{\rho}_g \circ \mathsf{\Sigma} \circ \gamma_{\calA}.
\]

For each $X \in \calA$, note that $(\calC, \calW \cap \calF)$ is a complete cotorsion pair in $\calA$.
Hence we may choose short exact sequences
\[
\rho_g X \rightarrowtail W_{\rho_g X} \twoheadrightarrow \mathsf{\Sigma} \rho_g X
\quad \text{and} \quad
\rho_g X \rightarrowtail \rho_g W_X \twoheadrightarrow \rho_g \mathsf{\Sigma} X,
\]
where $W_{\rho_g X}, \rho_g W_X \in \calW \cap \calF$ and $\mathsf{\Sigma} \rho_g X, \rho_g \mathsf{\Sigma} X \in \calC$.
Let $P_X$ be the pushout of $W_{\rho_g X}$ and $\rho_g W_X$ over $\rho_g X$, and denote by
\[
w_X : P_X \twoheadrightarrow \mathsf{\Sigma} \rho_g X
\quad \text{and} \quad
v_X : P_X \twoheadrightarrow \rho_g \mathsf{\Sigma} X
\]
the induced cokernel maps.
Both $w_X$ and $v_X$ are weak equivalences.
By Lemma \ref{The Connecting Mor}, we obtain the connecting morphism
\[
\gamma_{\calA}(v_X) \circ (\gamma_{\calA}(w_X))^{-1} \in \Hom_{\Ho(\calM)}(\mathsf{\Sigma} \rho_g X, \rho_g \mathsf{\Sigma} X),
\]
which is therefore an isomorphism in $\Ho(\calM)$.
The dual of Lemma \ref{fundamental Lemma} implies that the constructions of these short exact sequences are unique up to $\omega$-homotopy, where $\omega = \calC \cap \calW \cap \calF$.
Since $\gamma_{\calA}$ identifies homotopic morphisms, this connecting morphism does not depend on the choices of the sequences.
We claim that these isomorphisms constitute the desired natural isomorphism
\[
\{\gamma_{\calA}(v_X) \circ (\gamma_{\calA}(w_X))^{-1} : \mathsf{\Sigma} \rho_g X \to \rho_g \mathsf{\Sigma} X \}_{X \in \calA} :
\mathsf{\Sigma} \circ \widehat{\rho}_g \circ \gamma_{\calA} \to \widehat{\rho}_g \circ \mathsf{\Sigma} \circ \gamma_{\calA}.
\]

To this end, we have to show that any morphism $f : X \to Y$ in $\calA$ induces a commutative diagram in $\Ho(\calM)$:
\[
\xymatrix@C=3CM{
   \mathsf{\Sigma} \rho_g X
\ar[d]_-{\mathsf{\Sigma} \circ \widehat{\rho}_g \circ \gamma_{\calA}(f)}
\ar[r]^-{\gamma_{\calA}(v_X) \circ (\gamma_{\calA}(w_X))^{-1}}
&  \rho_g \mathsf{\Sigma} X
\ar[d]^-{\widehat{\rho}_g \circ \mathsf{\Sigma} \circ \gamma_{\calA}(f)}           \\
   \mathsf{\Sigma} \rho_g Y
\ar[r]^-{\gamma_{\calA}(v_Y) \circ (\gamma_{\calA}(w_Y))^{-1}}
&  \rho_g \mathsf{\Sigma} Y                                                                  }
\]
However, we have
\begin{align*}
     \mathsf{\Sigma} \circ \widehat{\rho}_g \circ \gamma_{\calA} (f)
& =  \mathsf{\Sigma} \circ \gamma_{\calA} (\rho_g f)
  =  \gamma_{\calA} \mathsf{\Sigma} (\rho_g f)
  =  \gamma_{\calA} (\mathsf{\Sigma} \rho_g f), \ \text{and}    \\
     \widehat{\rho}_g \circ \mathsf{\Sigma} \circ \gamma_{\calA} (f)
& =  \widehat{\rho}_g \circ \gamma_{\calA} \mathsf{\Sigma} (f)
  =  \widehat{\rho}_g \circ \gamma_{\calA} (\mathsf{\Sigma} f)
  =  \gamma_{\calA} (\rho_g \mathsf{\Sigma} f).
\end{align*}
Thus it suffices to show that
$$\gamma_{\calA} (\rho_g \mathsf{\Sigma} f) \circ (\gamma_{\calA}(v_X) \circ (\gamma_{\calA}(w_X))^{-1}) = (\gamma_{\calA}(v_Y) \circ (\gamma_{\calA}(w_Y))^{-1}) \circ \gamma_{\calA} (\mathsf{\Sigma} \rho_g f).$$

Now, because $(\calC, \calW \cap \calF)$ is a complete cotorsion pair in $\calA$, we have commutative diagrams:
\[
\xymatrix{
   \rho_g X
\ar[d]_-{\rho_g f}
\ar@{>->}[r]
&  W_{\rho_g X}
\ar[d]
\ar@{->>}[r]
&  \mathsf{\Sigma} \rho_g X
\ar[d]_-{\mathsf{\Sigma} \rho_g f}                      \\
   \rho_g Y
\ar@{>->}[r]
&  W_{\rho_g Y}
\ar@{->>}[r]
&  \mathsf{\Sigma} \rho_g Y                                                   }
\quad \raisebox{-3.8ex}{\text{and}} \quad
\xymatrix{
   \rho_g X
\ar[d]^-{\rho_g f}
\ar@{>->}[r]
&  \rho_g W_X
\ar[d]
\ar@{->>}[r]
&  \rho_g \mathsf{\Sigma} X
\ar[d]^-{\rho_g \mathsf{\Sigma} f}                      \\
   \rho_g Y
\ar@{>->}[r]
&  \rho_g W_Y
\ar@{->>}[r]
&  \rho_g \mathsf{\Sigma} Y                                                 }
\]
These constructions are unique up to $\omega$-homotopy, and $\gamma$ identifies homotopic morphisms.
Therefore, by the naturality of Lemma \ref{The Connecting Mor}, the square
\[
\xymatrix@C=3CM{
   \mathsf{\Sigma} \rho_g X
\ar[d]_-{\gamma_{\calA} (\mathsf{\Sigma} \rho_g f)}
\ar[r]^-{\gamma_{\calA}(v_X) \circ (\gamma_{\calA}(w_X))^{-1}}
&  \rho_g \mathsf{\Sigma} X
\ar[d]^-{\gamma_{\calA} (\rho_g \mathsf{\Sigma} f)}           \\
   \mathsf{\Sigma} \rho_g Y
\ar[r]^-{\gamma_{\calA}(v_Y) \circ (\gamma_{\calA}(w_Y))^{-1}}
&  \rho_g \mathsf{\Sigma} Y                                                                  }
\]
in $\Ho(\calM)$ commutes.
This completes the proof of naturality.

Next, the functor $\widehat{\rho}_g$ is additive because $\rho_g$ is.
Since $\rho_g$ is exact and preserves pushouts, $\widehat{\rho}_g$ preserves standard exact triangles in $\Ho(\calM)$, and hence it preserves all exact triangles.
Thus, $(\widehat{\rho}_g, \varphi_g)$ is a triangle functor.

Finally, by the same method used to construct the natural isomorphisms above, together with Lemmas \ref{natural iso} and \ref{The Connecting Mor}, one verifies that the family $\{\varphi_g\}_{g \in G}$ satisfies the associativity and unit conditions required in \ref{G fun}.
Consequently, $(\mathsf{\Sigma}, \varphi)$ is a $G$-functor.
\qedhere
\end{bfhpg}

\begin{rmk} \label{gamma is G fun}
The commutative diagram (\ref{gamma diagram}) also implies that the localization functor
\[
(\gamma_{\calA}, \Id) : (\calA, \rho, \theta) \to (\Ho(\calM), \widehat{\rho}, \widehat{\theta})
\]
is a $G$-functor.
Hence we obtain the equivariant functor $\gamma_{\calA}^G : \calA^G \to \Ho(\calM)^G$ between the corresponding equivariant categories, defined as follows:
\begin{rqm}
\item[$\bullet$]
for each object $(X, \phi)$ in $\calA^G$, $\gamma_{\calA}^G(X, \phi) = (X, \{\gamma_{\calA}(\phi_g)\})$;

\item[$\bullet$]
for each morphism $f : (X, \phi) \to (Y, \psi)$ in $\calA^G$, $\gamma_{\calA}^G(f) = \gamma_{\calA}(f)$.
\end{rqm}
\end{rmk}

\begin{bfhpg}[Proof of Lemma \ref{admiss for St(C cap F)}]\label{pf of lem 3.7}
By Lemma \ref{(Q,W,R) G-invariant}, all three classes $\calC$, $\calW$ and $\calF$ are $G$-invariant.
Hence, so is $\calC \cap \calF$.
It follows that the $G$-action $(\calA, \rho, \theta)$ on $\calA$ restricts to a $G$-action $(\calC \cap \calF, \rho, \theta)$ on $\calC \cap \calF$.
We will induce this $G$-action on the stable category $\St_\omega(\calC \cap \calF)$, and then show that it is admissible.
The proof is divided into two steps.

\vspace{1ex}\textbf{Step 1}. Construct the $G$-action $(\widetilde{\rho}, \widetilde{\theta})$.

For each element $g \in G$, define the functor $\widetilde{\rho}_g : \St_\omega(\calC \cap \calF) \to \St_\omega(\calC \cap \calF)$ as follows:
\begin{rqm}
\item[$\bullet$]
for each object $X \in \St_\omega(\calC \cap \calF)$, set $\widetilde{\rho}_g X = \rho_g X$,
\item[$\bullet$]
for each $[f] : X \to Y$ in $\St_\omega(\calC \cap \calF)$, set
$\widetilde{\rho}_g([f]) = [\rho_g f] : \rho_g X \to \rho_g Y$.
\end{rqm}
It is routine to check that $\widetilde{\rho}_g$ is an auto-equivalence, since $\rho_g$ is as well.

Next, for each $g, h \in G$, we construct the natural isomorphism
\[
\widetilde{\theta}_{g, h} : \widetilde{\rho}_g \circ \widetilde{\rho}_h \to \widetilde{\rho}_{hg}.
\]
For each object $X \in \St_\omega(\calC \cap \calF)$, set $\widetilde{\theta}_{g, h, X} = [\theta_{g, h, X}] : \rho_g \rho_h X \to \rho_{hg} X$, which is clearly an isomorphism in $\St_\omega(\calC \cap \calF)$.
For each $[f] : X \to Y$ in $\St_\omega(\calC \cap \calF)$, we have to show that the diagram
\[
\xymatrix@C=1CM{
   \rho_g \rho_h X
\ar[d]_-{[\theta_{g,h,X}]}
\ar[r]^-{[\rho_g \rho_h f]}
&  \rho_g \rho_h Y
\ar[d]^-{[\theta_{g,h,Y}]}      \\
   \rho_{hg} X
\ar[r]^-{[\rho_{hg} f]}
&  \rho_{hg} Y                                       }
\]
is commutative in $\St_\omega(\calC \cap \calF)$.
But it is clear, since $\rho_{hg} f \circ \theta_{g, h, X} = \theta_{g, h, Y} \circ \rho_g \rho_h f$ in $\calA$.

Using the 2-cocycle condition for $(\rho, \theta)$,
one can directly verify that the 2-cocycle condition for $(\widetilde{\rho}, \widetilde{\theta})$ holds.
This implies that $(\widetilde{\rho}, \widetilde{\theta})$ is a $G$-action on $\St_\omega(\calC \cap \calF)$.

\vspace{1ex}\textbf{Step 2}. Prove that $(\widetilde{\rho}, \widetilde{\theta})$ is admissible.

For each $g \in G$, we first construct a natural isomorphism
\[
\varphi_g: \mathsf{\Sigma} \circ \widetilde{\rho}_g \to \widetilde{\rho}_g \circ \mathsf{\Sigma},
\]
where $\mathsf{\Sigma}: \St_\omega(\calC \cap \calF) \to \St_\omega(\calC \cap \calF)$ is the shift functor.

Take $X \in \St_\omega(\calC \cap \calF)$.
Note that $(\calC \cap \calF, \omega)$ is a complete cotorsion pair in the exact category $\calC \cap \calF$.
There exists a commutative diagram
\[
\xymatrix{
   \rho_g X
\ar@{=}[d]
\ar@{>->}[r]
&  W_{\rho_g X}
\ar@{-->}[d]
\ar@{->>}[r]
&  \mathsf{\Sigma} \rho_g X
\ar@{-->}[d]_-{ \varphi_{g, X}} \\
   \rho_g X
\ar@{>->}[r]
&  \rho_g W_X
\ar@{->>}[r]
&  \rho_g \mathsf{\Sigma} X                               }
\]
in $\calC \cap \calF$, where $W_{\rho_g X}, \rho_g W_X \in \omega$.
Hence, define $\varphi_g(X) = [\varphi_{g, X}]$.
The uniqueness in the dual of Lemma \ref{fundamental Lemma} implies that $[\varphi_{g, X}]$ is an isomorphism in $\St_\omega(\calC \cap \calF)$.

For each $[f] : X \to Y$ in $\St_\omega(\calC \cap \calF)$, we claim that the square
\[
\xymatrix@C=1.5CM{
   \mathsf{\Sigma} \rho_g X
\ar[d]_-{[\varphi_{g, X}]}
\ar[r]^-{[\mathsf{\Sigma} \rho_g f]}
&  \mathsf{\Sigma} \rho_g Y
\ar[d]^-{[\varphi_{g, Y}]}           \\
   \rho_g \mathsf{\Sigma} X
\ar[r]^-{[\rho_g \mathsf{\Sigma} f]}
&  \rho_g \mathsf{\Sigma} Y                                      }
\]
is commutative in $\St_\omega(\calC \cap \calF)$.
Indeed, we have the following diagram with exact rows:
\[
\xymatrix@C=0.7CM @R=0.6CM{
&     \rho_g X         \ar@{>->}[rr]   \ar@{=}'[d][dd]      \ar[dl]_{\rho_g f }
&  &  W_{\rho_g X}     \ar@{->>}[rr]   \ar'[d][dd]          \ar[dl]
&  &  \mathsf{\Sigma} \rho_g X  \ar[dd]_{\varphi_{g,X}}              \ar[dl]_{\mathsf{\Sigma} \rho_g f}             \\
      \rho_g Y         \ar@{>->}[rr]   \ar@{=}[dd]
&  &  W_{\rho_g Y}     \ar@{->>}[rr]   \ar[dd]
&  &  \mathsf{\Sigma} \rho_g Y  \ar[dd]_(0.3){\varphi_{g,Y}}                                               \\
   &  \rho_g X         \ar@{>->}'[r][rr]                    \ar[dl]_{\rho_g f }
&  &  \rho_g W_X       \ar@{->>}'[r][rr]                    \ar[dl]
&  &  \rho_g \Sigma X                                       \ar[dl]_{\rho_g \mathsf{\Sigma} f}   \\
      \rho_g Y         \ar@{>->}[rr]
&  &  \rho_g W_Y       \ar@{->>}[rr]
&  &  \rho_g \mathsf{\Sigma} Y                     }
\]
It is easy to see that the squares in the upper part and the lower part of the diagram are commutative.
The uniqueness in the dual of Lemma \ref{fundamental Lemma} implies that $[\rho_g \mathsf{\Sigma} f] \circ [\varphi_{g,X}] = [\varphi_{g,Y}] \circ [\mathsf{\Sigma} \rho_g f]$.
Thus, $\varphi_g : \mathsf{\Sigma} \circ \widetilde{\rho}_g \to \widetilde{\rho}_g \circ \mathsf{\Sigma}$ is a natural isomorphism, as desired.

Now we show that $(\widetilde{\rho}_g, \varphi_g)$ is a triangle functor.
The functor $\widetilde{\rho}_g$ is additive because $\rho_g$ is.
Moreover, $\widetilde{\rho}_g$ is exact, so it preserves pushouts.
It follows that $\widetilde{\rho}_g$ preserves standard exact triangles in $\St_\omega(\calC \cap \calF)$, and hence all exact triangles.

Finally, using the uniqueness of Lemma \ref{fundamental Lemma} and its dual, one can check routinely that the family $\{\varphi_g\}_{g \in G}$ of natural isomorphisms satisfies the associativity and unit conditions in \ref{G fun}.
Hence $(\mathsf{\Sigma}, \varphi)$ is a $G$-functor.
\qedhere
\end{bfhpg}


\section*{Acknowledgments}
\noindent
We sincerely thank Miltiadis Karakikes for his prompt and insightful comments and suggestions, which helped us improve the manuscript.

\bibliographystyle{amsplain-nodash}

\def\cprime{$'$}
  \providecommand{\arxiv}[2][AC]{\mbox{\href{http://arxiv.org/abs/#2}{\sf
  arXiv:#2 [math.#1]}}}
  \providecommand{\oldarxiv}[2][AC]{\mbox{\href{http://arxiv.org/abs/math/#2}{\sf
  arXiv:math/#2
  [math.#1]}}}\providecommand{\MR}[1]{\mbox{\href{http://www.ams.org/mathscinet-getitem?mr=#1}{#1}}}
  \renewcommand{\MR}[1]{\mbox{\href{http://www.ams.org/mathscinet-getitem?mr=#1}{#1}}}
\providecommand{\bysame}{\leavevmode\hbox to3em{\hrulefill}\thinspace}
\providecommand{\MR}{\relax\ifhmode\unskip\space\fi MR }
\providecommand{\MRhref}[2]{%
  \href{http://www.ams.org/mathscinet-getitem?mr=#1}{#2}
}
\providecommand{\href}[2]{#2}

\end{document}